\documentclass[a4paper]{amsart}

\usepackage[mathscr]{eucal}
\usepackage{amsmath,amsfonts,amssymb}
\usepackage{mathrsfs}
\usepackage[dvips]{graphicx}

\newcommand{\mathsym}[1]{{}}

\theoremstyle{plain}

\newtheorem{theorem}{Theorem}[section]

\newtheorem{proposition}{Proposition}[section]
\newtheorem{lemma}{Lemma}[section]

\newtheorem{conjecture}{Conjecture}[section]

\begin{document}
\title{A detailed note\\ on the zeros of Eisenstein series for $\Gamma_0^* (5)$ and $\Gamma_0^* (7)$}
\author{Junichi Shigezumi}
\date{December 18, 2006.}

\maketitle \vspace{-0.17in}
\begin{center}
Graduate School of Mathematics, Kyushu University\\
Hakozaki 6-10-1, Higashi-ku, Fukuoka 812-8581, Japan\\
{\it E-mail address} : j.shigezumi@math.kyushu-u.ac.jp \vspace{-0.05in}
\end{center} \quad

\begin{quote}
{\small\bfseries Abstract.}
The present paper provides the details omitted from the more concise study ``On the zeros of Eisenstein series for $\Gamma_0^* (5)$ and $\Gamma_0^* (7)$.''

We locate almost all of the zeros of the Eisenstein series associated with the Fricke groups of level $5$ and $7$ in their fundamental domains by applying and extending the method of F. K. C. Rankin and H. P. F. Swinnerton-Dyer (1970). We also use the arguments of some terms of the Eisenstein series in order to improve existing error bounds.\\ \vspace{-0.15in}

\noindent
{\small\bfseries Key Words and Phrases.}
Eisenstein series, Fricke group, locating zeros, modular forms.\\ \vspace{-0.15in}

\noindent
2000 {\it Mathematics Subject Classification}. Primary 11F11; Secondary 11F12.\vspace{0.05in}
\end{quote}

\section{Introduction}

In a previous, more concise presentation of this material in ``On the zeros of Eisenstein series for $\Gamma_0^* (5)$ and $\Gamma_0^* (7)$.'' \cite{SJ}, some parts of some of the proofs were omitted for brevity. While some of these proofs rely on elementary methods, others require more complex methods. The present study thus presents the details of these proofs, and it is hoped it will be of interest to those who read the original paper. 

F. K. C. Rankin and H. P. F. Swinnerton-Dyer considered the problem of locating the zeros of the Eisenstein series $E_k(z)$ in the standard fundamental domain $\mathbb{F}$ \cite{RSD}. They proved that all of the zeros of $E_k(z)$ in $\mathbb{F}$ lie on the unit circle. They also stated towards the end of their study that ``This method can equally well be applied to Eisenstein series associated with subgroups of the modular group.'' However, it seems unclear how widely this claim holds. 

Subsequently, T. Miezaki, H. Nozaki, and the present author considered the same problem for the Fricke group $\Gamma_0^{*}(p)$ (See \cite{K}, \cite{Q}), and proved that all of the zeros of the Eisenstein series $E_{k, p}^{*}(z)$ in a certain fundamental domain lie on a circle whose radius is equal to $1 / \sqrt{p}$, $p = 2, 3$ \cite{MNS}.

The Fricke group $\Gamma_0^{*}(p)$ is not a subgroup of $\text{SL}_2(\mathbb{Z})$, but it is commensurable with $\text{SL}_2(\mathbb{Z})$. For a fixed prime $p$, we define $\Gamma_0^{*}(p) := \Gamma_0(p) \cup \Gamma_0(p) \: W_p$, where $\Gamma_0(p)$ is a congruence subgroup of $\text{SL}_2(\mathbb{Z})$.

Let $k \geqslant 4$ be an even integer. For $z \in \mathbb{H} := \{z \in \mathbb{C} \: ; \: Im(z)>0 \}$, let
\begin{equation}
E_{k, p}^{*}(z) := \frac{1}{p^{k / 2}+1} \left(p^{k/2} E_k(p z) + E_k(z) \right) \label{def:e*}
\end{equation}
be the Eisenstein series associated with $\Gamma_0^{*}(p)$. ({\it cf.} \cite{SG})

\renewcommand{\thefootnote}{\fnsymbol{footnote}}

Henceforth, we assume that $p = 5$ or $7$. The region\footnote[2]{In \cite{SJ}, there is a mistake on the definition of $\mathbb{F}^{*}(p)$. The definition in this paper is correct. We thank Prof. Rainer Schulze-Pillot for pointing it out.}
\begin{align}
\mathbb{F}^{*}(p) &:= \left\{|z| \geqslant 1 / \sqrt{p}, \: |z + 1/2| \geqslant 1 / (2 \sqrt{p}), \: - 1 / 2 \leqslant Re(z) \leqslant 0\right\} \notag\\
 &\qquad \qquad \bigcup \left\{|z| > 1 / \sqrt{p}, \: |z -1/2| > 1 / (2 \sqrt{p}), \: 0 \leqslant Re(z) < 1 / 2 \right\} \label{fd-0sp}
\end{align}
is a fundamental domain for $\Gamma_0^{*}(p)$. ({\it cf.} \cite{SH}, \cite{SE}) Define $A_p^{*} := \mathbb{F}^{*}(p) \cap \{z \in \mathbb{C} \: ; \: |z| = 1 / \sqrt{p} \; \text{or} \; |z \pm 1/2| = 1 / (2 \sqrt{p})\}$.

In the present paper, we will apply the method of F. K. C. Rankin and H. P. F. Swinnerton-Dyer (RSD Method) to the Eisenstein series associated with $\Gamma_0^{*}(5)$ and $\Gamma_0^{*}(7)$. We have the following conjectures:

\begin{conjecture}
Let $k \geqslant 4$ be an even integer. Then all of the zeros of $E_{k, 5}^{*}(z)$ in $\mathbb{F}^{*}(5)$ lie on the arc $A_5^{*}$. \label{conj-g0s5}
\end{conjecture}

\begin{conjecture}
Let $k \geqslant 4$ be an even integer. Then all of the zeros of $E_{k, 7}^{*}(z)$ in $\mathbb{F}^{*}(7)$ lie on the arc $A_7^{*}$. \label{conj-g0s7}
\end{conjecture}

First, we prove that all but at most $2$ zeros of $E_{k, p}^{*}(z)$ in $\mathbb{F}^{*}(p)$ lie on the arc $A_p^{*}$ (See Subsection \ref{subsec-g0s5-ab2} and \ref{subsec-g0s7-ab2}). Second, if $(24 / (p + 1)) \mid k$, we prove that all of the zeros of $E_{k, p}^{*}(z)$ in $\mathbb{F}^{*}(p)$ lie on $A_p^{*}$ (See Subsection \ref{subsec-g0s5-40} and \ref{subsec-g0s7-60}).

We can then prove that if $(24 / (p + 1)) \nmid k$, all but one of the zeros of $E_{k, p}^{*}(z)$ in $\mathbb{F}^{*}(p)$ lie on $A_p^{*}$. Furthermore, let $\alpha_5 \in [0, \pi]$ (resp. $\alpha_7 \in [0, \pi]$) be the angle that satisfies $\tan\alpha_5 = 2$ (resp. $\tan\alpha_7 = 5 / \sqrt{3}$), and let $\alpha_{p, k} \in [0, \pi]$ be the angle that satisfies $\alpha_{p, k} \equiv k (\pi / 2 + \alpha_p) / 2 \pmod{\pi}$. Then, since $\alpha_p$ is an irrational multiple of $\pi$, $\alpha_{p, k}$ appear uniformly in the interval $[0, \pi]$ for all even integers $k \geqslant 4$. In Subsection \ref{subsec-g0s5-41}, we prove that all of the zeros of $E_{k, 5}^{*}(z)$ in $\mathbb{F}^{*}(5)$ are on $A_5^{*}$ if $\alpha_{5, k} < (116/180) \pi$ or $(117/180) \pi < \alpha_{5, k}$. That is, we prove about $179/180$ of Conjecture \ref{conj-g0s5}. Similarly, in Subsection \ref{subsec-g0s7-61} and \ref{subsec-g0s7-62}, we prove that all of the zeros of $E_{k, 7}^{*}(z)$ in $\mathbb{F}^{*}(7)$ are on $A_7^{*}$ if ``$\alpha_{7, k} < (127.68/180) \pi$ or $(128.68/180) \pi < \alpha_{7, k}$ for $k \equiv 2 \pmod{6}$'' or ``$\alpha_{7, k} < (108.5/180) \pi$ or $(109.5/180) \pi < \alpha_{7, k}$ for $k \equiv 4 \pmod{6}$''. Thus we can also prove about $179/180$ of Conjecture \ref{conj-g0s7}.

In \cite{RSD}, we considered a bound for the error terms $R_1$ (see bound (\ref{eq-fkt}) therein) in terms only of their absolute values. However, in the present paper, we also use the arguments of some of the terms in the series. We can then approach the exact value of the Eisenstein series.

\begin{figure}[hbtp]
\begin{center}
$\begin{subarray}{c}\text{\includegraphics[width=2in]{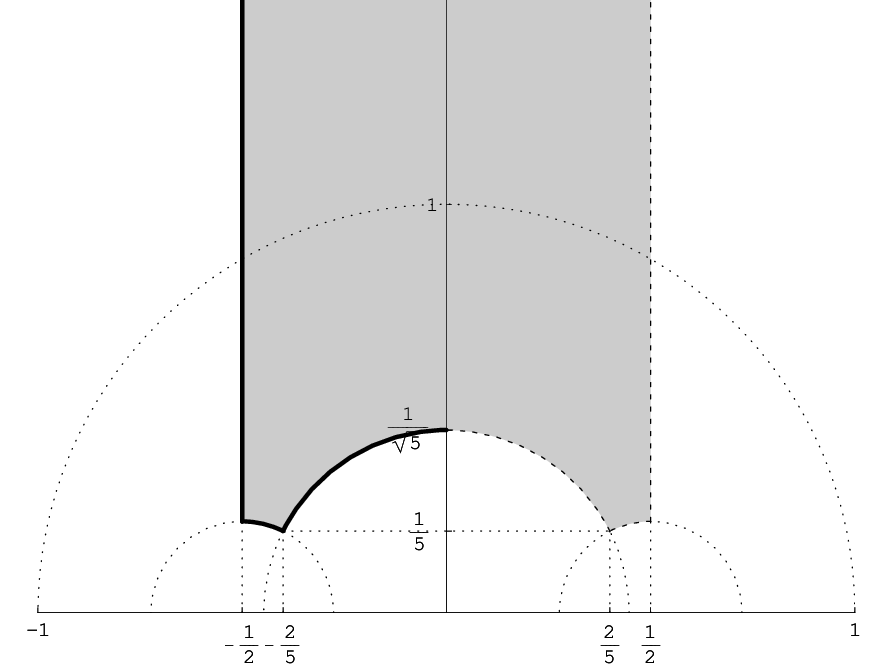}} \\ \Gamma_0^{*}(5)\end{subarray}$ \qquad
$\begin{subarray}{c}\text{\includegraphics[width=2in]{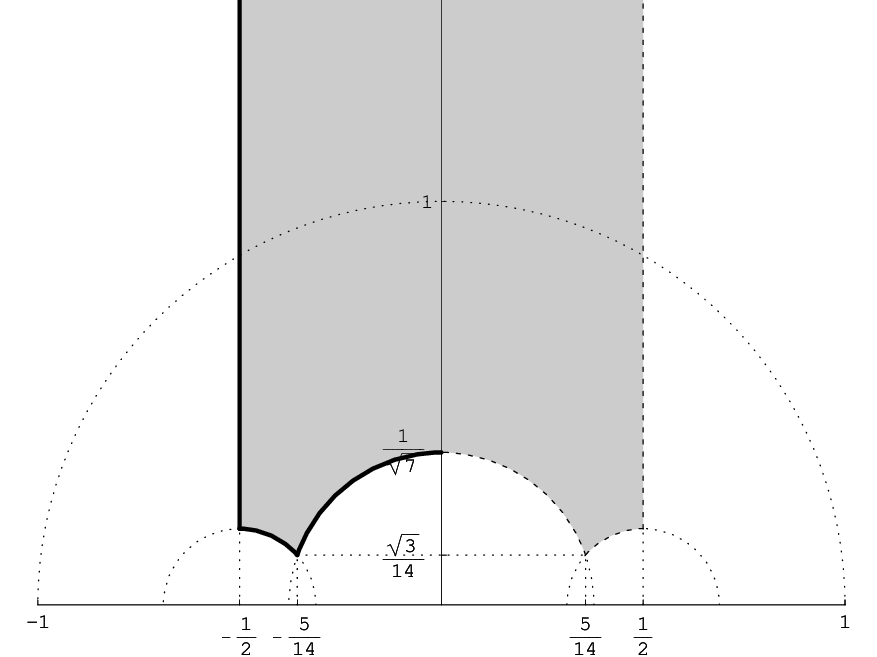}} \\ \Gamma_0^{*}(7)\end{subarray}$
\end{center}
\caption{Fundamental Domains of $\Gamma_0^{*}(5)$ and $\Gamma_0^{*}(7)$}
\label{fd-0sp}
\end{figure}

\section{General Theory}

\subsection{Preliminaries}
We refer the reader to \cite{SH} and \cite{SE}, and to \cite{SJM} for further details.

Let $\Gamma$ be a discrete subgroup of $\text{SL}_2(\mathbb{R})$. We define $P := \{ \left(\begin{smallmatrix} 1 & x \\ 0 & 1 \end{smallmatrix}\right) \: ; \: x \in \mathbb{R} \}$, and we assume that $\Gamma \cap P \setminus \left\{\pm I \right\} \ne \phi$. Let $v_p(f)$ be the order of a modular function $f$ at a point $p$.\\

\subsubsection{Fundamental Domain.}
Let $h := min\{x > 0 \: ; \: \left(\begin{smallmatrix} 1 & x \\ 0 & 1 \end{smallmatrix}\right) \in \Gamma\}$ be the {\it width} of $\Gamma$. Then, we define
\begin{equation*}
\mathbb{F}_{0,\Gamma} := \left\{z \in \mathbb{H} \: ; \: - h / 2 < Re(z) < h / 2 \: , \: |c z + d| > 1 \: \text{for all} \: \gamma = \left(\begin{smallmatrix} a & b \\ c & d \end{smallmatrix}\right) \in \Gamma \setminus P \right\}. \label{def-fd*}
\end{equation*}
Furthermore, we have a fundamental domain of $\Gamma$, $\mathbb{F}_{\Gamma}$ which satisfies
\begin{equation*}
\mathbb{F}_{0,\Gamma} \subset \mathbb{F}_{\Gamma} \subset \overline{\mathbb{F}_{0,\Gamma}}.
\end{equation*}
In this method, we have only to consider the following condition:
\begin{equation*}
\tag{$C$} |c z + d| > 1 \quad \text{for all} \: \gamma = \left(\begin{smallmatrix} a & b \\ c & d \end{smallmatrix}\right) \in \Gamma \setminus P.
\end{equation*}
Now, let $\Gamma = \Gamma_0^{*}(p)$ for $p = 5, \: 7$. Then, the condition $(C)$ is equivalent to 
\begin{equation*}
\tag{$C_p$}
|z| > 1 / \sqrt{p} \quad \text{and} \quad |z \pm 1 / 2| > 1 / 2 \sqrt{p}.
\end{equation*}

In conclusion, we have that $\mathbb{F}^{*}(p)$ is a fundamental domain for $\Gamma_0^{*}(p)$ for $p = 5, \: 7$.\\ \quad\\

\subsubsection{Eisenstein series.}
Let $\Gamma_{\infty} := \Gamma \cap P$, and let
\begin{equation}
E_{k, \infty}^{\Gamma} := e \sum_{\gamma \in \Gamma_{\infty} \setminus \Gamma} (c z + d)^{-k} \qquad z \in \mathbb{H}
\end{equation}
be the Eisenstein series associated with $\Gamma$ for the cusp $\infty$, where $e$ is a fixed number which is often chosen so that the constant term of the Fourier expansion of $E_{k, \infty}^{\Gamma}$ at $\infty$ is equal to $1$.

Now, let $\Gamma = \Gamma_0^{*}(p)$ for a prime $p$. In order to consider $\Gamma_{\infty} \setminus \Gamma_0^{*}(p)$, we need to consider $\gamma = \left( \begin{smallmatrix} a & b \\ c & d \end{smallmatrix} \right) \in \Gamma_0(p)$ and $\gamma W_p = \left( \begin{smallmatrix} b \sqrt{p} & a / \sqrt{p} \\ d \sqrt{p} & c / \sqrt{p} \end{smallmatrix} \right) \in \Gamma_0(p) W_p$. Thus, we have only to consider the pairs $(c, d)$ and $(d \sqrt{p}, c / \sqrt{p})$. Then, we have
\begin{equation*}
E_{k, \infty}^{\Gamma_0^{*}(p)} = \frac{1}{2} \sum_{\begin{subarray}{c} (c, d) = 1\\ p \mid c \end{subarray}}(c z + d)^{-k}
+ \frac{p^{k / 2}}{2} \sum_{\begin{subarray}{c} (c, d) = 1\\ p \mid d \end{subarray}}(c (p z) + d)^{-k}
\end{equation*}
as the Eisenstein series associated with $\Gamma_0^{*}(p)$ for the cusp $\infty$. Furthermore, it is easy to show that $E_{k, p}^{*} = E_{k, \infty}^{\Gamma_0^{*}(p)}$ ({cf.} Eq. (\ref{def:e*})). We can use each form as a definition.\\

\subsubsection{$\Gamma_0^{*}(5)$}

We define
\begin{align*}
A_{5, 1}^{*} &:= \{ z \: ; \: |z| = 1 / \sqrt{5}, \: \pi / 2 < Arg(z) < \pi / 2 + \alpha_5 \},\\
A_{5, 2}^{*} &:= \{ z \: ; \: |z + 1 / 2| = 1 / (2 \sqrt{5}), \: \alpha_5 < Arg(z) < \pi / 2 \}.
\end{align*}
Then, $A_5^{*} = A_{5, 1}^{*} \cup A_{5, 2}^{*} \cup \{ i / \sqrt{5}, \: \rho_{5, 1}, \: \rho_{5, 2} \}$, where $\rho_{5, 1} := - 1 / 2 + i / \left(2 \sqrt{5}\right)$ and $\rho_{5, 2} := - 2 / 5 + i / 5$.

Let $f$ be a modular form for $\Gamma_0^{*}(5)$ of weight $k$, and let $k \equiv 2 \pmod{4}$. Then, we have
\begin{align*}
f(i / \sqrt{5}) &= f(W_5 \: i / \sqrt{5}) = i^k f(i / \sqrt{5}) = - f(i / \sqrt{5}),\\
f(\rho_{5, 1}) &= f(\left(\begin{smallmatrix} 1 & - 1 \\ 0 & 1 \end{smallmatrix}\right) \left(\begin{smallmatrix} -2 & 1 \\ -5 & 2 \end{smallmatrix}\right) W_5 \: \rho_{5, 1}) = i^k f(\rho_{5, 1}) = - f(\rho_{5, 1}),\\
f(\rho_{5, 2}) &= f(\left(\begin{smallmatrix} -2 & 1 \\ -5 & 2 \end{smallmatrix}\right) \: \rho_{5, 2}) = i^k f(\rho_{5, 2}) = - f(\rho_{5, 2}).
\end{align*}
 We have further that $f(i / \sqrt{5}) = f(\rho_{5, 1}) = f(\rho_{5, 2}) = 0$, and so $v_{i / \sqrt{5}} (f) \geqslant 1$, $v_{\rho_{5, 1}} (f) \geqslant 1$, and $v_{\rho_{5, 2}} (f) \geqslant 1$.

Let $k$ be an even integer such that $k \equiv 0 \pmod{4}$. Then, we have
\begin{align*}
E_{k, 5}^{*}\left(\frac{i}{\sqrt{5}}\right) &= \frac{2 \cdot 5^{k / 2}}{5^{k / 2} + 1} E_k(\sqrt{5} i) \ne 0\\
E_{k, 5}^{*}(\rho_{5, 1}) &= \frac{2 \cdot 5^{k / 2}}{5^{k / 2} + 1} E_k\left(- \frac{1}{2} + \frac{\sqrt{5}}{2} i\right) \ne 0\\
E_{k, 5}^{*}(\rho_{5, 2}) &= \frac{1}{5^{k / 2} + 1}(5^{k / 2} + (2 + i)^k) E_k(i) \ne 0.
\end{align*}
Thus, $v_{i / \sqrt{5}} (E_{k, 5}^{*}) = v_{\rho_{5, 1}} (E_{k, 5}^{*}) = v_{\rho_{5, 2}} (E_{k, 5}^{*}) = 0$.\\

\subsubsection{$\Gamma_0^{*}(7)$}

We define
\begin{align*}
A_{7, 1}^{*} &:= \{ z \: ; \: |z| = 1 / \sqrt{7}, \: \pi / 2 < Arg(z) < \pi / 2 + \alpha_7 \},\\
A_{7, 2}^{*} &:= \{ z \: ; \: |z + 1 / 2| = 1 / (2 \sqrt{7}), \: \alpha_7 - \pi / 6 < Arg(z) < \pi / 2 \}.
\end{align*}
Then, we have $A_7^{*} = A_{7, 1}^{*} \cup A_{7, 2}^{*} \cup \{ i / \sqrt{7}, \: \rho_{7, 1}, \: \rho_{7, 2} \}$, where $\rho_{7, 1} := - 1 / 2 + i / \left(2 \sqrt{7}\right)$ and $\rho_{7, 2} := - 5 / 14 + \sqrt{3} i / 14$.

Let $f$ be a modular form for $\Gamma_0^{*}(7)$ of weight $k$, and let $k \equiv 2 \pmod{4}$. Then, we have
\begin{align*}
f(i / \sqrt{7}) &= f(W_7 \: i / \sqrt{7}) = i^k f(i / \sqrt{7}) = - f(i / \sqrt{7}),\\
f(\rho_{7, 1}) &= f(\left(\begin{smallmatrix} 1 & 1 \\ 0 & 1 \end{smallmatrix}\right) \left(\begin{smallmatrix} -3 & -1 \\ 7 & 2 \end{smallmatrix}\right) W_7 \: \rho_{7, 1}) = i^k f(\rho_{7, 1}) = - f(\rho_{7, 1}).
\end{align*}
Thus, $f(i / \sqrt{7}) = f(\rho_{7, 1}) = 0$, and so $v_{i / \sqrt{7}} (f) \geqslant 1$ and $v_{\rho_{7, 1}} (f) \geqslant 1$. On the other hand, let $k \not\equiv 0 \pmod{6}$. Then, we have
\begin{align*}
f(\rho_{7, 2}) = f(\left(\begin{smallmatrix} -3 & -1 \\ 7 & 2 \end{smallmatrix}\right) \: \rho_{7, 2}) = (e^{i 2 \pi / 3})^k f(\rho_{7, 2}).
\end{align*}
Thus $f(\rho_{7, 2}) = 0$, and so $v_{\rho_{7, 2}} (f) \geqslant 1$.

Let $k$ be an even integer such that $k \equiv 0 \pmod{4}$. Then, we have
\begin{align*}
E_{k, 7}^{*}\left(\frac{i}{\sqrt{7}}\right) &= \frac{2 \cdot 7^{k / 2}}{7^{k / 2} + 1} E_k(\sqrt{7} i) \ne 0,\\
E_{k, 7}^{*}(\rho_{7, 1}) &= \frac{2 \cdot 7^{k / 2}}{7^{k / 2} + 1} E_k\left(- \frac{1}{2} + \frac{\sqrt{7}}{2} i\right) \ne 0.
\end{align*}
Thus, $v_{i / \sqrt{7}} (E_{k, 7}^{*}) = v_{\rho_{7, 1}} (E_{k, 7}^{*}) = 0$. On the other hand, let $k$ be an even integer such that $k \equiv 0 \pmod{6}$. Then, we have
\begin{align*}
E_{k, 7}^{*}(\rho_{7, 2}) = \frac{1}{7^{k / 2} + 1}\left(7^{k / 2} + \left(\frac{5 + \sqrt{3} i}{2}\right)^k\right) E_k(\rho) \ne 0.
\end{align*}
Thus, $v_{\rho_{7, 2}} (E_{k, 7}^{*}) = 0$.\\

\subsection{Valence Formula}

In order to determine the location of zeros of $E_{k,p}^{*}(z)$ in $\mathbb{F}^{*}(p)$, we need the valence formula for $\Gamma_0^{*}(p)$.

\begin{proposition}\label{prop-vf-g0s5}
Let $f$ be a modular function of weight $k$ for $\Gamma_0^{*}(5)$, which is not identically zero. We have
\begin{equation}
v_{\infty}(f) + \frac{1}{2} v_{i / \sqrt{5}}(f) + \frac{1}{2} v_{\rho_{5,1}} (f) + \frac{1}{2} v_{\rho_{5,2}} (f) + \sum_{\begin{subarray}{c} p \in \Gamma_0^{*}(5) \setminus \mathbb{H} \\ p \ne i / \sqrt{5}, \rho_{5, 1}, \rho_{5, 2}\end{subarray}} v_p(f) = \frac{k}{4}.
\end{equation}
\end{proposition}

\begin{figure}[hbtp]
\begin{center}
\includegraphics[width=1.7in]{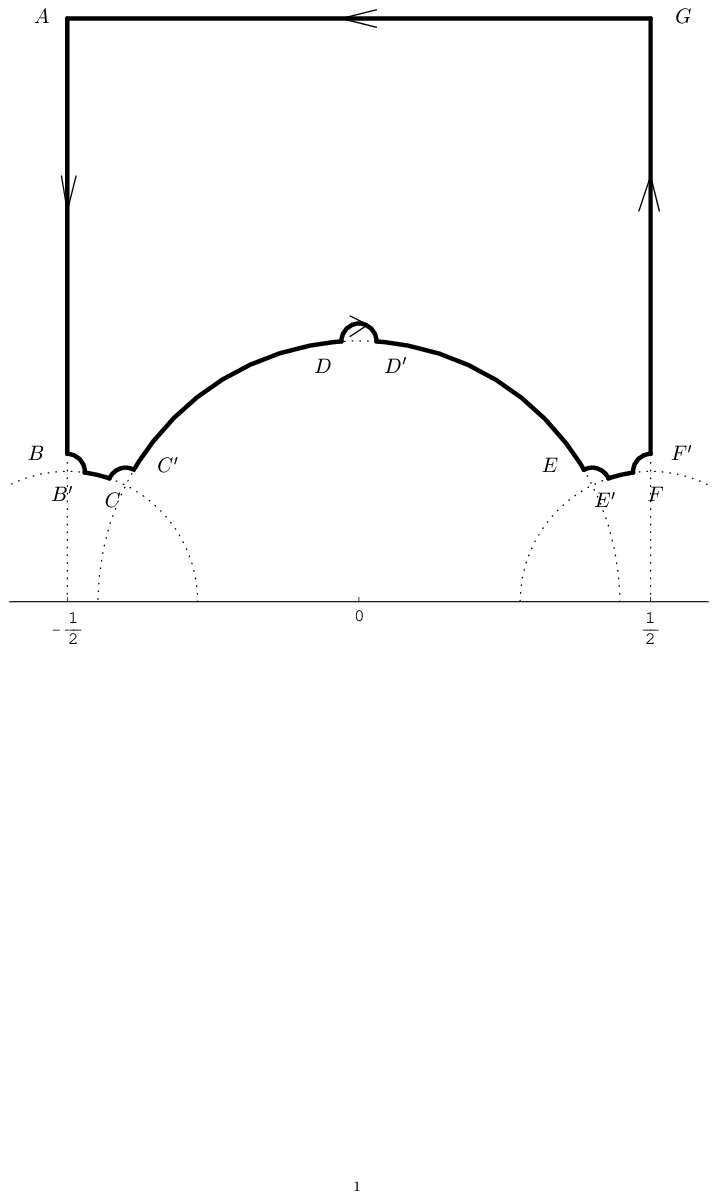}
\end{center}
\caption{}
\label{vf-0s5}
\end{figure}

\begin{proof}
Let $f$ be a nonzero modular function of weight $k$ for $\Gamma_0^{*}(5)$, and let $\mathscr{C}$ be the contour of its fundamental domain $\mathbb{F}^{*}(5)$ represented in Figure \ref{vf-0s5}, whose interior contains every zero and pole of $f$ except for $i / \sqrt{5}$, $\rho_{5, 1}$, and $\rho_{5, 2}$. By the {\it Residue theorem}, we have
\begin{equation*}
\frac{1}{2 \pi i} \int_{\mathscr{C}} \frac{df}{f} = \sum_{\begin{subarray}{c} p \in \Gamma_0^{*}(5) \setminus \mathbb{H} \\ p \ne i / \sqrt{5}, \rho_{5, 1}, \rho_{5, 2}\end{subarray}} v_p(f).
\end{equation*}
Similar to the proof of valence formula for $\text{SL}_2(\mathbb{Z})$, (See \cite{SE})
\begin{trivlist}
\item[(i)]
For the arc $G A$, we have
\begin{align*}
\frac{1}{2 \pi i} \int_{G}^{A} \frac{df}{f} = - v_{\infty}(f).
\end{align*}

\item[(ii)]
For the arcs $B B'$, $C C'$, $D D'$, $E E'$, and $F F'$, in the limit as the radii of each arc tends to $0$, we have
\begin{align*}
\frac{1}{2 \pi i} \int_{B}^{B'} \frac{df}{f} = \frac{1}{2 \pi i} \int_{F}^{F'} \frac{df}{f} &\to - \frac{1}{4} v_{\rho_{5, 1}}(f),\\
\frac{1}{2 \pi i} \int_{C}^{C'} \frac{df}{f} = \frac{1}{2 \pi i} \int_{E}^{E'} \frac{df}{f} &\to - \frac{1}{4} v_{\rho_{5, 2}}(f),\\
\frac{1}{2 \pi i} \int_{D}^{D'} \frac{df}{f} &\to - \frac{1}{2} v_{i / \sqrt{5}}(f).
\end{align*}

\item[(iii)]
For the arcs $A B$ and $F' G$, since $f(T z) = f(z)$ for $T = \left(\begin{smallmatrix} 1 & 1 \\ 0 & 1 \end{smallmatrix}\right)$,
\begin{equation*}
\frac{1}{2 \pi i} \int_{A}^{B} \frac{df}{f} + \frac{1}{2 \pi i} \int_{F'}^{G} \frac{df}{f}
 = 0.
\end{equation*}

\item[(iv)]
For the arcs $C' D$ and $D' E$, since $f(W_5 z) = (\sqrt{5} z)^k f(z)$, we have
\begin{equation*}
\frac{df(W_5 z)}{f(W_5 z)} = k \frac{dz}{z} + \frac{df(z)}{f(z)}.
\end{equation*}
In the limit as the radii of the arcs $C C'$, $D D'$, $E E'$ tend to $0$, 
\begin{equation*}
\frac{1}{2 \pi i} \int_{C'}^{D} \frac{df(z)}{f(z)} + \frac{1}{2 \pi i} \int_{D'}^{E} \frac{df(z)}{f(z)}
 = \frac{1}{2 \pi i} \int_{C'}^{D} \left( - k \frac{dz}{z} \right) \to k \frac{\theta_1}{2 \pi},
\end{equation*}
where $\tan\theta_1 = 2$.

Similarly, for the arcs $B' C$ and $E' F$, since $f\left( \left(\begin{smallmatrix} -2 & 1 \\ -5 & 2 \end{smallmatrix}\right) W_5 z\right) = (2 \sqrt{5} z + \sqrt{5})^k f(z)$, we have
\begin{equation*}
\frac{df\left( \left(\begin{smallmatrix} -2 & 1 \\ -5 & 2 \end{smallmatrix}\right) W_5 z\right)}{f\left( \left(\begin{smallmatrix} -2 & 1 \\ -5 & 2 \end{smallmatrix}\right) W_5 z\right)}
 = k \frac{dz}{z + 1 / 2} + \frac{df(z)}{f(z)}.
\end{equation*}
In the limit as the radii of the arcs $C C'$, $D D'$, $E E'$ tend to $0$, 
\begin{equation*}
\frac{1}{2 \pi i} \int_{B'}^{C} \frac{df(z)}{f(z)} + \frac{1}{2 \pi i} \int_{E'}^{F} \frac{df(z)}{f(z)}
 = \frac{1}{2 \pi i} \int_{B'}^{C} \left( - k \frac{dz}{z + 1 / 2} \right) \to k \frac{\theta_2}{2 \pi},
\end{equation*}
where $\tan\theta_1 = 1 / 2$.

Thus, since $\theta_1 + \theta_2 = \pi / 2$,
\begin{equation*}
k \frac{\theta_1}{2 \pi} + k \frac{\theta_2}{2 \pi} = \frac{k}{4}.
\end{equation*}
\end{trivlist}
\end{proof}

\begin{proposition}\label{prop-vf-g0s7}
Let $f$ be a modular function of weight $k$ for $\Gamma_0^{*}(7)$, which is not identically zero. We have
\begin{equation}
v_{\infty}(f) + \frac{1}{2} v_{i / \sqrt{7}}(f) + \frac{1}{2} v_{\rho_{7,1}} (f) + \frac{1}{3} v_{\rho_{7,2}} (f) + \sum_{\begin{subarray}{c} p \in \Gamma_0^{*}(7) \setminus \mathbb{H} \\ p \ne i / \sqrt{7}, \rho_{7, 1}, \rho_{7, 2}\end{subarray}} v_p(f) = \frac{k}{3}.
\end{equation}
\end{proposition}

The proof of above proposition is very similar to that for Proposition \ref{prop-vf-g0s5}.\\

\subsection{Some Eisenstein series of low weights}
The location of the zeros of the Eisenstein series $E_{k, 5}^{*}$ in $\mathbb{F}^{*}(5)$ for $4 \geqslant k \geqslant 10$ are follows (See Proposition \ref{prop-ek5s-lw})$:$
\begin{center}
\begin{tabular}{ccccccc}
$k$ & $v_{\infty}$ & $v_{i / \sqrt{5}}$ & $v_{\rho_{5, 1}}$ & $v_{\rho_{5, 2}}$ & $V_{5, 1}^{*}$ & $V_{5, 2}^{*}$\\
\hline
$4$ & $0$ & $0$ & $0$ & $0$ & $1$ & $0$\\
$6$ & $0$ & $1$ & $1$ & $1$ & $0$ & $0$\\
$8$ & $0$ & $0$ & $0$ & $0$ & $1$ & $1$\\
$10$ & $0$ & $1$ & $1$ & $1$ & $1$ & $0$\\
\hline
\end{tabular}
\end{center}
where $V_{5, n}^{*}$ denotes the number of simple zeros of the Eisenstein series $E_{k, 5}^{*}$ on the arc $A_{5, n}^{*}$ for $n = 1, 2$.

Similarly, the location of the zeros of the Eisenstein series $E_{k, 7}^{*}$ in $\mathbb{F}^{*}(7)$ for $k = 4, 6,$ and $12$ are follows (See Proposition \ref{prop-ek7s-lw})$:$
\begin{center}
\begin{tabular}{cccccc}
$k$ & $v_{\infty}$ & $v_{i / \sqrt{7}}$ & $v_{\rho_{7, 1}}$ & $v_{\rho_{7, 2}}$ & $V_7^{*}$\\
\hline
$4$ & $0$ & $0$ & $0$ & $1$ & $1$\\
$6$ & $0$ & $1$ & $1$ & $0$ & $1$\\
$12$ & $0$ & $0$ & $0$ & $0$ & $4$\\
\hline
\end{tabular}
\end{center}
where $V_7^{*}$ denote the number of simple zeros of the Eisenstein series $E_{k, 7}^{*}$ on $A_{7, 1}^{*} \cup A_{7, 2}^{*}$.\\

\subsection{The space of modular forms}

\subsubsection{$\Gamma_0^{*}(5)$}
Let $M_{k, 5}^{*}$ be the space of modular forms for $\Gamma_0^{*} (5)$ of weight $k$, and let $M_{k, 5}^{* 0}$ be the space of cusp forms for $\Gamma_0^{*} (5)$ of weight $k$. Consider the map $M_{k, 5}^{*} \ni f \mapsto f(\infty) \in \mathbb{C}$. The kernel of this map is $M_{k, 5}^{* 0}$, so $\dim(M_{k, 5}^{*} / M_{k, 5}^{* 0}) \leqslant 1$, and $M_{k, 5}^{*} = \mathbb{C} E_{k, 5}^{*} \oplus M_{k, 5}^{* 0}$.

Note that $\Delta_5 = \eta^4(z) \eta^4(5 z)$ is a cusp form for $\Gamma_0^{*} (5)$ of weight $4$, where $\eta(z)$ is the {\it Dedekind's $\eta$-function}. We have the following theorem:

\begin{theorem}Let $k$ be an even integer.
\def\labelenumi{(\arabic{enumi})}
\begin{enumerate}
\item
 For $k < 0$ and $k = 2$, $M_{k, 5}^{*} = 0$.
\item
 For $k = 0$ and $6$, we have $M_{k, 5}^{* 0} = 0$, and $\dim(M_{k, 5}^{*}) = 1$ with a base $E_{k, 5}^{*}$.
\item
 $M_{k, 5}^{* 0} = \Delta_5 M_{k - 4, 5}^{*}$.
\end{enumerate}
\def\labelenumi{\arabic{enumi}.}\label{th-mod_sp_5}
\end{theorem}
The proof of above theorem is very similar to that for $\text{SL}_2(\mathbb{Z})$. Furthermore, for an even integer $k \geqslant 4$, $\dim(M_{k, 5}^{*}) = (k - 2) / 4$ if $k \equiv 2 \pmod{4}$, and $\dim(M_{k, 5}^{*}) = k / 4 + 1$ if $k \equiv 0 \pmod{4}$. We have $M_{k, 5}^{*} = \mathbb{C} E_{k - 4 n, 5}^{*} (E_{4, 5}^{*})^n \oplus M_{k, 5}^{* 0}$. Then,
\begin{align*}
M_{4 n, 5}^{*} &= \mathbb{C} (E_{4, 5}^{*})^n \oplus \mathbb{C} (E_{4, 5}^{*})^{n - 1} \Delta_5 \oplus \cdots \oplus \mathbb{C} \Delta_5^n,\\
M_{4 n + 6, 5}^{*} &= E_{6, 5}^{*} ((E_{4, 5}^{*})^n \oplus \mathbb{C} (E_{4, 5}^{*})^{n - 1} \Delta_5 \oplus \cdots \oplus \mathbb{C} \Delta_5^n)
\end{align*}
Thus, for every $p \in \mathbb{H}$ and for every $f \in M_{k, 5}^{*}$, $v_p(f) \geqslant v_p(E_{k - 4 n, 5}^{*})$.

Finally, we have the following proposition:
\begin{proposition} \label{prop-bd_ord_5}
Let $k \geqslant 4$ be an even integer. For every $f \in M_{k, 5}^{*}$, we have
\begin{equation}
\begin{split}
v_{i / \sqrt{5}}(f) \geqslant s_k, \quad v_{\rho_{5, 1}}(f) \geqslant s_k, \quad v_{\rho_{5, 2}}(f) \geqslant s_k\\
(s_k=0, 1 \; \text{such that} \; 2 s_k \equiv k \pmod{4}).
\end{split}
\end{equation}
\end{proposition}\quad

\subsubsection{$\Gamma_0^{*}(7)$}
Let $M_{k, 7}^{*}$ be the space of modular forms for $\Gamma_0^{*}(7)$ of weight $k$, and let $M_{k, 7}^{* 0}$ be the space of cusp forms for $\Gamma_0^{*}(7)$ of weight $k$. Then, we have $M_{k, 7}^{*} = \mathbb{C} E_{k, 7}^{*} \oplus M_{k, 7}^{* 0}$.

Note that $\Delta_7 = \eta^6(z) \eta^6(7 z)$ is a cusp form for $\Gamma_0(7)$ of weight $6$, and $(\Delta_7)^2$ is a cusp form for $\Gamma_0^{*}(7)$ of weight $12$. We also have ${E_{2, 7}}'(z) = (7 E_2(7 z) - E_2(z)) / 6$ which is a modular form for $\Gamma_0(7)$ with $v_{\rho_{7, 2}}({E_{2, 7}}') = 2$ and $v_{p}({E_{2, 7}}') = 0$ for every $p \ne \rho_{7, 2}$. Furthermore, because we have ${E_{2, 7}}'(W_7 z) = - (\sqrt{7} z)^2 {E_{2, 7}}'(z)$, $({E_{2, 7}}')^2$ is a modular form for $\Gamma_0^{*}(7)$ of weight $4$.

We have the following theorem:
\begin{theorem}Let $k$ be an even integer.
\def\labelenumi{(\arabic{enumi})}
\begin{enumerate}
\item
 For $k < 0$ and $k = 2$, $M_{k, 7}^{*} = 0$. We have $M_{0, 7}^{*} = \mathbb{C}$.
\item
 For $k = 4, 6$, we have $M_{k, 7}^{* 0} = \mathbb{C} \Delta_{7, k}$.
\item
 For $k = 8, 10$, we have $M_{k, 7}^{* 0} = \mathbb{C} \Delta_{7, k}^0 \oplus \mathbb{C} \Delta_{7, k}^1$.
\item
 We have $M_{12, 7}^{* 0} = \mathbb{C} \Delta_{7, 12}^0 \oplus \mathbb{C} \Delta_{7, 12}^1 \oplus \mathbb{C} \Delta_{7, 12}^2 \oplus \mathbb{C} \Delta_{7, 12}^3$.
\item
 We have $M_{14, 7}^{* 0} = \mathbb{C} \Delta_{7, 14}^0 \oplus \mathbb{C} \Delta_{7, 14}^1 \oplus \mathbb{C} \Delta_{7, 14}^2$.
\item
 $M_{k, 7}^{* 0} = M_{12, 7}^{* 0} M_{k - 12, 7}^{*}$.
\end{enumerate}
\def\labelenumi{\arabic{enumi}.}
where $\Delta_{7, 4} := (5 / 16) (({E_{2, 7}}')^2 - E_{4, 7}^{*})$, $\Delta_{7, 10}^{0} := (559 / 690) ((41065 / 137592) (E_{4, 7}^{*} E_{6, 7}^{*} - E_{10, 7}^{*}) - E_{6, 7}^{*} \Delta_{7, 4})$, $\Delta_{7, 6} := \Delta_{7, 10}^0 / \Delta_{7, 4}$, $\Delta_{7, 8}^0 := (\Delta_{7, 4})^2$, $\Delta_{7, 8}^1 := E_{4, 7}^{*} \Delta_{7, 4}$, $\Delta_{7, 10}^1 := E_{6, 7}^{*} \Delta_{7, 4}$, $\Delta_{7, 12}^0 := (\Delta_7)^2$, $\Delta_{7, 12}^1 := (\Delta_{7, 4})^3$, $\Delta_{7, 12}^2 := E_{4, 7}^{*} (\Delta_{7, 4})^2$, $\Delta_{7, 12}^3 := (E_{4, 7}^{*})^2 \Delta_{7, 4}$, $\Delta_{7, 14}^0 := \Delta_{7, 4} \Delta_{7, 10}^0$, $\Delta_{7, 14}^1 := E_{6, 7}^{*} (\Delta_{7, 4})^2$, and $\Delta_{7, 14}^2 := E_{4, 7}^{*} E_{6, 7}^{*} \Delta_{7, 4}$.
\label{th-mod_sp_7}
\end{theorem}

The proof of this theorem is similar to that for Theorem \ref{th-mod_sp_5}. Regarding the orders of zeros of the basis for $M_{k, 7}^{*}$, we have the following table:

\begin{allowdisplaybreaks}
\begin{center}
\begin{tabular}{rcccccccrcccccc}
$k$ & $f$ & $v_0$ & $v_1$ & $v_2$ & $v_3$ & $V_7^{*}$ &\quad& $k$ & $f$ & $v_0$ & $v_1$ & $v_2$ & $v_3$ & $V_7^{*}$\\
\hline
$4$ & $E_{4, 7}^{*}$ & $0$ & $0$ & $0$ & $1$ & $1$ && $12$ & $E_{12, 7}^{*}$ & $0$ & $0$ & $0$ & $0$ & $4$\\
 & $({E_{2, 7}}')^2$ & $0$ & $0$ & $0$ & $4$ & $0$ && & $\Delta_{7, 12}^0$ & $4$ & $0$ & $0$ & $0$ & $0$\\
 & $\Delta_{7, 4}$ & $1$ & $0$ & $0$ & $1$ & $0$ && & $\Delta_{7, 12}^1$ & $3$ & $0$ & $0$ & $3$ & $0$\\
$6$ & $E_{6, 7}^{*}$ & $0$ & $1$ & $1$ & $0$ & $1$ && & $\Delta_{7, 12}^2$ & $2$ & $0$ & $0$ & $3$ & $1$\\
 & $\Delta_{7, 6}$ & $1$ & $1$ & $1$ & $0$ & $0$ && & $\Delta_{7, 12}^3$ & $1$ & $0$ & $0$ & $3$ & $2$\\
$8$ & $(E_{4, 7}^{*})^2$ & $0$ & $0$ & $0$ & $2$ & $2$ && $14$ & $(E_{4, 7}^{*})^2 E_{6, 7}^{*}$ & $0$ & $1$ & $1$ & $2$ & $3$\\
 & $\Delta_{7, 8}^0$ & $2$ & $0$ & $0$ & $2$ & $0$ && & $\Delta_{7, 14}^0$ & $3$ & $1$ & $1$ & $2$ & $0$\\
 & $\Delta_{7, 8}^1$ & $1$ & $0$ & $0$ & $2$ & $1$ && & $\Delta_{7, 14}^1$ & $2$ & $1$ & $1$ & $2$ & $1$\\
$10$ & $E_{4, 7}^{*} E_{6, 7}^{*}$ & $0$ & $1$ & $1$ & $1$ & $2$ && & $\Delta_{7, 14}^2$ & $1$ & $1$ & $1$ & $2$ & $2$\\
 & $\Delta_{7, 10}^0$ & $2$ & $1$ & $1$ & $1$ & $0$ &&&&&&&&\\
 & $\Delta_{7, 10}^1$ & $1$ & $1$ & $1$ & $1$ & $1$ &&&&&&&&\\
\hline
\end{tabular}
\end{center}
\end{allowdisplaybreaks}
where $V_7^{*}$ denotes the number of simple zeros of the Eisenstein series $E_{k, 7}^{*}$ on $A_{7, 1}^{*} \cup A_{7, 2}^{*}$, and let $v_0 := v_{\infty}$, $v_1 := v_{i / \sqrt{7}}$, $v_2 := v_{\rho_{7, 1}}$, and $v_3 := v_{\rho_{7, 2}}$.
\quad

Define $m_7(k) := \left\lfloor \frac{k}{3} - \frac{t}{2} \right\rfloor$, where $t=0, 2$ is chosen so that $t \equiv k \pmod{4}$. Then, we have $\dim(M_{k, 7}^{* 0}) = m_7(k)$ and $\dim(M_{k, 7}^{*}) = m_7(k) + 1$. We have $M_{k, 7}^{*} = \mathbb{C} E_{k - 12 n, 7}^{*} (E_{4, 7}^{*})^{3 n} \oplus M_{k, 7}^{* 0}$. Then
\begin{align*}
M_{k, 7}^{*}
 &= E_{k - 12 n, 7}^{*} \left\{\mathbb{C} (E_{4, 7}^{*})^{3 n} \oplus (E_{4, 7}^{*})^{3 (n - 1)} M_{12, 7}^{* 0} \oplus (E_{4, 7}^{*})^{3 (n - 2)} (M_{12, 7}^{* 0})^2 \oplus \cdots \oplus (M_{12, 7}^{* 0})^n \right\}\\
 &\qquad \oplus M_{k - 12 n, 7}^{* 0} (M_{12, 7}^{* 0})^n
\end{align*}
Thus, for every $p \in \mathbb{H}$ and for every $f \in M_{k, 7}^{*}$, $v_p(f) \geqslant v_p(E_{k - 12 n, 7}^{*})$.

Finally, we have the following proposition:
\begin{proposition}\label{prop-bd_ord_7}
Let $k \geqslant 4$ be an even integer. For every $f \in M_{k, 7}^{*}$, we have
\begin{equation}
\begin{split}
v_{i / \sqrt{7}}(f) \geqslant s_k, \quad v_{\rho_{7, 1}}(f) \geqslant s_k \quad
 &(s_k=0, 1 \; \text{such that} \; 2 s_k \equiv k \pmod{4}),\\
v_{\rho_{7, 2}}(f) \geqslant t_k \quad
 &(s_k=0, 1, 2 \; \text{such that} \; - 2 t_k \equiv k \pmod{6}).
\end{split}
\end{equation}
\end{proposition}

\newpage

\section{The method of Rankin and Swinnerton-Dyer}

\subsection{RSD Method}
Let $k \geqslant 4$ be an even integer. For $z \in \mathbb{H}$, we have
\begin{equation}
E_k(z) = \frac{1}{2} \sum_{(c,d)=1}(c z + d)^{-k}. \label{def:e}
\end{equation}
Moreover, we have $\mathbb{F} = \left\{|z| \geqslant 1, \: - 1 / 2 \leqslant Re(z) \leqslant 0\right\} \cup \left\{|z| > 1, \: 0 \leqslant Re(z) < 1 / 2 \right\}$.

At the beginning of the proof in \cite{RSD}, F. K. C. Rankin and H. P. F. Swinnerton-Dyer considered the following series:
\begin{equation}
F_k(\theta) := e^{i k \theta / 2} E_k\left(e^{i \theta}\right), \label{def:f}
\end{equation}
which is real for all $\theta \in [0, \pi]$. Considering the four terms with $c^2 + d^2 = 1$, they proved that
\begin{equation}
F_k(\theta) = 2 \cos(k \theta / 2) + R_1, \label{eq-fkt}
\end{equation}
where $R_1$ denotes the remaining terms of the series. Moreover they showed $|R_1| < 2$ for all $k \geqslant 12$. If $\cos (k \theta / 2)$ is $+1$ or $-1$, then $F_k(2 m \pi / k)$ is positive or negative, respectively. We can then show the existence of the zeros. In addition, we can prove for all of the zeros by {\it Valence Formula} and the theory on the space of modular forms for $\text{SL}_2(\mathbb{Z})$.\\

\subsection{The function: $F_{k, p, n}^{*}$}
We expect that all of the zeros of the Eisenstein series $E_{k, p}^{*}(z)$ in $\mathbb{F}^{*}(p)$ lie on the arcs $e^{i \theta} / \sqrt{p}$ and $e^{i \theta} / (2 \sqrt{p}) - 1 / 2$, which form the boundary of the fundamental domain in Figure \ref{fd-0sp}.

We define
\begin{gather}
F_{k, p, 1}^{*}(\theta) := e^{i k \theta / 2} E_{k,p}^{*}\left(e^{i \theta} / \sqrt{p}\right). \label{def:f*1}\\
F_{k, p, 2}^{*}(\theta) := e^{i k \theta / 2} E_{k,p}^{*}\left(e^{i \theta} / 2 \sqrt{p} - 1 / 2\right). \label{def:f*p2}
\end{gather}

We consider an expansion of $F_{k, p, 1}^{*}(\theta)$. Then, we have
\begin{align*}
F_{k, p, 1}^{*}(z) &= \frac{e^{i k \theta / 2}}{2} \sum_{\begin{subarray}{c} (c, d) = 1\\ p \mid c \end{subarray}}(c e^{i \theta} / \sqrt{p} + d)^{-k}
+ \frac{p^{k / 2} e^{i k \theta / 2}}{2} \sum_{\begin{subarray}{c} (c, d) = 1\\ p \mid d \end{subarray}}(c (p e^{i \theta} / \sqrt{p}) + d)^{-k}\\
 &= \frac{1}{2} \sum_{\begin{subarray}{c} (c, d) = 1\\ p \nmid d \end{subarray}}(d e^{- i \theta / 2} + \sqrt{p} c' e^{i \theta / 2})^{-k}
+ \frac{1}{2} \sum_{\begin{subarray}{c} (c, d) = 1\\ p \nmid c \end{subarray}}(c e^{i \theta / 2} + \sqrt{p} d' e^{- i \theta / 2})^{-k}.
\end{align*}

Similarly, we consider an expansion of $F_{k, p, 2}^{*}(\theta)$. When $p \mid c$, then we can write $c = c' p$ for $\exists c' \in \mathbb{Z}$, and have that $p \nmid d$. Further, when $p \mid d$, then we have $p \nmid c$ and $d = d' p$ for $\exists d' \in \mathbb{Z}$. Thus
\begin{align*}
F_{k, p, 2}^{*}(z) &= \frac{e^{i k \theta / 2}}{2} \sum_{\begin{subarray}{c} (c, d) = 1\\ p \mid c \end{subarray}}\left( c \left( \frac{e^{i \theta}}{2 \sqrt{p}} - \frac{1}{2} \right) + d \right)^{-k}
+ \frac{p^{k / 2} e^{i k \theta / 2}}{2} \sum_{\begin{subarray}{c} (c, d) = 1\\ p \mid d \end{subarray}}\left( c p \left( \frac{e^{i \theta}}{2 \sqrt{p}} - \frac{1}{2} \right) + d \right)^{-k}\\
 &= \frac{1}{2} \sum_{\begin{subarray}{c} (c, d) = 1\\ p \nmid d \end{subarray}}\left( \frac{2 d - c' p}{2} e^{- i \theta / 2} + \frac{c'}{2} \sqrt{p} e^{i \theta / 2} \right)^{-k}
+ \frac{1}{2} \sum_{\begin{subarray}{c} (c, d) = 1\\ p \nmid c \end{subarray}}\left( \frac{c}{2} e^{i \theta / 2} + \frac{2 d' - c}{2} \sqrt{p} e^{- i \theta / 2} \right)^{-k}.
\end{align*}
Now, we divide the terms into two cases, namely those terms for which $2 \mid c$ and those terms for which $2 \nmid c$. Note that the parities of $c$ and $c'$ are the same.

For the case $2 \mid c$, we can write $c' = 2 c''$ and $c = 2 c'''$ for $\exists c'', c''' \in \mathbb{Z}$. Then
\begin{equation*}
\frac{1}{2} \sum_{\begin{subarray}{c} (c, d) = 1\\ p \nmid d \end{subarray}}\left( (d - c'' p) e^{- i \theta / 2} + c'' \sqrt{p} e^{i \theta / 2} \right)^{-k}
+ \frac{1}{2} \sum_{\begin{subarray}{c} (c, d) = 1\\ p \nmid c \end{subarray}}\left( c''' e^{i \theta / 2} + (d' - c''') \sqrt{p} e^{- i \theta / 2} \right)^{-k}.
\end{equation*}
Then, we have $(d - c'' p, c'') = 1$, $p \nmid d - c'' p$, $2 \mid (d - c'') c''$, and $(c''', d' - c''') = 1$, $p \nmid c'''$, $2 \mid c'''(d' - c''')$.

For the other case $2 \nmid c$,
\begin{equation*}
\frac{2^k}{2} \sum_{\begin{subarray}{c} (c, d) = 1\\ p \nmid d \end{subarray}}\left( (2 d - c' p) e^{- i \theta / 2} + c' \sqrt{p} e^{i \theta / 2} \right)^{-k}
+ \frac{2^k}{2} \sum_{\begin{subarray}{c} (c, d) = 1\\ p \nmid c \end{subarray}}\left( c e^{i \theta / 2} + (2 d' - c) \sqrt{p} e^{- i \theta / 2} \right)^{-k}.
\end{equation*}
Then, we have $(2 d - c' p, c') = 1$, $p \nmid 2 d - c' p$, $2 \nmid (2 d - c') c'$, and $(c, 2 d' - c) = 1$, $p \nmid c$, $2 \nmid c(d' - c)$.

In conclusion, we have the following expressions:
\begin{equation}
F_{k, p, 1}^{*}(\theta)
 = \frac{1}{2} \sum_{\begin{subarray}{c} (c,d)=1\\ p \nmid c\end{subarray}}(c e^{i \theta / 2} + \sqrt{p} d e^{-i \theta / 2})^{-k}
 + \frac{1}{2} \sum_{\begin{subarray}{c} (c,d)=1\\ p \nmid c\end{subarray}}(c e^{-i \theta / 2} + \sqrt{p} d e^{i \theta / 2})^{-k}.
\end{equation}
\begin{equation}
\begin{split}
F_{k, p, 2}^{*}(\theta)
 &= \frac{1}{2} \sum_{\begin{subarray}{c} (c, d) = 1\\ p \nmid c \; 2 \mid c d \end{subarray}}\left( c e^{i \theta / 2} + d \sqrt{p} e^{- i \theta / 2} \right)^{-k}
+ \frac{1}{2} \sum_{\begin{subarray}{c} (c, d) = 1\\ p \nmid c \; 2 \mid c d \end{subarray}}\left( c e^{- i \theta / 2} + d \sqrt{p} e^{i \theta / 2} \right)^{-k}\\
&+ \frac{2^k}{2} \sum_{\begin{subarray}{c} (c, d) = 1\\ p \nmid c \; 2 \nmid c d \end{subarray}}\left( c e^{i \theta / 2} + d \sqrt{p} e^{- i \theta / 2} \right)^{-k}
+ \frac{2^k}{2} \sum_{\begin{subarray}{c} (c, d) = 1\\ p \nmid c \; 2 \nmid c d \end{subarray}}\left( c e^{- i \theta / 2} + d \sqrt{p} e^{i \theta / 2} \right)^{-k}.
\end{split}
\end{equation}

The above expressions can thus be used as definitions. Note that $(c e^{i \theta / 2} + \sqrt{p} d e^{-i \theta / 2})^{-k}$ and $(c e^{-i \theta / 2} + \sqrt{p} d e^{i \theta / 2})^{-k}$ are a complex conjugate pair for any pair $(c, \: d)$. Thus, we have the following proposition:

\begin{proposition}
$F_{k, p, 1}^{*}(\theta)$ is real, for all $\theta \in [0, \pi]$. \label{prop-f*1}
\end{proposition}
\begin{proposition}
$F_{k, p, 2}^{*}(\theta)$ is real, for all $\theta \in [0, \pi]$. \label{prop-f*2}
\end{proposition}

Now, we have
\begin{equation*}
\rho_{5, 2} = - \frac{2}{5} + \frac{i}{5} = \frac{e^{i (\pi / 2 + \alpha_5)}}{\sqrt{5}} = \frac{e^{i \alpha_5}}{2 \sqrt{5}} - \frac{1}{2},
\end{equation*}
and
\begin{align*}
F_{k, 5, 1}^{*} (\pi / 2 + \alpha_5) &= e^{i k (\pi / 2 + \alpha_5) / 2} E_{k, 5}^{*} (\rho_{5, 2}),\\
F_{k, 5, 2}^{*} (\alpha_5) &= e^{i k \alpha_5 / 2} E_{k, 5}^{*} (\rho_{5, 2}).
\end{align*}
Next, we define
\begin{equation*}
F_{k, 5}^{*} (\theta)=
\begin{cases}
F_{k, 5, 1}^{*}(\theta) & \pi / 2 \leqslant \theta \leqslant \pi / 2 + \alpha_5\\
F_{k, 5, 2}^{*}(\theta - \pi / 2) & \pi / 2 + \alpha_5 \leqslant \theta \leqslant \pi
\end{cases}.
\end{equation*}
Then, $F_{k, 5}^{*}$ is continuous in the interval $[\pi / 2, \pi]$. Note that $F_{k, 5, 1}^{*} (\pi / 2 + \alpha_5) = e^{i (\pi / 2) k / 2} F_{k, 5, 2}^{*} (\alpha_5)$.

Similarly, we have
\begin{equation*}
\rho_{7, 2} = - \frac{5}{14} + \frac{\sqrt{3} i}{14} = \frac{e^{i (\pi / 2 + \alpha_7)}}{\sqrt{7}} = \frac{e^{i (\alpha_7 - \pi / 6)}}{2 \sqrt{7}} - \frac{1}{2}.
\end{equation*}
Thus, we define
\begin{equation*}
F_{k, 7}^{*} (\theta)=
\begin{cases}
F_{k, 7, 1}^{*}(\theta) & \pi / 2 \leqslant \theta \leqslant \pi / 2 + \alpha_7\\
F_{k, 7, 2}^{*}(\theta - 2 \pi / 3) & \pi / 2 + \alpha_7 \leqslant \theta \leqslant 7 \pi / 6
\end{cases}.
\end{equation*}
Then, $F_{k, 7}^{*}$ is continuous in the interval $[\pi / 2, 7 \pi / 6]$. Note that $F_{k, 7, 1}^{*} (\pi / 2 + \alpha_7) = e^{i (2 \pi / 3) k / 2} F_{k, 7, 2}^{*} (\alpha_7 - \pi / 6)$.\\

\subsection{Application of RSD Method}
Note that we denote by $N := c^2 + d^2$. Let $v_{k}(c, d, \theta) := |c e^{i \theta / 2} + \sqrt{p} d e^{-i \theta / 2}|^{-k}$, then $v_{k}(c, d, \theta) = 1 / \left(c^2 + p d^2 + 2 \sqrt{p} c d \cos\theta \right)^{k/2}$, and $v_{k}(c, d, \theta)=v_{k}(-c, -d, \theta)$.

First, we consider the case of $N = 1$. For this case, we can write:
\begin{gather}
F_{k, p, 1}^{*}(\theta) = 2 \cos(k \theta /2) + R_{p, 1}^{*},\\
F_{k, p, 2}^{*}(\theta) = 2 \cos(k \theta /2) + R_{p, 2}^{*}
\end{gather}
where $R_{p, 1}^{*}$ and $R_{p, 2}^{*}$ denote terms for which $N > 1$ of $F_{k, p, 1}^{*}$ and $F_{k, p, 2}^{*}$, respectively.

\subsubsection{For $\Gamma_0^{*}(5)$}

For $R_{5, 1}^{*}$, we will consider the following cases: $N = 2$, $5$, $10$, $13$, $17$, and $N \geqslant 25$. Considering $- 2 / \sqrt{5} \leqslant \cos\theta \leqslant 0$, we have the following:
\begin{allowdisplaybreaks}
\begin{align*}
&\text{When $N = 2$,}&
v_k(1, 1, \theta) &\leqslant (1 / 2)^{k/2}, \quad
&v_k(1, - 1, \theta) &\leqslant (1 / 6)^{k/2}.\\
&\text{When $N = 5$,}&
v_k(1, 2, \theta) &\leqslant (1 / 13)^{k/2},
&v_k(1, - 2, \theta) &\leqslant (1 / 21)^{k/2},\\
&&v_k(2, 1, \theta) &\leqslant 1,
&v_k(2, - 1, \theta) &\leqslant (1 / 3)^{k}.\\
&\text{When $N = 10$,}&
v_k(1, 3, \theta) &\leqslant (1 / 34)^{k/2},
&v_k(1, - 3, \theta) &\leqslant (1 / 46)^{k/2},\\
&&v_k(3, 1, \theta) &\leqslant (1 / 2)^{k/2},
&v_k(3, - 1, \theta) &\leqslant (1 / 14)^{k/2}.\\
&\text{When $N = 13$,}&
v_k(2, 3, \theta) &\leqslant (1 / 5)^{k},
&v_k(2, - 3, \theta) &\leqslant (1 / 7)^{k},\\
&&v_k(3, 2, \theta) &\leqslant (1 / 5)^{k/2},
&v_k(3, - 2, \theta) &\leqslant (1 / 29)^{k/2}.\\
&\text{When $N = 17$,}&
v_k(1, 4, \theta) &\leqslant (1 / 65)^{k/2},
&v_k(1, - 4, \theta) &\leqslant (1 / 9)^{k},\\
&&v_k(4, 1, \theta) &\leqslant (1 / 5)^{k/2},
&v_k(4, - 1, \theta) &\leqslant (1 / 21)^{k/2}.\\
&\text{When $N \geqslant 25$,}&
|c e^{i \theta / 2} \pm \sqrt{5} d &e^{-i \theta / 2}|^2 \geqslant N / 6,
\end{align*}
\end{allowdisplaybreaks}
and the remaining problem for $N \geqslant 25$ concerns the number of terms with $c^2 + d^2 = N$. Because $5 \nmid c$, the number of $|c|$ is not more than $(4 / 5) N^{1/2} + 1$. Thus, the number of terms with $c^2 + d^2 = N$ is not more than $4 ((4 / 5) N^{1/2} + 1) \leqslant 4 N^{1/2}$ for $N \geqslant 25$. Then,
\begin{align*}
|R_{5, 1}^{*}|_{N \geqslant 25}
\leqslant 4 \sqrt{6} \sum_{N = 25}^{\infty} \left(\frac{1}{6} N\right)^{(1-k)/2}
\leqslant \frac{384 \sqrt{6}}{k - 3} \left(\frac{1}{2}\right)^{k}
\end{align*}
Thus,
\begin{align}
|R_{5, 1}^{*}|
 &\leqslant 2
 + 4 \left(\frac{1}{2}\right)^{k/2} + 2 \left(\frac{1}{3}\right)^{k/2}
 + \cdots + 2 \left(\frac{1}{9}\right)^{k}
 + \frac{384 \sqrt{6}}{k - 3} \left(\frac{1}{2}\right)^{k}, \label{r*51bound0}
\end{align}

On the other hand, for $R_{5, 2}^{*}$, we will consider the following cases: $N = 2$, $5$, $10$,... $29$, and $N \geqslant 34$. Considering $0 \leqslant \cos\theta \leqslant 1 / \sqrt{5}$, we have
\begin{allowdisplaybreaks}
\begin{align*}
&\text{When $N = 2$,}&
2^k \cdot v_k(1, 1, \theta) &\leqslant (2 / 3)^{k/2}, \quad
&2^k \cdot v_k(1, - 1, \theta) &\leqslant 1.\\
&\text{When $N = 5$,}&
v_k(1, 2, \theta) &\leqslant (1 / 21)^{k/2},
&v_k(1, - 2, \theta) &\leqslant (1 / 17)^{k/2},\\
&&v_k(2, 1, \theta) &\leqslant (1 / 3)^{k},
&v_k(2, - 1, \theta) &\leqslant (1 / 5)^{k/2}.\\
&\text{When $N = 10$,}&
2^k \cdot v_k(1, 3, \theta) &\leqslant (2 / 23)^{k/2},
&2^k \cdot v_k(1, - 3, \theta) &\leqslant (1 / 10)^{k/2},\\
&&2^k \cdot v_k(3, 1, \theta) &\leqslant (2 / 7)^{k/2},
&2^k \cdot v_k(3, - 1, \theta) &\leqslant (1 / 2)^{k/2}.\\
&\text{When $N = 13$,}&
v_k(2, 3, \theta) &\leqslant (1 / 7)^{k},
&v_k(2, - 3, \theta) &\leqslant (1 / 37)^{k/2},\\
&&v_k(3, 2, \theta) &\leqslant (1 / 29)^{k/2},
&v_k(3, - 2, \theta) &\leqslant (1 / 17)^{k/2}.\\
&\text{When $N = 17$,}&
v_k(1, 4, \theta) &\leqslant (1 / 9)^{k},
&v_k(1, - 4, \theta) &\leqslant (1 / 73)^{k/2},\\
&&v_k(4, 1, \theta) &\leqslant (1 / 21)^{k/2},
&v_k(4, - 1, \theta) &\leqslant (1 / 13)^{k/2}.\\
&\text{When $N = 25$,}&
v_k(3, 4, \theta) &\leqslant (1 / 89)^{k/2},
&v_k(3, - 4, \theta) &\leqslant (1 / 65)^{k/2},\\
&&v_k(4, 3, \theta) &\leqslant (1 / 61)^{k/2},
&v_k(4, - 3, \theta) &\leqslant (1 / 37)^{k/2}.\\
&\text{When $N = 26$,}&
2^k \cdot v_k(1, 5, \theta) &\leqslant (2 / 63)^{k/2},
&2^k \cdot v_k(1, - 5, \theta) &\leqslant (1 / 29)^{k/2},\\
&\text{When $N = 29$,}&
v_k(2, 5, \theta) &\leqslant (1 / 129)^{k/2},
&v_k(2, - 5, \theta) &\leqslant (1 / 109)^{k/2},\\
&\text{When $N \geqslant 34$,}&
|c e^{i \theta / 2} \pm \sqrt{5} d &e^{-i \theta / 2}|^2 \geqslant N / 2,
\end{align*}
\end{allowdisplaybreaks}
and the number of terms with $c^2 + d^2 = N$ is not more than $4 N^{1/2}$ for $N \geqslant 34$. Then,
\begin{align*}
|R_{5, 2}^{*}|_{N \geqslant 34}
\leqslant 8 \sqrt{2} \sum_{N=34}^{\infty} \left(\frac{1}{8} N\right)^{(1-k)/2}
\leqslant \frac{2112 \sqrt{33}}{k-3} \left(\frac{8}{33}\right)^{k/2}.
\end{align*}
\begin{equation}
|R_{5, 2}^{*}|
 \leqslant 2
 + 2 \left(\frac{2}{3}\right)^{k/2} + 2 \left(\frac{1}{2}\right)^{k/2}
 + \cdots + 2 \left(\frac{1}{129}\right)^{k/2}
 + \frac{2112 \sqrt{33}}{k-3} \left(\frac{8}{33}\right)^{k/2}. \label{r*52bound0}
\end{equation}

Recalling the approach of the previous subsection (RSD Method), we want to show that $|R_{5, 1}^{*}| < 2$ and $|R_{5, 2}^{*}| < 2$. However, the right-hand sides of both bounds are greater than 2. Note that the case in which $(c, d) = \pm (2, 1)$ (resp. $(c, d) = \pm (1, -1)$) yields a bound equal to $2$ for $|R_{5, 1}^{*}|$ (resp. $|R_{5, 2}^{*}|$).\\

\subsubsection{$\Gamma_0^{*}(7)$}
For $R_{7, 1}^{*}$, we will consider the following cases: $N = 2$, $5$, $10$,..., $61$, and $N \geqslant 65$. Considering $- 5 / (2 \sqrt{7}) \leqslant \cos\theta \leqslant 0$, we have
\begin{allowdisplaybreaks}
\begin{align*}
&\text{When $N = 2$,}&
v_k(1, 1, \theta) &\leqslant (1 / 3)^{k/2}, \quad
&v_k(1, - 1, \theta) &\leqslant (1 / 8)^{k/2}.\\
&\text{When $N = 5$,}&
v_k(1, 2, \theta) &\leqslant (1 / 19)^{k/2},
&v_k(1, - 2, \theta) &\leqslant (1 / 29)^{k/2},\\
&&v_k(2, 1, \theta) &\leqslant 1,
&v_k(2, - 1, \theta) &\leqslant (1 / 11)^{k/2}.\\
&\text{When $N = 10$,}&
v_k(1, 3, \theta) &\leqslant (1 / 7)^{k},
&v_k(1, - 3, \theta) &\leqslant (1 / 8)^{k},\\
&&v_k(3, 1, \theta) &\leqslant 1,
&v_k(3, - 1, \theta) &\leqslant (1 / 4)^{k}.\\
&\text{When $N = 13$,}&
v_k(2, 3, \theta) &\leqslant (1 / 39)^{k/2},
&v_k(2, - 3, \theta) &\leqslant (1 / 69)^{k/2},\\
&&v_k(3, 2, \theta) &\leqslant (1 / 7)^{k/2},
&v_k(3, - 2, \theta) &\leqslant (1 / 37)^{k/2}.\\
&\text{When $N = 17$,}&
v_k(1, 4, \theta) &\leqslant (1 / 93)^{k/2},
&v_k(1, - 4, \theta) &\leqslant (1 / 113)^{k/2},\\
&&v_k(4, 1, \theta) &\leqslant (1 / 3)^{k/2},
&v_k(4, - 1, \theta) &\leqslant (1 / 23)^{k/2}.\\
&\text{When $N = 25$,}&
v_k(3, 4, \theta) &\leqslant (1 / 61)^{k/2},
&v_k(3, - 4, \theta) &\leqslant (1 / 11)^{k},\\
&&v_k(4, 3, \theta) &\leqslant (1 / 19)^{k/2},
&v_k(4, - 3, \theta) &\leqslant (1 / 79)^{k/2}.\\
&\text{When $N = 26$,}&
v_k(1, 5, \theta) &\leqslant (1 / 151)^{k/2},
&v_k(1, - 5, \theta) &\leqslant (1 / 176)^{k/2},\\
&&v_k(5, 1, \theta) &\leqslant (1 / 7)^{k/2},
&v_k(5, - 1, \theta) &\leqslant (1 / 32)^{k/2}.\\
&\text{When $N = 29$,}&
v_k(2, 5, \theta) &\leqslant (1 / 129)^{k/2},
&v_k(2, - 5, \theta) &\leqslant (1 / 179)^{k/2},\\
&&v_k(5, 2, \theta) &\leqslant (1 / 3)^{k/2},
&v_k(5, - 2, \theta) &\leqslant (1 / 53)^{k/2}.\\
&\text{When $N = 34$,}&
v_k(3, 5, \theta) &\leqslant (1 / 109)^{k/2},
&v_k(3, - 5, \theta) &\leqslant (1 / 184)^{k/2},\\
&&v_k(5, 3, \theta) &\leqslant (1 / 13)^{k/2},
&v_k(5, - 3, \theta) &\leqslant (1 / 88)^{k/2}.\\
&\text{When $N = 37$,}&
v_k(1, 6, \theta) &\leqslant (1 / 223)^{k/2},
&v_k(1, - 6, \theta) &\leqslant (1 / 253)^{k/2},\\
&&v_k(6, 1, \theta) &\leqslant (1 / 13)^{k/2},
&v_k(6, - 1, \theta) &\leqslant (1 / 43)^{k/2}.\\
&\text{When $N = 41$,}&
v_k(4, 5, \theta) &\leqslant (1 / 91)^{k/2},
&v_k(4, - 5, \theta) &\leqslant (1 / 191)^{k/2},\\
&&v_k(5, 4, \theta) &\leqslant (1 / 37)^{k/2},
&v_k(5, - 4, \theta) &\leqslant (1 / 137)^{k/2}.\\
&\text{When $N = 50$,}&
v_k(1, 7, \theta) &\leqslant (1 / 309)^{k/2},
&v_k(1, - 7, \theta) &\leqslant (1 / 344)^{k/2},\\
&\text{When $N = 53$,}&
v_k(2, 7, \theta) &\leqslant (1 / 277)^{k/2},
&v_k(2, - 7, \theta) &\leqslant (1 / 347)^{k/2},\\
&\text{When $N = 58$,}&
v_k(3, 7, \theta) &\leqslant (1 / 247)^{k/2},
&v_k(3, - 7, \theta) &\leqslant (1 / 352)^{k/2},\\
&\text{When $N = 61$,}&
v_k(5, 6, \theta) &\leqslant (1 / 127)^{k/2},
&v_k(5, - 6, \theta) &\leqslant (1 / 277)^{k/2},\\
&&v_k(6, 5, \theta) &\leqslant (1 / 61)^{k/2},
&v_k(6, - 5, \theta) &\leqslant (1 / 211)^{k/2}.\\
&\text{When $N \geqslant 65$,}&
|c e^{i \theta / 2} \pm \sqrt{7} d &e^{-i \theta / 2}|^2 \geqslant N / 11,
\end{align*}
\end{allowdisplaybreaks}
and the remaining problem for $N \geqslant 65$ concerns the number of terms with $c^2 + d^2 = N$. Because $7 \nmid c$, the number of $|c|$ is not more than $(6 / 7) N^{1/2} + 1$. Thus, the number of terms with $c^2 + d^2 = N$ is not more than $4 ((6 / 7) N^{1/2} + 1) \leqslant (55 / 14) N^{1/2}$ for $N \geqslant 65$. Then,
\begin{align*}
|R_{7, 1}^{*}|_{N \geqslant 65}
\leqslant \frac{55 \sqrt{11}}{14} \sum_{N=65}^{\infty} \left(\frac{1}{11} N\right)^{(1-k)/2}
\leqslant \frac{28160}{7 (k - 3)} \left(\frac{11}{64}\right)^{k/2}.
\end{align*}
Thus,
\begin{equation}
|R_{7, 1}^{*}|
 \leqslant 4
 + 6 \left(\frac{1}{3}\right)^{k/2} + 4 \left(\frac{1}{7}\right)^{k/2}
 + \cdots + 2 \left(\frac{1}{352}\right)^{k/2}
 + \frac{28160}{7(k-3)} \left(\frac{11}{64}\right)^{k/2}, \label{r*71bound0}
\end{equation}

On the other hand, for $R_{7, 2}^{*}$, we will consider the following cases: $N = 2$, $5$, $10$,..., $89$, and $N \geqslant 97$. Considering $0 \leqslant \cos\theta \leqslant 2 / \sqrt{7}$, we have
\begin{allowdisplaybreaks}
\begin{align*}
&\text{When $N = 2$,}&
2^k \cdot v_k(1, 1, \theta) &\leqslant (1 / 2)^{k/2}, \quad
&2^k \cdot v_k(1, - 1, \theta) &\leqslant 1.\\
&\text{When $N = 5$,}&
v_k(1, 2, \theta) &\leqslant (1 / 29)^{k/2},
&v_k(1, - 2, \theta) &\leqslant (1 / 21)^{k/2},\\
&&v_k(2, 1, \theta) &\leqslant (1 / 11)^{k/2},
&v_k(2, - 1, \theta) &\leqslant (1 / 3)^{k/2}.\\
&\text{When $N = 10$,}&
2^k \cdot v_k(1, 3, \theta) &\leqslant (1 / 4)^{k},
&2^k \cdot v_k(1, - 3, \theta) &\leqslant (1 / 13)^{k/2},\\
&&2^k \cdot v_k(3, 1, \theta) &\leqslant (1 / 2)^{k},
&2^k \cdot v_k(3, - 1, \theta) &\leqslant 1.\\
&\text{When $N = 13$,}&
v_k(2, 3, \theta) &\leqslant (1 / 69)^{k/2},
&v_k(2, - 3, \theta) &\leqslant (1 / 45)^{k/2},\\
&&v_k(3, 2, \theta) &\leqslant (1 / 37)^{k/2},
&v_k(3, - 2, \theta) &\leqslant (1 / 14)^{k/2}.\\
&\text{When $N = 17$,}&
v_k(1, 4, \theta) &\leqslant (1 / 113)^{k/2},
&v_k(1, - 4, \theta) &\leqslant (1 / 97)^{k/2},\\
&&v_k(4, 1, \theta) &\leqslant (1 / 23)^{k/2},
&v_k(4, - 1, \theta) &\leqslant (1 / 7)^{k/2}.\\
&\text{When $N = 25$,}&
v_k(3, 4, \theta) &\leqslant (1 / 11)^{k},
&v_k(3, - 4, \theta) &\leqslant (1 / 73)^{k/2},\\
&&v_k(4, 3, \theta) &\leqslant (1 / 79)^{k/2},
&v_k(4, - 3, \theta) &\leqslant (1 / 6)^{k/2}.\\
&\text{When $N = 26$,}&
2^k \cdot v_k(1, 5, \theta) &\leqslant (1 / 44)^{k/2},
&2^k \cdot v_k(1, - 5, \theta) &\leqslant (1 / 39)^{k/2},\\
&&2^k \cdot v_k(5, 1, \theta) &\leqslant (1 / 8)^{k/2},
&2^k \cdot v_k(5, - 1, \theta) &\leqslant (1 / 3)^{k/2}.\\
&\text{When $N = 29$,}&
v_k(2, 5, \theta) &\leqslant (1 / 179)^{k/2},
&v_k(2, - 5, \theta) &\leqslant (1 / 139)^{k/2},\\
&&v_k(5, 2, \theta) &\leqslant (1 / 53)^{k/2},
&v_k(5, - 2, \theta) &\leqslant (1 / 13)^{k/2}.\\
&\text{When $N = 34$,}&
2^k \cdot v_k(3, 5, \theta) &\leqslant (1 / 46)^{k/2},
&2^k \cdot v_k(3, - 5, \theta) &\leqslant (1 / 31)^{k/2},\\
&&2^k \cdot v_k(5, 3, \theta) &\leqslant (1 / 22)^{k/2},
&2^k \cdot v_k(5, - 3, \theta) &\leqslant (1 / 7)^{k/2}.\\
&\text{When $N = 37$,}&
v_k(1, 6, \theta) &\leqslant (1 / 252)^{k/2},
&v_k(1, - 6, \theta) &\leqslant (1 / 228)^{k/2},\\
&&v_k(6, 1, \theta) &\leqslant (1 / 43)^{k/2},
&v_k(6, - 1, \theta) &\leqslant (1 / 19)^{k/2}.\\
&\text{When $N = 41$,}&
v_k(4, 5, \theta) &\leqslant (1 / 191)^{k/2},
&v_k(4, - 5, \theta) &\leqslant (1 / 111)^{k/2},\\
&&v_k(5, 4, \theta) &\leqslant (1 / 137)^{k/2},
&v_k(5, - 4, \theta) &\leqslant (1 / 57)^{k/2}.\\
&\text{When $N = 50$,}&
2^k \cdot v_k(1, 7, \theta) &\leqslant (1 / 86)^{k/2},
&2^k \cdot v_k(1, - 7, \theta) &\leqslant (1 / 79)^{k/2},\\
&\text{When $N = 53$,}&
v_k(2, 7, \theta) &\leqslant (1 / 347)^{k/2},
&v_k(2, - 7, \theta) &\leqslant (1 / 291)^{k/2},\\
&\text{When $N = 58$,}&
2^k \cdot v_k(3, 7, \theta) &\leqslant (1 / 88)^{k/2},
&2^k \cdot v_k(3, - 7, \theta) &\leqslant (1 / 67)^{k/2},\\
&\text{When $N = 61$,}&
v_k(5, 6, \theta) &\leqslant (1 / 277)^{k/2},
&v_k(5, - 6, \theta) &\leqslant (1 / 157)^{k/2},\\
&&v_k(6, 5, \theta) &\leqslant (1 / 211)^{k/2},
&v_k(6, - 5, \theta) &\leqslant (1 / 91)^{k/2}.\\
&\text{When $N = 65$,}&
v_k(1, 8, \theta) &\leqslant (1 / 449)^{k/2},
&v_k(1, - 8, \theta) &\leqslant (1 / 417)^{k/2},\\
&&v_k(8, 1, \theta) &\leqslant (1 / 71)^{k/2},
&v_k(8, - 1, \theta) &\leqslant (1 / 39)^{k/2}.\\
&&v_k(4, 7, \theta) &\leqslant (1 / 359)^{k/2},
&v_k(4, - 7, \theta) &\leqslant (1 / 247)^{k/2}.\\
&\text{When $N = 73$,}&
v_k(3, 8, \theta) &\leqslant (1 / 457)^{k/2},
&v_k(3, - 8, \theta) &\leqslant (1 / 361)^{k/2},\\
&&v_k(8, 3, \theta) &\leqslant (1 / 127)^{k/2},
&v_k(8, - 3, \theta) &\leqslant (1 / 31)^{k/2}.\\
&\text{When $N = 74$,}&
2^k \cdot v_k(5, 7, \theta) &\leqslant (1 / 92)^{k/2},
&2^k \cdot v_k(5, - 7, \theta) &\leqslant (1 / 57)^{k/2},\\
&\text{When $N = 82$,}&
2^k \cdot v_k(1, 9, \theta) &\leqslant (1 / 142)^{k/2},
&2^k \cdot v_k(1, - 9, \theta) &\leqslant (1 / 133)^{k/2},\\
&&2^k \cdot v_k(9, 1, \theta) &\leqslant (1 / 22)^{k/2},
&2^k \cdot v_k(9, - 1, \theta) &\leqslant (1 / 13)^{k/2}.\\
&\text{When $N = 85$,}&
v_k(2, 9, \theta) &\leqslant (1 / 571)^{k/2},
&v_k(2, - 9, \theta) &\leqslant (1 / 499)^{k/2},\\
&&v_k(9, 2, \theta) &\leqslant (1 / 109)^{k/2},
&v_k(9, - 2, \theta) &\leqslant (1 / 37)^{k/2}.\\
&&v_k(6, 7, \theta) &\leqslant (1 / 379)^{k/2},
&v_k(6, - 7, \theta) &\leqslant (1 / 211)^{k/2}.\\
&\text{When $N = 89$,}&
v_k(5, 8, \theta) &\leqslant (1 / 473)^{k/2},
&v_k(5, - 8, \theta) &\leqslant (1 / 313)^{k/2},\\
&&v_k(8, 5, \theta) &\leqslant (1 / 239)^{k/2},
&v_k(8, - 5, \theta) &\leqslant (1 / 79)^{k/2}.\\
&\text{When $N \geqslant 97$,}&
|c e^{i \theta / 2} \pm \sqrt{7} d &e^{-i \theta / 2}|^2 \geqslant N / 3,
\end{align*}
\end{allowdisplaybreaks}
and the number of terms with $c^2 + d^2 = N$ is not more than $(244 / 63) N^{1/2}$ for $N \geqslant 97$. Then,
\begin{align*}
|R_{7, 2}^{*}|_{N \geqslant 97}
\leqslant \frac{488 \sqrt{3}}{63} \sum_{N=97}^{\infty} \left(\frac{1}{12} N\right)^{(1-k)/2}
\leqslant \frac{62464 \sqrt{6}}{21 (k - 3)} \left(\frac{1}{8}\right)^{k/2}.
\end{align*}
Thus,
\begin{equation}
|R_{7, 2}^{*}|
 \leqslant 4
 + 2 \left(\frac{1}{2}\right)^{k/2} + 2 \left(\frac{1}{3}\right)^{k/2}
 + \cdots + 2 \left(\frac{1}{571}\right)^{k/2}
 + \frac{62464 \sqrt{6}}{21(k - 3)} \left(\frac{1}{8}\right)^{k/2}. \label{r*72bound0}
\end{equation}

We want to show that $|R_{7, 1}^{*}| < 2$ and $|R_{7, 2}^{*}| < 2$. But the right-hand sides of both bounds are much greater than 2. Note that the cases $(c, d) = \pm (2, 1)$ and $\pm (3, 1)$ (resp. $(c, d) = \pm (1, -1)$ and $\pm (3, -1)$) yield a bound equal to $4$ for $|R_{7, 1}^{*}|$ (resp. $|R_{7, 2}^{*}|$).\\

\subsection{Arguments of some terms}

\subsubsection{$\Gamma_0^{*}(5)$}
In the previous subsection, the point was the cases of $(c, d) = \pm (2, 1)$ and $(c, d) = \pm (1, -1)$ for the bounds $|R_{5, 1}^{*}|$ and $|R_{5, 2}^{*}|$, respectively.

For the case in which $(c, d) = \pm (2, 1)$, we have
\begin{equation*}
2 e^{i \theta_1 / 2} + \sqrt{5} e^{- i \theta_1 / 2} = (\sqrt{5} + 2) \cos \theta_1 / 2 - i (\sqrt{5} - 2) \sin \theta_1 / 2.
\end{equation*}
Let ${\theta_1}' := 2 Arg \left\{ 2 e^{i \theta_1 / 2} + \sqrt{5} e^{- i \theta_1 / 2} \right\}$, then we have
\begin{equation*}
\tan {\theta_1}' / 2 = \frac{- (\sqrt{5} - 2) \sin \theta_1 / 2}{(\sqrt{5} + 2) \cos \theta_1 / 2}
 = - (9 - 4 \sqrt{5}) \tan \theta_1 / 2.
\end{equation*}
Furthermore, since we have $\tan (\pi / 2 \pm \alpha_5) / 2 = \sqrt{5} \pm 2$, for a positive number $d_1$,
\begin{align*}
&- (9 - 4 \sqrt{5}) \tan (\pi / 2 + \alpha_5 - (t \pi / k)) / 2 - \tan (- \pi + \pi / 2 + \alpha_5 + d_1 (t \pi / k)) / 2\\
= &- \frac{(9 - 4 \sqrt{5}) (\sqrt{5} + 2) - (9 - 4 \sqrt{5}) \tan (t \pi / k) / 2}{1 + (\sqrt{5} + 2) \tan (t \pi / k) / 2}
 + \frac{(\sqrt{5} - 2) - \tan d_1 (t \pi / k) / 2}{1 + (\sqrt{5} - 2) \tan d_1 (t \pi / k) / 2}\\
= &\frac{(10 - 4 \sqrt{5}) \left\{ \tan (t \pi / k) / 2 - \tan d_1(t \pi / k) / 2 - 4 \tan (t \pi / k) / 2 \tan d_1 (t \pi / k) / 2 \right\}}{\left\{ 1 + (\sqrt{5} + 2) \tan (t \pi / k) / 2 \right\} \left\{ 1 + (\sqrt{5} - 2) \tan d_1 (t \pi / k) / 2 \right\}}.
\end{align*}
Let $D := \tan (t \pi / k) / 2 - \tan d_1(t \pi / k) / 2 - 4 \tan (t \pi / k) / 2 \tan d_1 (t \pi / k) / 2$, which is a factor of the numerator of the above fraction. When $t$ is small enough or $k$ is large enough, then the denominator of the above fraction is positive, thus the sign of the fraction is equal to that of $D$.

When $d_1 > 1$, then we have $D < 0$. On the other hand, when $d_1 < 1$ and $(t \pi / k) / 2 < \pi / 2$ (note that the latter condition is necessary condition for the condition``$t$ is small enough or $k$ is large enough''), then we have $\tan d_1 (t \pi / k) / 2 < d_1 \tan (t \pi / k) / 2$, thus we have $D > \tan (t \pi / k) / 2 \{ 1 - d_1 - 4 d_1 \tan (t \pi / k) / 2 \}$. Therefore, when $d_1 < 1 / (1 + 4 \tan (t / 2) (\pi / k))$, we have $D > 0$.

In conclusion, we have
\begin{align*}
\tan \left( - \pi + \pi / 2 + \alpha_5 + d_1 (t \pi / k) \right) / 2 &< - (9 - 4 \sqrt{5}) \tan \left(\pi / 2 + \alpha_5 - (t \pi / k) \right) / 2\\
 &< \tan \left( - \pi + \pi / 2 + \alpha_5 + (t \pi / k) \right) / 2,
\end{align*}
where $d_1 < 1 / (1 + 4 \tan (t / 2) (\pi / k))$. Thus,
\begin{align*}
&\theta_1 = \pi / 2 + \alpha_5 - (t \pi / k)\\
 &\qquad \Rightarrow - \pi + \pi / 2 + \alpha_5 + d_1 (t \pi / k) < {\theta_1}' < - \pi + \pi / 2 + \alpha_5 + (t \pi / k),\\
&k \theta_1 / 2 = k (\pi / 2 + \alpha_5) / 2 - (t / 2) \pi \\
 &\qquad \Rightarrow - (k / 2) \pi + k (\pi / 2 + \alpha_5) / 2 + d_1 (t / 2) \pi < k {\theta_1}' / 2 < - (k / 2) \pi + k (\pi / 2 + \alpha_5) / 2 + (t / 2) \pi.
\end{align*}

Similarly, for the case $(c, d) = \pm (1, -1)$, let ${\theta_2}' := 2 Arg \left\{ - e^{i \theta_2 / 2} + \sqrt{5} e^{- i \theta_2 / 2} \right\}$, then we have $\tan {\theta_2}' / 2 = - ((3 + \sqrt{5}) / 2) \tan \theta_2 / 2$. For a positive number $d_2$, we have
\begin{align*}
&- ((3 + \sqrt{5}) / 2) \tan (\alpha_5 + (t \pi / k)) / 2 - \tan (- \pi + \alpha_5 - d_2 (t \pi / k)) / 2\\
= &\frac{(5 + \sqrt{5}) \left\{ \tan d_2 (t \pi / k) / 2 - \tan (t \pi / k) / 2 + \tan (t \pi / k) / 2 \tan d_2 (t \pi / k) / 2 \right\}}{2 \left\{ 1 - ((\sqrt{5} - 1) / 2) \tan (t \pi / k) / 2 \right\} \left\{ 1 - ((\sqrt{5} + 1) / 2) \tan d_2 (t \pi / k) / 2 \right\}}.
\end{align*}
Thus,
\begin{align*}
\tan \left( - \pi + \alpha_5 - (t \pi / k)\right) / 2 &< - ((3 + \sqrt{5}) / 2) \tan \left( \alpha_5 + (t \pi / k) \right) / 2\\
 &< \tan \left( - \pi + \alpha_5 - d_2 (t \pi / k) \right) / 2,
\end{align*}
where $d_2 < 1 / (1 + \tan (t / 2) (\pi / k))$. Furthermore, we have
\begin{align*}
&k \theta_2 / 2 = k \alpha_5 / 2 + (t / 2) \pi \\
 &\qquad \Rightarrow - (k / 2) \pi + k \alpha_5 / 2 - (t / 2) \pi < k {\theta_2}' / 2 < - (k / 2) \pi + k \alpha_5 / 2 - d_2 (t / 2) \pi.
\end{align*}

Note that
\begin{equation*}
- (k / 2) \pi \equiv
\begin{cases}
0 & (k \equiv 0 \pmod{4})\\
\pi & (k \equiv 2 \pmod{4})
\end{cases} \quad \pmod{2 \pi},
\end{equation*}
and both $d_1$ and $d_2$ tend to $1$ in the limit as $k$ tends to $\infty$ or in the limit as $t$ tends to $0$.

Recall that $\alpha_{5, k} \equiv k (\pi / 2 + \alpha_5) / 2 \pmod{\pi}$, then we can write $k (\pi / 2 + \alpha_5) / 2 = \alpha_{5, k} + m \pi$ for some integer $m$. We define ${\alpha_{5, k}}' \equiv k {\theta_1}' / 2 - m \pi \pmod {2 \pi}$ for $\theta_1 = \pi / 2 + \alpha_5 - (t \pi / k)$.

Similarly, we define $\beta_{5, k} \equiv k \alpha_5 / 2 \pmod{\pi}$ and ${\beta_{5, k}}' \equiv k {\theta_2}' / 2 - (k \alpha_5 / 2 - \beta_{5, k} ) \pmod {2 \pi}$ for $\theta_2 = \alpha_5 + (t \pi / k)$.\\

\subsection{$\Gamma_0^{*}(7)$}

Similar to the previous subsubsection, we consider the arguments of some terms such that $(c, d) = \pm (2, 1)$ and $\pm (3, 1)$ for $|R_{7, 1}^{*}|$, and $(c, d) = \pm (1, -1)$ and $\pm (3, -1)$ for $|R_{7, 2}^{*}|$.

Let ${\theta_{1, 1}}' := 2 Arg \left\{ 2 e^{i \theta_1 / 2} + \sqrt{7} e^{- i \theta_1 / 2} \right\}$, ${\theta_{1, 2}}' := 2 Arg \left\{ 3 e^{i \theta_1 / 2} + \sqrt{7} e^{- i \theta_1 / 2} \right\}$,\\
 ${\theta_{2, 1}}' := 2 Arg \left\{ - e^{i \theta_2 / 2} + \sqrt{7} e^{- i \theta_2 / 2} \right\}$, and ${\theta_{2, 2}}' := 2 Arg \left\{ - 3 e^{i \theta_2 / 2} + \sqrt{7} e^{- i \theta_2 / 2} \right\}$.\\
Then, we have $\tan {\theta_{1, 1}}' / 2 = - ((11 - 4 \sqrt{7}) / 3) \tan \theta_1 / 2$, $\tan {\theta_{1, 2}}' / 2 = - (8 - 3 \sqrt{7}) \tan \theta_1 / 2$, $\tan {\theta_{2, 1}}' / 2 = - ((4 + \sqrt{7}) / 3) \tan \theta_1 / 2$, and $\tan {\theta_{2, 2}}' / 2 = - (8 + 3 \sqrt{7}) \tan \theta_1 / 2$.

When $\theta_1 = \pi / 2 + \alpha_7 - (t \pi / k) \pi$, we have
\begin{allowdisplaybreaks}
\begin{align*}
&- ((11 - 4 \sqrt{7}) / 3) \tan (\pi / 2 + \alpha_7 - (t \pi / k) \pi) / 2 - \tan (2 \pi / 3 + \pi / 2 + \alpha_7 + d_{1, 1} (t \pi / k)) / 2\\
= &\frac{4 \sqrt{7} (2 \sqrt{7} - 1) \left\{ 3 \tan (t \pi / k) / 2 - \tan d_{1, 1} (t \pi / k) / 2 - 2 \sqrt{3} \tan (t \pi / k) / 2 \tan d_{1, 1} (t \pi / k) / 2 \right\}}{27 \left\{ 1 + ((2 \sqrt{7} + 5) / \sqrt{3}) \tan (t \pi / k) / 2 \right\} \left\{ 1 + ((2 \sqrt{7} - 1) / (3 \sqrt{3})) \tan d_{1, 1} (t \pi / k) / 2 \right\}},\\
&- (8 - 3 \sqrt{7}) \tan (\pi / 2 + \alpha_7 - (t \pi / k) \pi) / 2 - \tan (- 2 \pi / 3 + \pi / 2 + \alpha_7 - d_{1, 2} (t \pi / k)) / 2\\
= &\frac{2 \sqrt{7} (\sqrt{7} - 2) \left\{ \tan d_{1, 2} (t \pi / k) / 2 - 2 \tan (t \pi / k) / 2 + \sqrt{3} \tan (t \pi / k) / 2 \tan d_{1, 2} (t \pi / k) / 2 \right\}}{3 \left\{ 1 + ((\sqrt{7} - 2) / \sqrt{3}) \tan (t \pi / k) / 2 \right\} \left\{ 1 + ((2 \sqrt{7} + 5) / \sqrt{3}) \tan d_2 (t \pi / k) / 2 \right\}}.
\end{align*}
\end{allowdisplaybreaks}
Thus,
\begin{align*}
\tan \left( 2 \pi / 3 + \pi / 2 + \alpha_7 + d_{1, 1} (t \pi / k)\right) / 2 &< - ((11 - 4 \sqrt{7}) / 3) \tan \left( \pi / 2 + \alpha_7 - (t \pi / k) \right) / 2\\
 &< \tan \left( 2 \pi / 3 + \pi / 2 + \alpha_7 + 3 (t \pi / k) \right) / 2,\\
 \tan \left( - 2 \pi / 3 + \pi / 2 + \alpha_7 - 2 (t \pi / k)\right) / 2 &< - (8 - 3 \sqrt{7}) \tan \left( \pi / 2 + \alpha_7 - (t \pi / k) \right) / 2\\
 &< \tan \left( - 2 \pi / 3 + \pi / 2 + \alpha_7 - d_{1, 2} (t \pi / k) \right) / 2,
\end{align*}
where $d_{1, 1} < 3 / (1 + 2 \sqrt{3} \tan (t / 2) (\pi / k))$ and $d_{1, 2} < 2 / (1 + \sqrt{3} \tan (t / 2) (\pi / k))$. Furthermore, we have
\begin{align*}
&k \theta_1 / 2 = k (\pi / 2 + \alpha_7) / 2 - (t / 2) \pi \\
 &\qquad \Rightarrow
\begin{cases}
 (k / 3) \pi + k (\pi / 2 + \alpha_7) / 2 + d_{1, 1} (t / 2) \pi\\
 \qquad \qquad < k {\theta_{1, 1}}' / 2 < (k / 3) \pi + k (\pi / 2 + \alpha_7) / 2 + (3 t / 2) \pi,\\
 - (k / 3) \pi + k (\pi / 2 + \alpha_7) / 2 - t \pi\\
 \qquad \qquad < k {\theta_{1, 2}}' / 2 < - (k / 3) \pi + k (\pi / 2 + \alpha_7) / 2 - d_{1, 2} (t / 2) \pi.
\end{cases}
\end{align*}

On the other hand, when $\theta_2 = \alpha_7 - \pi / 6 + (t \pi / k) \pi$, then we have
\begin{allowdisplaybreaks}
\begin{align*}
&- ((4 + \sqrt{7}) / 3) \tan (\alpha_7 - \pi / 6 + (t \pi / k) \pi) / 2 - \tan (- 2 \pi / 3 + \alpha_7 - \pi / 6 - d_{2, 1} (t \pi / k)) / 2\\
= &\frac{(28 - 2 \sqrt{7}) \left\{ 2 \tan d_{2, 1} (t \pi / k) / 2 - 3 \tan (t \pi / k) / 2 + \sqrt{3} \tan (t \pi / k) / 2 \tan d_{2, 1} (t \pi / k) / 2 \right\}}{27 \left\{ 1 - ((\sqrt{7} - 2) / \sqrt{3}) \tan (t \pi / k) / 2 \right\} \left\{ 1 - ((2 \sqrt{7} - 1) / (3 \sqrt{3})) \tan d_{2, 1} (t \pi / k) / 2 \right\}},\\
&- (8 + 3 \sqrt{7}) \tan (\alpha_7 - \pi / 6 + (t \pi / k) \pi) / 2 - \tan (2 \pi / 3 + \alpha_7 - \pi / 6 + d_{2, 2} (t \pi / k)) / 2\\
= &\frac{(28 + 10 \sqrt{7}) \left\{ \tan (t \pi / k) / 2 - 2 \tan d_{2, 2} (t \pi / k) / 2 + 3 \sqrt{3} \tan (t \pi / k) / 2 \tan d_{2, 2} (t \pi / k) / 2 \right\}}{3 \left\{ 1 - ((\sqrt{7} - 2) / \sqrt{3}) \tan (t \pi / k) / 2 \right\} \left\{ 1 - ((2 \sqrt{7} + 5) / \sqrt{3}) \tan d_{2, 2} (t \pi / k) / 2 \right\}}.
\end{align*}
\end{allowdisplaybreaks}
Thus,
\begin{align*}
\tan \left( - 2 \pi / 3 + \alpha_7 - \pi / 6 - (3 / 2) (t \pi / k)\right) / 2 &< - ((4 + \sqrt{7}) / 3) \tan \left( \alpha_7 - \pi / 6 + (t \pi / k) \right) / 2\\
 &< \tan \left( - 2 \pi / 3 + \alpha_7 - \pi / 6 - d_{2, 1} (t \pi / k) \right) / 2,\\
 \tan \left( 2 \pi / 3 + \alpha_7 - \pi / 6 + d_{2, 2} (t \pi / k)\right) / 2 &< - (8 + 3 \sqrt{7}) \tan \left( \alpha_7 - \pi / 6 + (t \pi / k) \right) / 2\\
 &< \tan \left( 2 \pi / 3 + \alpha_7 - \pi / 6 + (1 / 2) (t \pi / k) \right) / 2,
\end{align*}
where $d_{2, 1} < 3 / (2 + \sqrt{3} \tan (t / 2) (\pi / k))$ and $d_{2, 2} < 1 / (2 + 3 \sqrt{3} \tan (t / 2) (\pi / k))$. Furthermore, we have
\begin{align*}
&k \theta_2 / 2 = k (\alpha_7 - \pi / 6) / 2 + (t / 2) \pi \\
 &\qquad \Rightarrow
\begin{cases}
 - (k / 3) \pi + k (\alpha_7 - \pi / 6) / 2 - (3 t / 4) \pi\\
 \qquad \qquad < k {\theta_{2, 1}}' / 2 < - (k / 3) \pi + k (\alpha_7 - \pi / 6) / 2 - d_{2, 1} (t / 2) \pi,\\
 (k / 3) \pi + k (\alpha_7 - \pi / 6) / 2 + d_{2, 2} (t / 2) \pi\\
 \qquad \qquad < k {\theta_{2, 2}}' / 2 < (k / 3) \pi + k (\alpha_7 - \pi / 6) / 2 + (t / 4) \pi.
\end{cases}
\end{align*}

Note that
\begin{equation*}
(k / 3) \pi \equiv
\begin{cases}
0 & (k \equiv 0 \pmod{6})\\
2 \pi / 3 & (k \equiv 2 \pmod{6})\\
4 \pi / 3 & (k \equiv 4 \pmod{6})
\end{cases}, \quad
- (k / 3) \pi \equiv
\begin{cases}
0 & (k \equiv 0 \pmod{6})\\
4 \pi / 3 & (k \equiv 2 \pmod{6})\\
2 \pi / 3 & (k \equiv 4 \pmod{6})
\end{cases},
\end{equation*}
modulo $2 \pi$. Furthermore, in the limit as $k$ tends to $\infty$ or in the limit as $t$ tends to $0$, $d_{1, 1}$ (resp. $d_{1, 2}$, $d_{2, 1}$, and $d_{2, 2}$) tends to $3$ (resp. $2$, $3 / 2$, and $1 / 2$).

Recall that $\alpha_{7, k} \equiv k (\pi / 2 + \alpha_7) / 2 \pmod{\pi}$. Then, we define ${\alpha_{7, k, n}}' \equiv k {\theta_{1, n}}' / 2 - (k (\pi / 2 + \alpha_7) / 2 - \alpha_{7, k} ) \pmod {2 \pi}$ for $n = 1, 2$ and for $\theta_1 = \pi / 2 + \alpha_7 - (t \pi / k)$.

Similarly, we define $\beta_{7, k} \equiv k (\alpha_7 - \pi / 6) / 2 \pmod{\pi}$ and ${\beta_{7, k, n}}' \equiv k {\theta_{2, n}}' / 2 - (k (\alpha_7 - \pi / 6) / 2 - \beta_{7, k} ) \pmod {2 \pi}$ for $n = 1, 2$ and for $\theta_2 = \alpha_7 - \pi / 6 + (t \pi / k)$.\\

\subsection{Absolute values of some terms}

\subsubsection{$\Gamma_0^{*}(5)$}

First, we observe that $|2 e^{i \theta_1 / 2} + \sqrt{5} e^{- i \theta_1 / 2}|^2 = 9 + 4 \sqrt{5} \cos \theta_1$. Let $f_1(k) := 9 + 4 \sqrt{5} \cos (\pi / 2 + \alpha_5 - (t \pi / k)) - (1 + 4 t \pi / k)$, then
\begin{align*}
\lim_{k \to \infty} f_1(k) &= 8 + 4 \sqrt{5} \cos (\pi / 2 + \alpha_5) = 0,\\
{f_1}'(k) &= - (4 t \pi / k^2) \{ \sqrt{5} \sin (\pi / 2 + \alpha_5 - (t \pi / k)) - 1 \}\\
 &\leqslant - (4 t \pi / k^2) \{ \sqrt{5} \sin (\pi / 2 + \alpha_5) - 1 \} = 0.
\end{align*}
Thus we have $f_1(k) \geqslant 0$ for every positive integer $k$. On the other hand, let $f_2(k) := e^{4 t \pi / k} - (9 + 4 \sqrt{5} \cos (\pi / 2 + \alpha_5 - (t \pi / k)))$, then
\begin{align*}
\lim_{k \to \infty} f_2(k) &= 1 - (9 + 4 \sqrt{5} \cos (\pi / 2 + \alpha_5)) = 0,\\
{f_2}'(k) &= (4 t \pi / k^2) \{ \sqrt{5} \sin (\pi / 2 + \alpha_5 - (t \pi / k)) - e^{4 t \pi / k} \}.
\end{align*}
Moreover, let $g(k) := \sqrt{5} \sin (\pi / 2 + \alpha_5 - (t \pi / k)) - e^{4 t \pi / k}$, then
\begin{align*}
\lim_{k \to \infty} g(k) &= \sqrt{5} \sin (\pi / 2 + \alpha_5) - 1 = 0,\\
g'(k) &= (4 t \pi / k^2) \{ e^{4 t \pi / k} + (\sqrt{5} / 4) \cos (\pi / 2 + \alpha_5 - (t \pi / k)) \}\\
 &\geqslant (4 t \pi / k^2) \{ 1 + (\sqrt{5} / 4) \cos (\pi / 2 + \alpha_5) \} = 0.
\end{align*}
Thus we have $f_2(k) \geqslant 0$ for every positive integer $k$. In conclusion, we have
\begin{equation*}
1 + 4 t \pi / k \leqslant 9 + 4 \sqrt{5} \cos (\pi / 2 + \alpha_5 - (t \pi / k)) \leqslant e^{4 t \pi / k}
\end{equation*}
for every positive integer $k$.

Similarly, for $|- e^{i \theta_2 / 2} + \sqrt{5} e^{- i \theta_2 / 2}|^2 / 4 = (6 - 2 \sqrt{5} \cos \theta_2) / 4$, we have
\begin{equation*}
1 + t \pi / k \leqslant (6 - 2 \sqrt{5} \cos (\alpha_5 + (t \pi / k))) / 4 \leqslant e^{t \pi / k}
\end{equation*}
for every positive integer $k$.\\

\subsubsection{$\Gamma_0^{*}(7)$}

First, we observe that $|2 e^{i \theta_1 / 2} + \sqrt{7} e^{- i \theta_1 / 2}|^2 = 11 + 4 \sqrt{7} \cos \theta_1$, $|3 e^{i \theta_1 / 2} + \sqrt{7} e^{- i \theta_1 / 2}|^2 = 16 + 6 \sqrt{7} \cos \theta_1$, $|- e^{i \theta_2 / 2} + \sqrt{7} e^{- i \theta_2 / 2}|^2 / 4 = (8 - 2 \sqrt{7} \cos \theta_2) / 4$, and $|- 3 e^{i \theta_2 / 2} + \sqrt{7} e^{- i \theta_2 / 2}|^2 / 4 = (16 - 6 \sqrt{7} \cos \theta_2) / 4$.

Similar to the previous subsubsection, we have
\begin{align*}
1 + 2 \sqrt{3} \: t \pi / k &\leqslant 11 + 4 \sqrt{7} \cos (\pi / 2 + \alpha_7 - (t \pi / k)) \leqslant e^{2 \sqrt{3} \: t \pi / k},\\
1 + 3 \sqrt{3} \: t \pi / k &\leqslant 16 + 6 \sqrt{7} \cos (\pi / 2 + \alpha_7 - (t \pi / k)) \leqslant e^{3 \sqrt{3} \: t \pi / k},\\
1 + (\sqrt{3} / 2) \: t \pi / k &\leqslant (8 - 2 \sqrt{7} \cos (\alpha_7 - \pi / 6 + (t \pi / k))) / 4, \\
1 + (3 \sqrt{3} / 2) \: t \pi / k &\leqslant (16 - 6 \sqrt{7} \cos (\alpha_7 - \pi / 6 + (t \pi / k))) / 4 \leqslant e^{(3 \sqrt{3} / 2) \: t \pi / k}
\end{align*}
for every positive integer $k$. Also, we can show
\begin{equation*}
(8 - 2 \sqrt{7} \cos (\alpha_7 - \pi / 6+ (t \pi / k))) / 4 \leqslant 1 + (\sqrt{3} / 2) (t \pi / k) + (1 / 2) (t \pi / k)^2
\end{equation*}
for every positive integer $k$.

Furthermore, assume that $t / k \leqslant 1 / 10$, which is satisfied when $k \geqslant 10$ for example. Let $f(s) := e^{(\sqrt{3} / 2) \: s \pi} - (8 - 2 \sqrt{7} \cos (\alpha_7 - \pi / 6 + s \pi) / 4$, then
\begin{align*}
f'(s) &:= \left( \sqrt{3} \pi / 2 \right) \left\{ e^{(\sqrt{3} / 2) \: s \pi} - (\sqrt{7} / \sqrt{3}) \sin (\alpha_7 - \pi / 6 + s \pi) \right\},\\
f''(s) &:= \left( \sqrt{3} \pi / 2 \right)^2 \left\{ e^{(\sqrt{3} / 2) \: s \pi} - (2 \sqrt{7} / 3) \cos (\alpha_7 - \pi / 6 + s \pi) \right\},\\
f'''(s) &:= \left( \sqrt{3} \pi / 2 \right)^3 \left\{ e^{(\sqrt{3} / 2) \: s \pi} + ((4 \sqrt{7}) / (3 \sqrt{3})) \sin (\alpha_7 - \pi / 6 + s \pi) \right\},
\end{align*}
and then
\begin{align*}
f'''(s) &> 0 \quad \text{for all} \quad s > 0,\\
f''(0) &< 0 \quad \text{and} \quad f''(1) > 0,\\
f'(0) &= 0 \quad \text{and} \quad f'(1) > 0,\\
f(0) &= 0 \quad \text{and} \quad f(1 / 10) = - 0.0038812... < 0.
\end{align*}
Thus, we have
\begin{equation*}
e^{(\sqrt{3} / 2) \: (t \pi / k)} \leqslant (8 - 2 \sqrt{7} \cos (\alpha_7 - \pi / 6 + (t \pi / k)) / 4
\end{equation*}
for every positive number $t$ and $k$ such that $t / k \leqslant 1 / 10$.\\

\subsection{Algorithm} \label{subsec-alg}

In this subsection, we consider the following bound:
\begin{equation}
|R_{p, n}^{*}| < 2 c_0 \quad \text{for every} \; k \geqslant k_0, \label{bound-c0}
\end{equation}
for some $c_0 > 0$ and an even integer $k_0$. Furthermore, we will proceed to present an algorithm to show the above bound.

Let $\Lambda$ be an index set, and let us write
\begin{equation*}
|R_{p, n}^{*}| \leqslant 2 {\textstyle \sum_{\lambda \in \Lambda}} \; e_{\lambda}^k \; v_k(c_{\lambda}, d_{\lambda}, \theta)
\end{equation*}
by an application of the RSD Method, where the factor ``$2$'' arises from the fact $v_k(c, d, \theta) = v_k(- c, - d, \theta)$. Furthermore, let $I$ be a finite subset of $\Lambda$ such that $e_i^k \; v_k(c_i, d_i, \theta)$ does not tend to $0$ in the limit as $k$ tends to $\infty$ for every $i \in I$, and assume $I \subset \mathbb{N}$. Then, we define $X_i := e_i^{- 2} \; v_k(c_i, d_i, \theta)^{- 2 / k}$ for every $i \in I$.

Assume that for every $i \in I$ and $k \geqslant k_0$ and for some ${c_i}'$ and $u_i$,
\begin{gather*}
\left| Re \left\{ e_i^k \; \left( c_i e^{i \theta / 2} + \sqrt{p} d_i e^{- i \theta / 2} \right)^{- k} \right\} \right| \leqslant {c_i}' X_i^{- k / 2},\\
 X_i^{- k / 2} \geqslant 1 + u_i (\pi / k), \qquad 2 {\textstyle \sum_{\lambda \in \Lambda \setminus I}} \; e_{\lambda}^k \; v_k(c_{\lambda}, d_{\lambda}, \theta) \leqslant b (1 / s)^{k / 2},
\end{gather*}
and let the number $t > 0$ be given.\\

\begin{trivlist}
\item[\bfseries Step 1.]``Determine the number $a_1$.''

First, in order to show the bound (\ref{bound-c0}), we wish to make use of the following bound:
\begin{equation}
{\textstyle \sum_{i \in I}} \; {c_i}' X_i^{- k / 2} < c_0 - a_1 (t \pi / k)^2 \label{bound-a10}
\end{equation}
for every $i \in I$ and $k \geqslant k_0$ and for some $a_1 > 0$.

To show the bound (\ref{bound-c0}) using the above bound (\ref{bound-a10}), we need $b (1 / s)^{k / 2} < a_1 (t \pi / k)^2$ for every $k \geqslant k_0$. Define $f(k) := s^{k / 2} / b - k^2 / (2 a_1 t^2 \pi^2)$. If we have $k_0 \log s > 4$ and
\begin{equation}
a_1 > (b k_0^2) \; / \; (2 s^{k_0/ 2} t^2 \pi^2), \label{bound-a1}
\end{equation}
then we have $f(k_0) > 0$, $f'(k_0) > 0$, and $f''(k_0) > 0$. In the present paper, we always have $k_0 \log s > 4$. Thus, it is enough to consider bound (\ref{bound-a1}).\\

\item[\bfseries Step 2.]``Determine the number $c_{0, i}$ and $a_{1, i}$.''

Second, in order to show the bound (\ref{bound-a10}), we wish to use the following bounds:
\begin{equation}
{c_i}' X_i^{- k / 2} < c_{0, i} - a_{1, i} (t \pi / k)^2 \label{bound-Xi}
\end{equation}
for every $i \in I$ and $k \geqslant k_0$ and for some $c_{0, i} > 0$, $a_{1, i} > 0$.

We determine $c_{0, i}$ and $a_{1, i}$ such that 
\begin{equation}
c_{0, i} > 0, \quad a_{1, i} > 0, \quad c_0 = {\textstyle \sum_{i \in I}} \; c_{0, i} , \quad \text{and} \quad a_1 = {\textstyle \sum_{i \in I}} \; a_{1, i}
\end{equation}\quad

\item[\bfseries Step 3.]``Determine a discriminant $Y_i$ for every $i \in I$.''

Finally, for the bound (\ref{bound-Xi}), we consider the following sufficient conditions:
\begin{equation*}
X_i^{k / 2} > c_i + a_{2, i} (t \pi / k)^2, \quad X_i > a_{3, i} + a_{4, i} (t \pi / k)^2.
\end{equation*}
For the former bound, it is enough to show that
\begin{equation}
c_i = {c_i}' / c_{0, i}, \quad a_{2, i} > c_i^2 (a_{1, i} / {c_i}') \; / \; \{1 - c_i (a_{1, i} / {c_i}') (t \pi / k_0)^2\}, \label{bound-a2}
\end{equation}
while for the latter bound, it is enough to show that
\begin{equation*}
a_{3, i} = c_i^{2 / k_0}, \quad a_{4, i} = \left( (2 a_{2, i}) / (c_i k_0) \right) c_i^{2 / k_0}.
\end{equation*}

Because we have
\begin{equation*}
c_i^{2 / k} \leqslant 1 + 2 (\log c_i) / k + 2 (\log c_i)^2 c_i^{2 / k} / k^2,
\end{equation*}
\begin{align}
X_i - \left( a_{3, i} + a_{4, i} \left( t \frac{\pi}{k} \right)^2 \right)
& \geqslant \frac{1}{k} \left\{ u_i \pi - 2 \log c_i - 2 (\log c_i)^2 c_i^{2 / k} \frac{1}{k_0} - \frac{2 a_{2, i} t^2 \pi^2}{c_i} c_i^{2 / k_0} \frac{1}{k_0^2} \right\} \notag \\
&\quad =: \frac{1}{k} \times Y_i.
\end{align}

In conclusion, if we have $Y_i > 0$, then we can show bounds (\ref{bound-Xi}), (\ref{bound-a10}), and (\ref{bound-c0}).
\end{trivlist}

Note that the above bounds are sufficient conditions, but are not always necessary.

\newpage

\section{$\Gamma_0^{*}(5)$ (For Conjecture \ref{conj-g0s5})}

To prove Conjecture \ref{conj-g0s5} is much more difficult than the proof of theorems for $\Gamma_0^{*}(2)$ and $\Gamma_0^{*}(3)$. Difficulties arise in particular due to the argument $Arg(\rho_{5, 2})$, which is not a rational multiple of $\pi$.\\

The proof of all the lemmas of this section are presented at the end of this section.\\

\subsection{$E_{k, 5}^{*}$ of low weights}
First, we have the following Lemma.

\begin{lemma}
Let $k \geqslant 4$ be an even integer.
\def\labelenumi{(\arabic{enumi})}
\begin{enumerate}
\item
\begin{equation*}
F_{k, 5}^{*} (\pi / 2)
\begin{cases}
> 0 & k \equiv 0 \pmod{8}\\
< 0 & k \equiv 4 \pmod{8}\\
= 0 & k \equiv 2 \pmod{4}
\end{cases}.
\end{equation*}

\item
\begin{equation*}
F_{k, 5}^{*} (\pi / 2 + \alpha_5) =
\begin{cases}
\frac{2 \cdot 5^{k / 2}}{5^{k / 2} + 1} \cos (k (\pi / 2 + \alpha_5) / 2) E_k(i) & k \equiv 0 \pmod{4}\\
0 & k \equiv 2 \pmod{4}
\end{cases}.
\end{equation*}
Furthermore, we have $E_k(i) > 0$ for every $k \geqslant 4$ such that $k \equiv 0 \pmod{4}$.

\item
For an even integer $k \geqslant 8$,
\begin{equation*}
F_{k, 5}^{*} (\pi)
\begin{cases}
> 0 & k \equiv 0 \pmod{8}\\
< 0 & k \equiv 4 \pmod{8}\\
= 0 & k \equiv 2 \pmod{4}
\end{cases}.
\end{equation*}

\item
For an even integer $k \geqslant 10$,
\begin{equation*}
F_{k, 5}^{*} (\theta) = 2 \cos(k \theta / 2) + R_{5, 5 \pi / 6}^{*} \quad \text{for} \: \theta \in [\pi / 2, 5 \pi / 6],
\end{equation*}
where $|R_{5, 5 \pi / 6}^{*}| < 2$.
\end{enumerate}
\def\labelenumi{\arabic{enumi}.}\label{lem-f5s1}
\end{lemma}

The proof of the above lemma is presented in Subsection \ref{subsec-pf51}. Now, we have the following proposition:

\begin{proposition}
The location of the zeros of the Eisenstein series $E_{k, 5}^{*}$ in $\mathbb{F}^{*}(5)$ for $4 \geqslant k \geqslant 10$ are follows$:$
\begin{center}
\begin{tabular}{ccccccc}
$k$ & $v_{\infty}$ & $v_{i / \sqrt{5}}$ & $v_{\rho_{5, 1}}$ & $v_{\rho_{5, 2}}$ & $V_{5, 1}^{*}$ & $V_{5, 2}^{*}$\\
\hline
$4$ & $0$ & $0$ & $0$ & $0$ & $1$ & $0$\\
$6$ & $0$ & $1$ & $1$ & $1$ & $0$ & $0$\\
$8$ & $0$ & $0$ & $0$ & $0$ & $1$ & $1$\\
$10$ & $0$ & $1$ & $1$ & $1$ & $1$ & $0$\\
\hline
\end{tabular}
\end{center}
where $V_{5, n}^{*}$ denotes the number of simple zeros of the Eisenstein series $E_{k, 5}^{*}$ on the arc $A_{5, n}^{*}$ for $n = 1, 2$. \label{prop-ek5s-lw}
\end{proposition}

\begin{proof}\quad
\begin{trivlist}
\item[($k = 4$)]
We have $F_{4, 5}^{*}(\pi / 2) < 0$ by Lemma \ref{lem-f5s1} (1), and we have $F_{4, 5}^{*}(\pi / 2 + \alpha_5) > 0$ by Lemma \ref{lem-f5s1} (2) because $\cos (2 (\pi / 2 + \alpha_5)) = 3 / 5 > 0$. Thus, $E_{4, 5}^{*}$ has at least one zero on $A_{5, 1}^{*}$. Furthermore, by the valence formula for $\Gamma_0^{*}(5)$ (Proposition \ref{prop-vf-g0s5}), $E_{4, 5}^{*}$ has no other zeros.

\item[($k = 6$)]
By previous subsubsection, we have $v_{i / \sqrt{5}} (E_{6, 5}^{*}) \geqslant 1$, $v_{\rho_{5, 1}} (E_{6, 5}^{*}) \geqslant 1$, and $v_{\rho_{5, 2}} (E_{6, 5}^{*}) \geqslant 1$. Furthermore, by the valence formula for $\Gamma_0^{*}(5)$, $E_{6, 5}^{*}$ has no other zeros.

\item[($k = 8$)]
We have $F_{8, 5}^{*}(\pi / 2) > 0$ by Lemma \ref{lem-f5s1} (1), and we have $F_{8, 5}^{*}(\pi / 2 + \alpha_5) > 0$ by Lemma \ref{lem-f5s1} (2) because $\cos (4 (\pi / 2 + \alpha_5)) = - 7 / 25 < 0$, and we have $F_{8, 5}^{*}(\pi) > 0$ by Lemma \ref{lem-f5s1} (3). Thus, $E_{8, 5}^{*}$ has at least two zeros on each arc $A_{5, 1}^{*}$ and $A_{5, 2}^{*}$. Furthermore, by the valence formula for $\Gamma_0^{*}(5)$, $E_{8, 5}^{*}$ has no other zeros.

\item[($k = 10$)]
We have $F_{10, 5}^{*}(3 \pi / 5) < 0$ and $F_{10, 5}^{*}(4 \pi / 5) > 0$ by Lemma \ref{lem-f5s1} (4). Thus, $E_{10, 5}^{*}$ has at least one zero on $A_{5, 1}^{*}$. In addition, by previous subsubsection, we have $v_{i / \sqrt{5}} (E_{10, 5}^{*}) \geqslant 1$, $v_{\rho_{5, 1}} (E_{10, 5}^{*}) \geqslant 1$, and $v_{\rho_{5, 2}} (E_{10, 5}^{*}) \geqslant 1$. Furthermore, by the valence formula for $\Gamma_0^{*}(5)$, $E_{10, 5}^{*}$ has no other zeros.
\end{trivlist}
\end{proof}\quad

\subsection{All but at most $2$ zeros}\label{subsec-g0s5-ab2}
\begin{lemma}We have the following bounds$:$\\
``We have $|R_{5, 1}^{*}| < 2 \cos({c_0}' \pi)$ for $\theta_1 \in [\pi / 2, \pi / 2 + \alpha_5 - t \pi / k]$''\\
\quad$(1)$ For $k \geqslant 12$, $({c_0}', t) = (1 / 3, 1 / 6)$.\\
\quad$(2)$ For $k \geqslant 58$, $({c_0}', t) = (33 / 80, 9 / 40)$.\\
``We have $|R_{5, 2}^{*}| < 2 \cos({c_0}' \pi)$ for $\theta_2 \in [\alpha_5 + t \pi / k, \pi / 2]$''\\
\quad$(3)$ For $k \geqslant 12$, $({c_0}', t) = (0, 1 / 2)$.\\
\quad$(4)$ For $k \geqslant 22$, $({c_0}', t) = (1 / 3, 1 / 2)$.\\
\quad$(5)$ For $k \geqslant 46$, $({c_0}', t) = (7 / 30, 1 / 5)$.
\label{lem-f5s2}
\end{lemma}

When $4 \mid k$, by the valence formula for $\Gamma_0^{*}(5)$ and Proposition \ref{prop-bd_ord_5}, we have at most $k / 4$ zeros on the arc $A_5^{*}$. We have $k / 4 + 1$ {\it integer points} ({\it i.e.} $\cos\left( k \theta / 2 \right) = \pm 1$) in the interval $[\pi / 2, \pi]$. If we prove that the sign of $F_{k, 5}^{*} (\theta)$ is equal to $\cos\left( k \theta / 2 \right)$ at every integer point, then we can prove Conjecture {\ref{conj-g0s5} for the case $4 \mid k$. By the above lemma's conditions (1) and (3), we can prove $|R_{5, 1}^{*}| < 2$ or $|R_{5, 2}^{*}| < 2$ at all of the integer points which satisfy $\theta_1 \in [\pi / 2, \pi / 2 + \alpha_5 - \pi / (6 k)]$ or $\theta_2 \in [\alpha_5 + \pi / (2 k), \pi / 2]$, respectively. Then, we have all but at most $2$ zeros on $A_5^{*}$.

On the other hand, when $4 \nmid k$, we have at most $(k - 6) / 4$ zeros on the arc $A_5^{*}$. When $0 < \alpha_{5, k} < \pi / 2$, between the first integer point and the last one for $A_{5, 1}^{*}$ we have $(k \alpha_5 / 2 - \pi / 2 - \alpha_{5, k}) / \pi + 1$ integer points, then we expect $(k \alpha_5 / 2 - \pi / 2 - \alpha_{5, k}) / \pi$ zeros on the arc $A_{5, 1}^{*}$. Similarly, we expect $(k (\pi / 2 - \alpha_5) / 2 - \pi / 2 - (\pi - \beta_{5, k})) / \pi$ zeros on the arc $A_{5, 2}^{*}$. Furthermore, we have $\alpha_{5, k} + (\pi - \beta_{5, k}) = \pi / 2$. Thus, in agreement with expectations we have $(k - 6) / 4$ zeros on the arc $A_5^{*}$. By the above lemma's conditions (1) and (3), we can prove $|R_{5, 1}^{*}| < 2$ at every integer point less than last one for $A_{5, 1}^{*}$, and prove $|R_{5, 2}^{*}| < 2$ at every integer point greater than the first one for $A_{5, 2}^{*}$. Then, we have all but at most $2$ zeros on $A_5^{*}$.

Finally, when $\pi / 2 < \alpha_{5, k} < \pi$, we expect $(k - 10) / 4$ zeros between adjacent integer points, which are proved by the above lemma's conditions (1) and (3). Then, we have all but at most one zero are on $A_5^{*}$.

Thus, we have the following proposition:
\begin{proposition}
Let $k \geqslant 4$ be an even integer. All but at most $2$ of the zeros of $E_{k, 5}^{*}(z)$ in $\mathbb{F}^{*}(5)$ are on the arc $A_5^{*}$. \label{prop-g0s5-ab2}
\end{proposition}\quad

\subsection{The case $4 \mid k$}\label{subsec-g0s5-40}
When $\pi / 12 < \alpha_{5, k} < 3 \pi / 4$, we can prove the proposition at all of the integer points.

Now, we can write
\begin{gather*}
F_{k, 5, 1}^{*}(\theta_1) = 2 \cos\left( k \theta_1 / 2 \right) + 2 Re(2 e^{- i \theta_1 / 2} + \sqrt{5} e^{i \theta_1 / 2})^{- k} + {R_{5, 1}^{*}}',\\
F_{k, 5, 2}^{*}(\theta_2) = 2 \cos\left( k \theta_2 / 2 \right) + 2^k \cdot 2 Re(e^{- i \theta_2 / 2} - \sqrt{5} e^{i \theta_2 / 2})^{- k} + {R_{5, 2}^{*}}'.
\end{gather*}

When $0 < \alpha_{5, k} < \pi / 12$, the last integer point of $F_{k, 5, 1}^{*}(\theta_1)$ is in the interval $[\pi / 2 + \alpha_5 - \pi / (6 k), \pi / 2 + \alpha_5]$. We have $|{R_{5, 1}^{*}}'| < 2$ for $\theta_1 \in [\pi / 2, \pi / 2 + \alpha_5]$. Furthermore, because $0 < {\alpha_{5, k}}' < \pi / 6$ for $0 < t < 1 / 6$, we have $Sign\{\cos(k \theta_1 / 2)\} = Sign\{Re(2 e^{- i \theta_1 / 2} + \sqrt{5} e^{i \theta_1 / 2})^{- k}\}$ for $\theta_1 \in [\pi / 2 + \alpha_5 - \pi / (6 k), \pi / 2 + \alpha_5]$.

When $3 \pi / 4 < \alpha_{5, k} < \pi$, the first integer point of $F_{k, 5, 2}^{*}(\theta_2)$ is in the interval $[\alpha_5, \alpha_5 + \pi / (2 k)]$. We have $|{R_{5, 2}^{*}}'| < 2$ and $Sign\{\cos(k \theta_2 / 2)\} = Sign\{Re(e^{- i \theta_2 / 2} - \sqrt{5} e^{i \theta_2 / 2})^{- k}\}$ for $\theta_2 \in [\alpha_5, \alpha_5 + \pi / (2 k)]$.

Thus, we have the following proposition:
\begin{proposition}
Let $k \geqslant 4$ be an integer which satisfies $4 \mid k$. All of the zeros of $E_{k, 5}^{*}(z)$ in $\mathbb{F}^{*}(5)$ are on the arc $A_5^{*}$. \label{prop-g0s5-40}
\end{proposition}\quad

\subsection{The case $4 \nmid k$}\label{subsec-g0s5-41}

\subsubsection{The case $0 < \alpha_{5, k} < \pi / 2$}
At points such that $k \theta_1 / 2 = k (\pi / 2 + \alpha_5) / 2 - \alpha_{5, k} - \pi / 3$, we have $|R_{5, 1}^{*}| < 1$ by Lemma \ref{lem-f5s2} (1), and we have $2 \cos(k \theta_1 / 2) = \pm 1$. Then, we have at least one zero between the second to last integer point for $A_{5, 1}^{*}$ and the point $k \theta_1 / 2$. Similarly, by Lemma \ref{lem-f5s2} (4), we have at least one zero between the second integer point and the point $k \theta_2 / 2 = k \alpha_5 / 2 + \beta_{5, k} + \pi / 3$.\\

\subsubsection{The case $\pi / 2 < \alpha_{5, k} < \pi$}

We have the following lemmas:

\begin{lemma}Let $4 \nmid k$. We have the following bounds$:$\\
``When $(x / 180) \pi < \alpha_{5, k} < (y / 180) \pi$, we have $|R_{5, 1}^{*}| < 2 \cos({c_0}' \pi)$ for $k \geqslant k_0$ and $\theta_1 = \pi / 2 + \alpha_5 - t \pi / k$.''
\def\labelenumi{(\arabic{enumi})}
\begin{enumerate}
\item $(x, y) = (121, 126)$, $({c_0}', k_0, t) = (29 / 72, 100, 3 / 20)$.
\item $(x, y) = (120, 121)$, $({c_0}', k_0, t) = (23 / 60, 18, 1 / 10)$.
\item $(x, y) = (118.8, 120)$, $({c_0}', k_0, t) = (39 / 100, 100, 1 / 10)$.
\item $(x, y) = (118.1, 118.8)$, $({c_0}', k_0, t) = (691 / 1800, 100, 2 / 25)$.
\item $(x, y) = (117.7, 118.1)$, $({c_0}', k_0, t) = (638 / 1800, 86, 1 / 15)$.
\item $(x, y) = (117.45, 117.7)$, $({c_0}', k_0, t) = (151 / 400, 100, 3 / 50)$.
\item $(x, y) = (117.27, 117.45)$, $({c_0}', k_0, t) = (747 / 1800, 100, 1 / 20)$.
\item $(x, y) = (117.15, 117.27)$, $({c_0}', k_0, t) = (223 / 600, 100, 9 / 200)$.
\item $(x, y) = (117.06, 117.15)$, $({c_0}', k_0, t) = (1109 / 3000, 100, 1 / 25)$.
\item $(x, y) = (117, 117.06)$, $({c_0}', k_0, t) = (37 / 100, 100, 1 / 25)$.
\end{enumerate}
\def\labelenumi{\arabic{enumi}.}\label{lem-f5s3}
\end{lemma}

\begin{lemma}Let $4 \nmid k$. We have the following bounds$:$\\
``When $(x / 180) \pi < \alpha_{5, k} < (y / 180) \pi$, we have $|R_{5, 2}^{*}| < 2 \cos({c_0}' \pi)$ for $k \geqslant k_0$ and $\theta_2 = \alpha_5 + t \pi / k$.''
\def\labelenumi{(\arabic{enumi})}
\begin{enumerate}
\item $(x, y) = (114, 115.4)$, $({c_0}', k_0, t) = (199 / 900, 100, 4 / 25$.
\item $(x, y) = (115.4, 115.8)$, $({c_0}', k_0, t) = (61 / 300, 100, 3 / 25)$.
\item $(x, y) = (115.8, 116)$, $({c_0}', k_0, t) = (7 / 36, 100, 1 / 10)$.
\end{enumerate}
\def\labelenumi{\arabic{enumi}.}\label{lem-f5s4}
\end{lemma}

Now, we expect one more zero between the last integer point for $A_{5, 1}^{*}$ and the first one for $A_{5, 2}^{*}$. Then, we consider the following cases.

\paragraph{(i)} ``The case $7 \pi / 10 < \alpha_{5, k} < \pi$''
\begin{itemize}
\item When $3 \pi / 4 < \alpha_{5, k} < \pi$, we can use Lemma \ref{lem-f5s2} (1).
\item When $7 \pi / 10 < \alpha_{5, k} < 3 \pi / 4$, we can use Lemma \ref{lem-f5s2} (2).
\end{itemize}
For each case, we consider the point such that $k \theta_1 / 2 = k (\pi / 2 + \alpha_5) / 2 - \alpha_{5, k} + \pi - {c_0}' \pi$. We have $\alpha_{5, k} - \pi + {c_0}' \pi > (t / 2) \pi$ and $|R_{5, 1}^{*}| < 2 \cos({c_0}' \pi)$, and we have $2 \cos(k \theta_1 / 2) = \pm 2 \cos({c_0}' \pi)$. Then, we have at least one zero between the second to the last integer point for $A_{5, 1}^{*}$ and the point $k \theta_1 / 2$.\\

\paragraph{(ii)} ``The case $\pi / 2 < \alpha_{5, k} < 19 \pi / 30$'' 
\begin{itemize}
\item When $\pi / 2 < \alpha_{5, k} < 7 \pi / 12$, we can use Lemma \ref{lem-f5s2} (4).
\item When $7 \pi / 12 < \alpha_{5, k} < 19 \pi / 30$, we can use Lemma \ref{lem-f5s2} (5).
\end{itemize}
Similar to the case (i), we consider the point such that $k \theta_2 / 2 = k \alpha_5 / 2 - \beta_{5, k} + {c_0}' \pi$ for each case.\\

\paragraph{(iii)} ``The case $13 \pi / 20 < \alpha_{5, k} < 7 \pi / 10$''
We can use Lemma \ref{lem-f5s3}. For each case, we consider the point such that $k \theta_1 / 2 = k (\pi / 2 + \alpha_5) / 2 - (t / 2) \pi$. We have $\alpha_{5, k} - \pi + {c_0}' \pi > (t / 2) \pi$ and $|R_{5, 1}^{*}| < 2 \cos({c_0}' \pi)$, and we have $|2 \cos(k \theta_1 / 2)| > 2 \cos({c_0}' \pi)$. Then, we have at least one zero between the second to last integer point for $A_{5, 1}^{*}$ and the point $k \theta_1 / 2$.\\

\paragraph{(iv)} ``The case $19 \pi / 30 < \alpha_{5, k} < 29 \pi / 45$''
Similar to the previous case, we can use Lemma \ref{lem-f5s4}. We consider the point such that $k \theta_2 / 2 = k \alpha_5 / 2 + (t / 2) \pi$ for each case.\\

In conclusion, we have the following proposition:
\begin{proposition}
Let $k \geqslant 4$ be an integer which satisfies $4 \nmid k$, and let $\alpha_{5, k} \in [0, \pi]$ be the angle that satisfies $\alpha_{5, k} \equiv k (\pi / 2 + \alpha_5) / 2 \pmod{\pi}$. If we have $\alpha_{5, k} < 29 \pi / 45$ or $13 \pi / 20 < \alpha_{5, k}$, then all of the zeros of $E_{k, 5}^{*}(z)$ in $\mathbb{F}^{*}(5)$ are on the arc $A_5^{*}$. Otherwise, all but at most one zero of $E_{k, 5}^{*}(z)$ in $\mathbb{F}^{*}(5)$ are on $A_5^{*}$\label{prop-g0s5-41}
\end{proposition}\quad

\subsection{The remaining case ``$4 \nmid k$ and $29 \pi / 45 < \alpha_{5, k} < 13 \pi / 20$''} \label{subsec-rest5}
In the previous subsection, we left one zero between the last integer point for $A_{5, 1}^{*}$ and the first one for $A_{5, 2}^{*}$ for the case of ``$4 \nmid k$ and $29 \pi / 45 < \alpha_{5, k} < 13 \pi / 20$''. In the Lemma \ref{lem-f5s3} and \ref{lem-f5s4}, the width $|x - y|$ becomes smaller as the intervals of the bounds approach the interval $[29 \pi / 45, 13 \pi / 20]$. However, it would appear that we require the width $|x - y|$ to become even smaller in this limit if we are to prove the proposition for the remaining interval $[29 \pi / 45, 13 \pi / 20]$. Furthermore, we may need infinite bounds such as we see in the lemmas \ref{lem-f5s3} and \ref{lem-f5s4}. Thus, we cannot prove the result for this remaining case in a similar manner. However, when $k$ is large enough, there is some possibility that we can prove the proposition for the remaining case.

Let $29 \pi / 45 < \alpha_{5, k} < 13 \pi / 20$, and let $t > 0$ be small enough. Then, we have $\pi / 2 < \alpha_{5, k} - (t / 2) \pi < \pi$ and $3 \pi / 2 < \pi + \alpha_{5, k} + d_1 (t / 2) \pi < {\alpha_{5, k}}' < \pi + \alpha_{5, k} + (t / 2) \pi < 2 \pi$. Moreover, we can easily show that $1 + 4 t (\pi / k) \leqslant v_k(2, 1, \theta_1)^{- 2 / k} \leqslant e^{4 t (\pi / k)}$. Thus, we have
\begin{align*}
&- \cos(\alpha_{5, k} - (t / 2) \pi) - \cos(\pi + \alpha_{5, k} + d_1 (t / 2) \pi) \cdot e^{- 2 \pi t}\\
 &\qquad > |\cos(k \theta_1 / 2)| - \left| Re\left\{ \left( 2 e^{i \theta_1 / 2} + \sqrt{5} e^{- i \theta_1 / 2} \right)^{- k} \right\} \right|\\
 &\qquad > - \cos(\alpha_{5, k} - (t / 2) \pi) - \cos(\pi + \alpha_{5, k} + (t / 2) \pi) \cdot (1 + 4 t (\pi / k))^{- k / 2}.
\end{align*}
We denote the upper bound by $A_1$ and the lower bound by $B_1$. Furthermore, we define ${A_1}' := A_1 / \cos(\pi + \alpha_{5, k} + d_1 (t / 2) \pi)$ and ${B_1}' := B_1 / \cos(\pi + \alpha_{5, k} + (t / 2) \pi)$. Then, we have
\begin{align*}
\frac{\partial}{\partial t} {A_1}' &= (\pi / 2) \frac{\sin(\alpha_{5, k} - (t / 2) \pi) \cos(\alpha_{5, k} + d_1 (t / 2) \pi) + d_1 \sin(\alpha_{5, k} + d_1 (t / 2) \pi) \cos(\alpha_{5, k} - (t / 2) \pi)}{\cos^2(\alpha_{5, k} + d_1 (t / 2) \pi)}\\
 &\quad+ 2 \pi e^{- 2 \pi t},\\
\frac{\partial}{\partial t} {B_1}' &= (\pi / 2) \frac{\sin(2 \alpha_{5, k})}{\cos^2(\alpha_{5, k} + (t / 2) \pi)} + 2 \pi (1 + 4 t (\pi / k))^{- k / 2}.
\end{align*}
Considering $\lim_{t \to 0} d_1 = 1$, we assume $d_1 |_{t=0} = 1$. First, we have $A_1 |_{t = 0} = B_1 |_{t = 0} = 0$. Second, we have $\frac{\partial}{\partial t} {A_1}' |_{t = 0} = \frac{\partial}{\partial t} {B_1}' |_{t = 0} = \pi (\tan \alpha_{k, 5} + 2)$. Finally, since $\pi (\tan (\pi - \alpha_5) + 2) = 0$, we have $B_1 > 0$ if $\alpha_{5, k} > \pi - \alpha_5$, and we have $A_1 < 0$ if $\alpha_{5, k} < \pi - \alpha_5$ for small enough $t$.

Similarly, we have $0 < \beta_{5, k} + (t / 2) \pi < \pi / 2$ and $\pi < \pi + \beta_{5, k} - (t / 2) \pi < {\beta_{5, k}}' < \pi + \beta_{5, k} - d_2 (t / 2) \pi < 3 \pi / 2$. Thus
\begin{align*}
&\cos(\beta_{5, k} + (t / 2) \pi) + \cos(\pi + \beta_{5, k} - d_2 (t / 2) \pi) \cdot e^{- \pi t / 2}\\
 &\qquad > |\cos(k \theta_2 / 2)| - \left| Re\left\{ 2^k \left( e^{i \theta_2 / 2} - \sqrt{5} e^{- i \theta_2 / 2} \right)^{- k} \right\} \right|\\
 &\qquad > \cos(\beta_{5, k} + (t / 2) \pi) + \cos(\pi + \beta_{5, k} - (t / 2) \pi) \cdot (1 + t (\pi / k))^{- k / 2}.
\end{align*}
We denote the upper bound by $A_2$ and the lower bound by $B_2$. Furthermore, we define ${A_2}' := A_2 / \cos(\beta_{5, k} - d_2 (t / 2) \pi)$ and ${B_2}' := B_2 / \cos(\beta_{5, k} - (t / 2) \pi)$. We assume $d_2 |_{t=0} = 1$. First, we have $A_2 |_{t = 0} = B_2 |_{t = 0} = 0$. Second, we have $\frac{\partial}{\partial t} {A_2}' |_{t = 0} = \frac{\partial}{\partial t} {B_2}' |_{t = 0} = \pi (1 / 2 - \tan \beta_{k, 5})$. Finally, since $\pi (1 / 2 - \tan (\pi / 2 - \alpha_5)) = 0$ and $\alpha_{5, k} = \beta_{5, k} + \pi / 2$, we have $B_2 > 0$ if $\alpha_{5, k} < \pi - \alpha_5$, and we have $A_2 < 0$ if $\alpha_{5, k} > \pi - \alpha_5$ for small enough $t$.

In conclusion, if $4 \nmid k$ is large enough, then $|{R_{5, 1}^{*}}'|$ and $|{R_{5, 2}^{*}}'|$ is small enough, and then we have one more zero on the arc $A_{5, 1}^{*}$ when $\alpha_{5, k} > \pi - \alpha_5$, and we have one more zero on the arc $A_{5, 2}^{*}$ when $\alpha_{5, k} < \pi - \alpha_5$. However, if $k$ is small, the proof is unclear.\\

The remaining subsections in this section detail the proofs of lemmas \ref{lem-f5s1}, ..., \ref{lem-f5s4}.\\

\subsection{Proof of Lemma \ref{lem-f5s1}} \label{subsec-pf51}\quad

\begin{trivlist}
\item[(1)]
Let $k \geqslant 4$ be an even integer divisible by $4$.

First, we consider the case $N = 1$. Then, we can write:
\begin{equation*}
F_{k, 5}^{*} (\pi / 2) = F_{k, 5, 1}^{*}(\pi / 2) = 2 \cos(k \pi / 4) + R_{5, \pi / 2}^{*}
\end{equation*}
where $R_{5, \pi / 2}^{*}$ denotes the remaining terms.

We have $v_{k}(c, d, \pi / 2) = 1 / (c^2 + 5 d^2)^{k/2}$. Similar to the application of the RSD Method, we will consider the following cases: $N = 2$, $5$, $10$, $13$, $17$, and $N \geqslant 25$. We have
\begin{allowdisplaybreaks}
\begin{align*}
&\text{When $N = 2$,}&
v_k(1, 1, \pi / 2) &\leqslant (1 / 6)^{k/2}.\\
&\text{When $N = 5$,}&
v_k(1, 2, \pi / 2) &\leqslant (1 / 21)^{k/2},
&v_k(2, 1, \pi / 2) &\leqslant (1 / 3)^{k}.\\
&\text{When $N = 10$,}&
v_k(1, 3, \pi / 2) &\leqslant (1 / 46)^{k/2},
&v_k(3, 1, \pi / 2) &\leqslant (1 / 14)^{k/2}.\\
&\text{When $N = 13$,}&
v_k(2, 3, \pi / 2) &\leqslant (1 / 7)^{k},
&v_k(3, 2, \pi / 2) &\leqslant (1 / 29)^{k/2}.\\
&\text{When $N = 17$,}&
v_k(1, 4, \pi / 2) &\leqslant (1 / 21)^{k/2},
&v_k(4, 1, \pi / 2) &\leqslant (1 / 3)^{2 k}.\\
&\text{When $N \geqslant 25$,}&
c^2 + 5 d^2 &\geqslant N,
\end{align*}
\end{allowdisplaybreaks}
and the number of terms with $c^2 + d^2 = N$ is not more than $(96 / 25) N^{1/2}$ for $N \geqslant 25$. Then
\begin{equation*}
|R_{5, \pi / 2}^{*}|_{N \geqslant 25}
\leqslant \frac{192}{25(k-3)} \left(\frac{1}{24}\right)^{(k - 3) / 2}.
\end{equation*}
Furthermore,
\begin{align*}
|R_{5, \pi / 2}^{*}|
 &\leqslant
 4 \left(\frac{1}{6}\right)^{k/2} + 4 \left(\frac{1}{3}\right)^{k}
 + \cdots + 4 \left(\frac{1}{3}\right)^{2 k}
 + \frac{192}{25(k-3)} \left(\frac{1}{24}\right)^{(k - 3) / 2},\\
 &\leqslant 1.77563... \quad (k \geqslant 4)
\end{align*}\quad

\item[(2)]
Let $k \geqslant 4$ be an even integer divisible by $4$.
\begin{align*}
F_{k, 5}^{*} (\pi / 2 + \alpha_5) &= e^{i k (\pi / 2 + \alpha_5) / 2} E_{k, 5}^{*}(\rho_{5, 2})\\
 &= \frac{e^{i k (\pi / 2 + \alpha_5) / 2}}{5^{k / 2} + 1}(5^{k / 2} + (2 + i)^k) E_k(i).
\end{align*}
Thus,
\begin{align*}
e^{i k (\pi / 2 + \alpha_5) / 2}(5^{k / 2} + (2 + i)^k) &= 5^{k / 2} e^{i k (\pi / 2 + \alpha_5) / 2}(1 + e^{- i k (\pi / 2 + \alpha_5)})\\
&= 2 \cdot 5^{k / 2} \cos (k (\pi / 2 + \alpha_5) / 2).
\end{align*}\quad

\item[(3)]
Let $k \geqslant 8$ be an even integer divisible by $4$.

First, we consider the case $N = 1$. Then, we can write:
\begin{equation*}
F_{k, 5}^{*} (\pi) = F_{k, 5, 2}^{*}(\pi / 2) = 2 \cos(k \pi / 4) + R_{5, \pi}^{*}.
\end{equation*}
where $R_{5, \pi}^{*}$ denotes the remaining terms. We will consider the following cases: $N = 2$, $5$, $10$, $13$, $17$, and $N \geqslant 25$. We have
\begin{allowdisplaybreaks}
\begin{align*}
&\text{When $N = 2$,}&
v_k(1, 1, \pi / 2) &\leqslant (2 / 3)^{k/2}.\\
&\text{When $N = 5$,}&
v_k(1, 2, \pi / 2) &\leqslant (1 / 21)^{k/2},
&v_k(2, 1, \pi / 2) &\leqslant (1 / 3)^{k}.\\
&\text{When $N = 10$,}&
v_k(1, 3, \pi / 2) &\leqslant (2 / 23)^{k/2},
&v_k(3, 1, \pi / 2) &\leqslant (2 / 7)^{k/2}.\\
&\text{When $N = 13$,}&
v_k(2, 3, \pi / 2) &\leqslant (1 / 7)^{k},
&v_k(3, 2, \pi / 2) &\leqslant (1 / 29)^{k/2}.\\
&\text{When $N = 17$,}&
v_k(1, 4, \pi / 2) &\leqslant (1 / 21)^{k/2},
&v_k(4, 1, \pi / 2) &\leqslant (1 / 3)^{2 k}.\\
&\text{When $N \geqslant 25$,}&
c^2 + 5 d^2 &\geqslant N,
\end{align*}
\end{allowdisplaybreaks}
and the number of terms with $c^2 + d^2 = N$ is not more than $(96 / 25) N^{1/2}$ for $N \geqslant 25$. Then
\begin{equation*}
|R_{5, \pi}^{*}|_{N \geqslant 25}
\leqslant \frac{1536}{25(k-3)} \left(\frac{1}{6}\right)^{(k - 3) / 2}.
\end{equation*}
Furthermore,
\begin{align*}
|R_{5, \pi}^{*}|
 &\leqslant
 4 \left(\frac{2}{3}\right)^{k/2}
 + \cdots + 4 \left(\frac{1}{3}\right)^{2 k}
 + \frac{1536}{25(k-3)} \left(\frac{1}{6}\right)^{(k - 3) / 2},\\
 &\leqslant 0.95701... \quad (k \geqslant 8)
\end{align*}\quad

\item[(4)]
Let $k \geqslant 10$ be an even integer.

First, we consider the case of $N = 1$. Since $5 \pi / 6 < \pi / 2 + \alpha_5$, we can write:
\begin{equation*}
F_{k, 5}^{*} (\theta) = F_{k, 5, 1}^{*}(\theta) = 2 \cos(k \theta / 2) + R_{5, 5 \pi / 6}^{*} \quad \text{for} \: \theta \in [\pi / 2, 5 \pi / 6],
\end{equation*}
where $R_{5, 5 \pi / 6}^{*}$ denotes the remaining terms. We will consider the following cases: $N = 2, 5, 10$, and $N \geqslant 13$. Considering $- \sqrt{3} / 2 \leqslant \cos\theta \leqslant 0$, we have
\begin{allowdisplaybreaks}
\begin{align*}
&\text{When $N = 2$,}&
v_k(1, 1, \theta) &\leqslant 1 / \left(6 - \sqrt{15}\right)^{k/2}, \quad
&v_k(1, - 1, \theta) &\leqslant (1 / 6)^{k/2}.\\
&\text{When $N = 5$,}&
v_k(1, 2, \theta) &\leqslant 1 / \left(21 - 2 \sqrt{15}\right)^{k/2},
&v_k(1, - 2, \theta) &\leqslant (1 / 21)^{k/2},\\
&&v_k(2, 1, \theta) &\leqslant 1 / \left(9 - 2 \sqrt{15}\right)^{k/2},
&v_k(2, - 1, \theta) &\leqslant (1 / 3)^{k}.\\
&\text{When $N = 10$,}&
v_k(1, 3, \theta) &\leqslant 1 / \left(46 - 3 \sqrt{15}\right)^{k/2},
&v_k(1, - 3, \theta) &\leqslant (1 / 46)^{k/2},\\
&&v_k(3, 1, \theta) &\leqslant 1 / \left(14 - 3 \sqrt{15}\right)^{k/2},
&v_k(3, - 1, \theta) &\leqslant (1 / 14)^{k/2}.\\
&\text{When $N \geqslant 13$,}&
|c e^{i \theta / 2} \pm \sqrt{5} d &e^{-i \theta / 2}|^2 \geqslant N / 5,
\end{align*}
\end{allowdisplaybreaks}
and the number of terms with $c^2 + d^2 = N$ is not more than $(21 / 5) N^{1/2}$ for $N \geqslant 13$. Then
\begin{equation*}
|R_{5, 5 \pi / 6}^{*}|_{N \geqslant 13}
\leqslant \frac{1008 \sqrt{3}}{5(k-3)} \left(\frac{5}{12}\right)^{(k - 3) / 2}.
\end{equation*}
Furthermore,
\begin{align*}
|R_{5, \pi / 2}^{*}|
 &\leqslant
 2 \left(\frac{1}{9 - 2 \sqrt{15}}\right)^{k/2}
 + \cdots + 2 \left(\frac{1}{46}\right)^{k/2}
 + \frac{1008 \sqrt{3}}{5(k-3)} \left(\frac{5}{12}\right)^{(k - 3) / 2},\\
 &\leqslant 1.34372... \quad (k \geqslant 10)
\end{align*}
\end{trivlist}
\begin{flushright}$\square$\end{flushright}\quad

In the proofs of the remaining lemmas, we will use the algorithm of subsection \ref{subsec-alg}. Furthermore, we have $X_1 = v_k(2, 1, \theta_1)^{- 2 / k} \geqslant 1 + 4 t (\pi / k)$ in the proof of lemmas \ref{lem-f5s3} and \ref{lem-f5s3}, and we have $X_1 = (1 / 4) \; v_k(1, -1, \theta_2)^{- 2 / k} \geqslant 1 + t (\pi / k)$ in the proof of lemmas \ref{lem-f5s4} and \ref{lem-f5s4}.\\

\subsection{Proof of Lemma \ref{lem-f5s2}}\quad\vspace{-0.1in} \\

\begin{trivlist}
\item[(3)]
Let $k \geqslant 12$ and $x = \pi / (2 k)$, then $0 \leqslant x \leqslant \pi / 24$, and then $1 - \cos x \geqslant (32/33) x^2$. Thus, we have
\begin{gather*}
\frac{1}{4} |e^{i \theta / 2} - \sqrt{5} e^{-i \theta / 2}|^2
 \geqslant \frac{1}{4} (6 - 2 \sqrt{5} \cos(\alpha_5 + x))
 \geqslant 1 + \frac{16}{11} x^2.\\
\frac{1}{2^k} |e^{i \theta / 2} + \sqrt{5} e^{-i \theta / 2}|^k \geqslant 1 + \frac{96}{11} x^2 \; (k \geqslant 12).\\
2^k \cdot 2 v_k(1, 1, \theta) \leqslant 2 - \frac{288 \pi^2}{\pi^2 + 66} \frac{1}{k^2}.
\end{gather*}
In inequality(\ref{r*52bound0}), replace $2$ with the bound $2 - \frac{288 \pi^2}{\pi^2 + 66} \frac{1}{k^2}$. Then
\begin{equation*}
|R_{5, 2}^{*}|
 \leqslant 2 - \frac{288 \pi^2}{\pi^2 + 66} \frac{1}{k^2}
 + 2 \left(\frac{2}{3}\right)^{k/2}
 + \cdots + 2 \left(\frac{1}{129}\right)^{k/2}
 + \frac{2112 \sqrt{33}}{k-3} \left(\frac{8}{33}\right)^{k/2}.
\end{equation*}
Furthermore, $(2 / 3)^{k / 2}$ is more rapidly decreasing in $k$ than $1 / k^2$, and for $k \geqslant 12$, we have
\begin{equation*}
|R_{5, 2}^{*}| \leqslant 1.9821...
\end{equation*}\quad

\noindent
We have $c_0 = c_{0, 1} \leqslant \cos({c_0}' \pi)$.

\item[(1)]
Let $(t, b, s, k_0) = (1 / 6, 59 / 10, 2, 12)$, then we can define $a_1 := 25$. Furthermore, let\\ $({c_1}', c_{0, 1}, a_{1, 1}, u_1) = (1, 1 / 2, 25, 2 / 3)$, then we can define $c_1 := 2$, $a_{2, 1} := 111$, and then we have $Y_1 = 0.38101... > 0$.

\item[(2)]
Let $(t, b, s, k_0) = (9 / 40, 5, 2, 58)$, then $a_1 := 1 / 100$. Furthermore, let\\ $({c_1}', c_{0, 1}, a_{1, 1}, u_1) = (1, 27 / 100, 1 / 100, 9 / 10)$, then $c_1 := 100 / 27$, $a_{2, 1} := 1$, and then $Y_1 = 0.14683... > 0$.\\

\item[(4)]
Let $(t, b, s, k_0) = (1 / 2, 11 / 5, 3 / 2, 22)$, then $a_1 := 5 / 2$. Furthermore, let\\ $({c_1}', c_{0, 1}, a_{1, 1}, u_1) = (1, 1 / 2, 5 / 2, 1 / 2)$, then $c_1 := 2$, $a_{2, 1} := 11$, and then $Y_1 = 0.078258... > 0$.

\item[(5)]
Let $(t, b, s, k_0) = (1 / 5, 41 /20, 3 / 2, 46)$, then $a_1 := 1 / 2$. Furthermore, let\\ $({c_1}', c_{0, 1}, a_{1, 1}, u_1) = (1, 37 / 50, 1 / 2, 1 / 5)$, then $c_1 := 50 / 37$, $a_{2, 1} := 1$, and then $Y_1 = 0.021834... > 0$.
\end{trivlist}
\begin{flushright}$\square$\end{flushright}\quad

\subsection{Proofs of Lemma \ref{lem-f5s3} and Lemma \ref{lem-f5s4}}\quad\vspace{-0.1in} \\

\noindent
{\it Proof of Lemma \ref{lem-f5s3}}\quad
We have $c_0 = c_{0, 1} \leqslant \cos({c_0}' \pi) = - \cos((x / 180) \pi - (t / 2) \pi)$. Furthermore, when $13 \pi / 20 \leqslant (x / 180) \pi < \alpha_{5, k} < (y / 180) \pi \leqslant 7 \pi / 10$, we have $3 \pi / 2 < \pi + (x / 180) \pi + d_1 (t / 2) \pi < {\alpha_{5, k}}' < \pi + (y / 180) \pi + (t / 2) \pi < 2 \pi$. Thus, we can define ${c_1}'$ such that ${c_1}' \geqslant \cos(\pi + (y / 180) \pi + (t / 2) \pi)$.
\begin{trivlist}
\item[(2)]
Let $(b, s) = (21 / 5, 2)$, then $a_1 := 27 / 2$. Furthermore, let\\ $({c_1}', c_{0, 1}, a_{1, 1}, u_1) = (643 / 1000, 179 / 500, 27 / 2, 2 / 5)$, then $c_1 := 643 / 358$, $a_{2, 1} := 69$, and then $Y_1 = 0.019768... > 0$.\vspace{-0.1in}\\

For the following items, we have $(b, s) = (41 / 10, 2)$, then we can define $a_1 := 1 / 100$.
\item[(1)]
Let $({c_1}', c_{0, 1}, a_{1, 1}, u_1) = (1521 / 2000, 3 / 10, 1 / 100, 3 / 5)$, then $c_1 := 507 / 200$, $a_{2, 1} := 1 / 10$, and then $Y_1 = 0.0069363... > 0$.

\item[(3)]
Let $({c_1}', c_{0, 1}, a_{1, 1}, u_1) = (63 / 100, 339 / 1000, 1 / 100, 2 / 5)$, then $c_1 := 210 / 113$, $a_{2, 1} := 1 / 10$, and then $Y_1 = 0.0094197... > 0$.

\item[(4)]
Let $({c_1}', c_{0, 1}, a_{1, 1}, u_1) = (2939 / 5000, 3567 / 10000, 1 / 100, 8 / 25)$, then $c_1 := 5878 / 3567$, $a_{2, 1} := 1 / 20$, and then $Y_1 = 0.0012860... > 0$.

\item[(5)]
Let $({c_1}', c_{0, 1}, a_{1, 1}, u_1) = (5607 / 10000, 3697 / 10000, 1 / 100, 4 / 15)$, then $c_1 := 5607 / 3697$, $a_{2, 1} := 1 / 20$, and then $Y_1 = 0.00069598... > 0$.

\item[(6)]
Let $({c_1}', c_{0, 1}, a_{1, 1}, u_1) = (2731 / 5000, 1877 / 5000, 1 / 100, 6 / 25)$, then $c_1 := 2731 / 1877$, $a_{2, 1} := 1 / 25$, and then $Y_1 = 0.0011623... > 0$.

\item[(7)]
Let $({c_1}', c_{0, 1}, a_{1, 1}, u_1) = (1323 / 2500, 387 / 1000, 1 / 100, 1 / 5)$, then $c_1 := 294 / 215$, $a_{2, 1} := 1 / 25$, and then $Y_1 = 0.00046395... > 0$.

\item[(8)]
Let $({c_1}', c_{0, 1}, a_{1, 1}, u_1) = (5199 / 10000, 3923 / 10000, 1 / 100, 9 / 50)$, then $c_1 := 5199 / 3923$, $a_{2, 1} := 1 / 25$, and then $Y_1 = 0.00067225... > 0$.

\item[(9)]
Let $({c_1}', c_{0, 1}, a_{1, 1}, u_1) = (5113 / 10000, 3981 / 10000, 1 / 100, 4 / 25)$, then $c_1 := 5113 / 3981$, $a_{2, 1} := 1 / 25$, and then $Y_1 = 0.032262... > 0$.

\item[(10)]
Let $({c_1}', c_{0, 1}, a_{1, 1}, u_1) = (51 / 100, 397 / 1000, 1 / 100, 4 / 25)$, then $c_1 := 510 / 397$, $a_{2, 1} := 1 / 25$, and then $Y_1 = 0.032358... > 0$. \qquad $\square$
\end{trivlist}\quad

\noindent
{\it Proof of Lemma \ref{lem-f5s4}}\quad
We have $\cos({c_0}' \pi) = \cos((y / 180) \pi - \pi / 2 + (t / 2) \pi)$. Furthermore, when $19 \pi / 30 \leqslant (x / 180) \pi < \alpha_{5, k} < (y / 180) \pi \leqslant 29 \pi / 45$, we have $2 \pi / 15 \leqslant (x / 180) \pi - \pi / 2 < \beta_{5, k} < (y / 180) \pi - \pi / 2 \leqslant 13 \pi / 90$ and $\pi < \pi + (x / 180) \pi - \pi / 2 - (t / 2) \pi < {\beta_{5, k}}' < \pi + (y / 180) \pi - \pi / 2 - d_2 (t / 2) \pi < 3 \pi / 2$. Thus, we can define ${c_1}'$ such that ${c_1}' \geqslant - \cos(\pi + (x / 180) \pi - \pi / 2 - (t / 2) \pi)$.

For each item, we have $(b, s) = (41 / 20, 3 / 2)$, then we can define $a_1 := 1 / 100$.
\begin{trivlist}
\item[(1)]
Let $({c_1}', c_{0, 1}, a_{1, 1}, u_1) = (493 / 500, 96 / 125, 1 / 100, 4 / 25)$, then $c_1 := 493 / 384$, $a_{2, 1} := 1 / 50$, and then $Y_1 = 0.0016658... > 0$.

\item[(2)]
Let $({c_1}', c_{0, 1}, a_{1, 1}, u_1) = (4839 / 5000, 2007 / 2500, 1 / 100, 3 / 25)$, then $c_1 := 1613 / 1338$, $a_{2, 1} := 1 / 50$, and then $Y_1 = 0.0024495... > 0$.

\item[(3)]
Let $({c_1}', c_{0, 1}, a_{1, 1}, u_1) = (24 / 25, 81 / 100, 1 / 100, 1 / 10)$, then $c_1 := 32 / 27$, $a_{2, 1} := 1 / 50$, and then $Y_1 = 0.0051974... > 0$. \qquad $\square$
\end{trivlist}\quad

\section{$\Gamma_0^{*}(7)$ (For Conjecture \ref{conj-g0s7})}
As for the case of $\Gamma_0^{*}(5)$, to prove Conjecture \ref{conj-g0s7} is difficult. Again, the difficulties arise from the argument $Arg(\rho_{7, 2})$.\\

We will prove the lemmas of this section in the final part of this section.

\subsection{$E_{k, 7}^{*}$ of low weights}
First, we have the following bounds.

\begin{lemma}
Let $k \geqslant 4$ be an even integer.
\def\labelenumi{(\arabic{enumi})}
\begin{enumerate}

\item
\begin{equation*}
F_{k, 7}^{*} (\pi / 2)
\begin{cases}
> 0 & k \equiv 0 \pmod{8}\\
< 0 & k \equiv 4 \pmod{8}\\
= 0 & k \equiv 2 \pmod{4}
\end{cases}.
\end{equation*}\quad

\item
\begin{equation*}
F_{4, 7}^{*} (5 \pi / 6) > 0.
\end{equation*}\quad

\item
For an even integer $k \geqslant 6$, we can write
\begin{equation*}
F_{k, 7}^{*} (\theta) = 2 \cos(k \theta / 2) + R_{7, 2 \pi / 3}^{*}, \ \text{and have} \ |R_{7, 2 \pi / 3}^{*}| < 2 \quad \text{for} \ \theta \in [\pi / 2, 2 \pi / 3].
\end{equation*}\quad

\item
For an even integer $k \geqslant 6$, we can write
\begin{equation*}
F_{k, 7}^{*} (\theta) = 2 \cos(k \theta / 2) + R_{7, \pi}^{*}, \ \text{and have} \ |R_{7, \pi}^{*}| < 2 \quad \text{for} \ \theta \in [\pi, 7 \pi / 6].
\end{equation*}\quad

\item
For an even integer $k \geqslant 8$, we can write
\begin{equation*}
F_{k, 7}^{*} (\theta) = 2 \cos(k \theta / 2) + R_{7, 5 \pi / 6}^{*}, \ \text{and have} \ |R_{7, 5 \pi / 6}^{*}| < 2 \quad \text{for} \ \theta \in [\pi / 2, 5 \pi / 6].
\end{equation*}
\end{enumerate}
\def\labelenumi{\arabic{enumi}.}\label{lem-f7s1}
\end{lemma}

We will present a proof of the above lemma in Subsection \ref{subsec-pf71}. Now, we have the following proposition:

\begin{proposition}
The location of the zeros of the Eisenstein series $E_{k, 7}^{*}$ in $\mathbb{F}^{*}(7)$ for $k = 4, 6,$ and $12$ are follows $:$
\begin{center}
\begin{tabular}{cccccc}
$k$ & $v_{\infty}$ & $v_{i / \sqrt{7}}$ & $v_{\rho_{7, 1}}$ & $v_{\rho_{7, 2}}$ & $V_7^{*}$\\
\hline
$4$ & $0$ & $0$ & $0$ & $1$ & $1$\\
$6$ & $0$ & $1$ & $1$ & $0$ & $1$\\
$12$ & $0$ & $0$ & $0$ & $0$ & $4$\\
\hline
\end{tabular}
\end{center}
where $V_7^{*}$ denote the number of simple zeros of the Eisenstein series $E_{k, 7}^{*}$ on $A_{7, 1}^{*} \cup A_{7, 2}^{*}$.\label{prop-ek7s-lw}
\end{proposition}

\begin{proof}\quad
\begin{trivlist}
\item[($k = 4$)]
We have $F_{4, 7}^{*}(\pi / 2) < 0$ by Lemma \ref{lem-f7s1} (1), and we have $F_{4, 7}^{*}(5 \pi / 6) > 0$ by Lemma \ref{lem-f7s1} (2). Thus, $E_{4, 7}^{*}$ has at least one zero on $A_{7, 1}^{*}$. In addition, by the previous subsubsections, we have $v_{\rho_{7, 2}} (E_{4, 7}^{*}) \geqslant 1$. Furthermore, by the valence formula for $\Gamma_0^{*}(7)$ (Proposition \ref{prop-vf-g0s7}), $E_{4, 7}^{*}$ has no other zeros.

\item[($k = 6$)]
We have $F_{6, 7}^{*}(2 \pi / 3) > 0$ by Lemma \ref{lem-f7s1} (3), and we have $F_{6, 7}^{*}(\pi) < 0$ by Lemma \ref{lem-f7s1} (4), thus $E_{6, 7}^{*}$ has at least one zero on $A_{7, 1}^{*} \cup A_{7, 2}^{*}$. By previous subsubsection, we have $v_{i / \sqrt{7}} (E_{6, 7}^{*}) \geqslant 1$ and $v_{\rho_{7, 1}} (E_{6, 7}^{*}) \geqslant 1$. Furthermore, by the valence formula for $\Gamma_0^{*}(7)$, $E_{6, 7}^{*}$ has no other zeros.

\item[($k = 12$)]
We have $F_{12, 7}^{*}(\pi / 2) < 0$, $F_{12, 7}^{*}(2 \pi / 3) > 0$, $F_{12, 7}^{*}(5 \pi / 6) < 0$ by Lemma \ref{lem-f7s1} (5), and we have $F_{12, 7}^{*}(\pi) > 0$, $F_{12, 7}^{*}(7 \pi / 6) < 0$ by Lemma \ref{lem-f7s1} (4). Thus, $E_{12, 7}^{*}$ has at least four zeros on $A_{7, 1}^{*} \cup A_{7, 2}^{*}$. Furthermore, by the valence formula for $\Gamma_0^{*}(7)$, $E_{12, 7}^{*}$ has no other zeros.
\end{trivlist}
\end{proof}\quad

\subsection{All but at most $2$ zeros}\label{subsec-g0s7-ab2}
We have the following lemma:

\begin{lemma}We have the following bounds$:$\\
``We have $|R_{7, 1}^{*}| < 2 \cos({c_0}' \pi)$ for $\theta_1 \in [\pi / 2, \pi / 2 + \alpha_7 - t \pi / k]$''\\
\quad$(1)$ For $k \geqslant 10$, $({c_0}', t) = (1 / 3, 1 / 3)$.\\
\quad$(2)$ For $k \geqslant 80$, $({c_0}', t) = (41 / 100, 8 / 25)$.\\
\quad$(3)$ For $k \geqslant 22$, $({c_0}', t) = (13 / 36, 1 / 3)$.\\
``We have $|R_{7, 2}^{*}| < 2 \cos({c_0}' \pi)$ for $\theta_2 \in [\alpha_7 - \pi / 6 + t \pi / k, \pi / 2]$''\\
\quad$(4)$ For $k \geqslant 8$, $({c_0}', t) = (1 / 6, 1 / 2)$.\\
\quad$(5)$ For $k = 26$, $k \geqslant 44$, $({c_0}', t) = (1 / 3, 2 / 3)$.\\
\quad$(6)$ For $k \geqslant 70$, $({c_0}', t) = (1 / 4, 1 / 2)$.\\
\quad$(7)$ For $k \geqslant 200$, $({c_0}', t) = (5 / 18, 1 / 2)$.\label{lem-f7s2}
\end{lemma}

\subsubsection{The case $6 \mid k$}
In agreement with expectations we have at most $k / 3$ (resp. $(k - 3) / 3$) zeros on the arc $A_{7, 1}^{*} \cup A_{7, 2}^{*}$ when $k \equiv 0 \pmod{12}$ (resp. $k \equiv 6 \pmod{12}$) between adjacent integer points ({\it i.e.} $\cos\left( k \theta / 2 \right) = \pm 1$). By Lemma \ref{lem-f7s2} (1) and (4), we can prove $|R_{7, 1}^{*}| < 2$ or $|R_{7, 2}^{*}| < 2$ at all of the integer points that satisfy $\theta_1 \in [\pi / 2, \pi / 2 + \alpha_7 - \pi / (3 k)]$ or $\theta_2 \in [\alpha_7 - \pi / 6 + \pi / (2 k), \pi / 2]$, respectively. Then, we have all but at most $2$ zeros on $A_7^{*}$.

\subsubsection{The case $k \equiv 2 \pmod{6}$}
We have at most $(k - 2) / 3 - 1$ (resp. $(k - 2) / 3$) zeros on the arc $A_{7, 1}^{*} \cup A_{7, 2}^{*}$ when $k \equiv 2 \pmod{12}$ (resp. $k \equiv 8 \pmod{12}$) between adjacent integer points.\\

\paragraph{\it The case $0 < \alpha_{5, k} < 2 \pi / 3$}
When $k \equiv 2 \pmod{12}$, between the first integer point and the last one for $A_{7, 1}^{*}$, we have $(k \alpha_7 / 2 - \pi / 2 - \alpha_{7, k}) / \pi + 1$ integer points, then we expect $(k \alpha_7 / 2 - \pi / 2 - \alpha_{7, k}) / \pi$ zeros on the arc $A_{7, 1}^{*}$. Similarly, we expect $(k (2 \pi / 3 - \alpha_7) / 2 - \pi / 2 - (\pi - \beta_{7, k})) / \pi$ zeros on the arc $A_{7, 2}^{*}$. Furthermore, we have $\alpha_{7, k} + (\pi - \beta_{7, k}) = 2 \pi / 3$. Thus, we expect $(k - 2) / 3 - 1$ zeros on the arc $A_{7, 1}^{*} \cup A_{7, 2}^{*}$. Also, when $k \equiv 8 \pmod{12}$, we expect $(k - 2) / 3$ zeros. By Lemma \ref{lem-f7s2} (1) and (4), we prove all but at most one of the zeros are on $A_7^{*}$.\\

\paragraph{\it The case $2 \pi / 3 < \alpha_{5, k} < \pi$}
We expect $(k - 2) / 3 - 2$ (resp. $(k - 2) / 3 - 1$) zeros between adjacent integer points for $A_{7, 1}^{*} \cup A_{7, 2}^{*}$ when $k \equiv 2 \pmod{12}$ (resp. $k \equiv 8 \pmod{12}$). By Lemma \ref{lem-f7s2} (1) and (4), we prove all but at most one of the zeros are on $A_7^{*}$.\\

\subsubsection{The case $k \equiv 4 \pmod{6}$}
We have at most $(k - 1) / 3$ (resp. $(k - 1) / 3 - 1$) zeros on the arc $A_{7, 1}^{*} \cup A_{7, 2}^{*}$ when $k \equiv 4 \pmod{12}$ (resp. $k \equiv 10 \pmod{12}$) between adjacent integer points.\\

\paragraph{\it The case $0 < \alpha_{5, k} < \pi / 3$}
When $k \equiv 4 \pmod{12}$, we expect $(k - 1) / 3$ zeros on the arc $A_{7, 1}^{*} \cup A_{7, 2}^{*}$ between adjacent integer points. Also, when $k \equiv 10 \pmod{12}$, we expect $(k - 1) / 3 - 1$ zeros. By Lemma \ref{lem-f7s2} (1) and (4), we prove all but at most two of the zeros are on $A_7^{*}$.\\

\paragraph{The case $\pi / 3 < \alpha_{5, k} < \pi$}
We expect $(k - 1) / 3 - 1$ (resp. $(k - 1) / 3 - 2$) zeros between adjacent integer points for $A_{7, 1}^{*} \cup A_{7, 2}^{*}$ when $k \equiv 4 \pmod{12}$ (resp. $k \equiv 10 \pmod{12}$). By Lemma \ref{lem-f7s2} (1) and (4), we prove all but at most one of the zeros are on $A_7^{*}$.\\

In conclusion, we have the following proposition:
\begin{proposition}
Let $k \geqslant 4$ be an even integer. All but at most two of the zeros of $E_{k, 7}^{*}(z)$ in $\mathbb{F}^{*}(7)$ are on the arc $A_7^{*}$. \label{prop-g0s7-ab2}
\end{proposition}\quad

\subsection{The case $6 \mid k$}\label{subsec-g0s7-60}
We can prove $|R_{7, 1}^{*}| < 2$ at all of the integer points for $A_{7, 1}^{*}$ when $\pi / 6 < \alpha_{7, k}$, and we can prove this bound at all of the integer points less than the last one when $\alpha_{7, k} < \pi / 6$. On the other hand, for $A_{7, 2}^{*}$, we can prove $|R_{7, 2}^{*}| < 2$ at all of the integer points when $\alpha_{7, k} < \pi / 4$, and we can prove this bound at all of the integer points greater than the first one when $\pi / 4 < \alpha_{7, k}$.

We can write
\begin{gather*}
F_{k, 7, 1}^{*}(\theta_1) = 2 \cos\left( k \theta_1 / 2 \right) + 2 Re(2 e^{- i \theta_1 / 2} + \sqrt{7} e^{i \theta_1 / 2})^{- k} + 2 Re(3 e^{- i \theta_1 / 2} + \sqrt{7} e^{i \theta_1 / 2})^{- k} + {R_{7, 1}^{*}}',\\
F_{k, 7, 2}^{*}(\theta_2) = 2 \cos\left( k \theta_2 / 2 \right) + 2^k \cdot 2 Re(e^{- i \theta_2 / 2} - \sqrt{7} e^{i \theta_2 / 2})^{- k} + 2^k \cdot 2 Re(3 e^{- i \theta_2 / 2} - \sqrt{7} e^{i \theta_2 / 2})^{- k} + {R_{7, 2}^{*}}'.
\end{gather*}\vspace{-0.1in}

When $0 < \alpha_{7, k} < \pi / 8$, the last integer point for $A_{7, 1}^{*}$ is in the interval $[\pi / 2 + \alpha_7 - \pi / (8 k), \pi / 2 + \alpha_7]$. We have $|{R_{7, 1}^{*}}'| < 2$ and $Sign\{\cos(k \theta_1 / 2)\} = Sign\{Re(2 e^{- i \theta_1 / 2} + \sqrt{7} e^{i \theta_1 / 2})^{- k}\} = Sign\{Re(3 e^{- i \theta_1 / 2} + \sqrt{7} e^{i \theta_1 / 2})^{- k}\}$ for $\theta_1 \in [\pi / 2 + \alpha_7 - \pi / (8 k), \pi / 2 + \alpha_7]$.

When $\pi / 8 < \alpha_{7, k} < \pi / 6$, we can use Lemma \ref{lem-f7s2} (1). Instead of the last integer point for $A_{7, 1}^{*}$, we consider the point such that $k \theta_1 / 2 = k (\pi / 2 + \alpha_7) / 2 - \alpha_{7, k} - \pi / 3$. Then, we have $|R_{7, 1}^{*}| < 1$ and $2 \cos(k \theta_1 / 2) = \pm 1$ at this point.

When $\pi / 4 < \alpha_{7, k} < 5 \pi / 6$, we can use Lemma \ref{lem-f7s2} (4). We consider the point such that $k \theta_2 / 2 = k (\alpha_7 - \pi / 6) / 2 + (\pi - \beta_{7, k}) + \pi / 6$. Then, we have $|R_{7, 2}^{*}| < \sqrt{3}$ and $2 \cos(k \theta_2 / 2) = \pm \sqrt{3}$ at this point.

When $5 \pi / 6 < \alpha_{7, k} < \pi$, the first integer point for $A_{7, 2}^{*}$ is in the interval $[\alpha_7 - \pi / 6 + \pi / (6 k), \pi / 2]$. We have$|{R_{7, 2}^{*}}'| < 2$ and $Sign\{\cos(k \theta_2 / 2)\} = Sign\{Re(e^{- i \theta_2 / 2} - \sqrt{7} e^{i \theta_2 / 2})^{- k}\} = Sign\{Re(3 e^{- i \theta_2 / 2} - \sqrt{7} e^{i \theta_2 / 2})^{- k}\}$ for $\theta_2 \in [\alpha_7 - \pi / 6 + \pi / (6 k), \pi / 2]$.

Thus, we have the following proposition:
\begin{proposition}
Let $k \geqslant 4$ be an integer which satisfies $6 \mid k$. All of the zeros of $E_{k, 7}^{*}(z)$ in $\mathbb{F}^{*}(7)$ are on the arc $A_7^{*}$. \label{prop-g0s7-60}
\end{proposition}\quad

\subsection{The case $k \equiv 2 \pmod{6}$}\label{subsec-g0s7-61}

\subsubsection{The case $0 < \alpha_{5, k} < 2 \pi / 3$}
At the integer points for $A_{7, 1}^{*}$, we can use Lemma \ref{lem-f7s2} (1). We can prove $|R_{7, 1}^{*}| < 2$ at all of the integer points when $\pi / 6 < \alpha_{7, k}$, and we can prove this bound at all of the integer points less than the last integer when $\alpha_{7, k} < \pi / 6$. When $\alpha_{7, k} < \pi / 6$, we consider the point $k \theta_1 / 2 = k (\pi / 2 + \alpha_7) / 2 - \alpha_{7, k} - \pi / 3$.

On the other hand, at the integer points for $A_{7, 2}^{*}$, we can use Lemma \ref{lem-f7s2} (4) and (5). We can prove $|R_{7, 2}^{*}| < 2$ at all of the integer points when $\alpha_{7, k} < 5 \pi / 12$ ({\it i.e.} $\beta_{7, k} < \pi / 4$), and we can prove this bound at all of the integer points greater than the first integer when $5 \pi / 12 < \alpha_{7, k}$. Furthermore, we consider the point $k \theta_2 / 2 = k (\alpha_7 - \pi / 6) / 2 + (\pi - \beta_{7, k}) + \pi / 6$ when $5 \pi / 12 < \alpha_{7, k} < 7 \pi / 12$ using Lemma \ref{lem-f7s2} (4), and we consider the point $k \theta_2 / 2 = k (\alpha_7 - \pi / 6) / 2 + (\pi - \beta_{7, k}) + \pi / 3$ when $7 \pi / 12 < \alpha_{7, k} < 2 \pi / 3$ using Lemma \ref{lem-f7s2} (5).\\

\subsubsection{The case $2 \pi / 3 < \alpha_{5, k} < \pi$}
We have the following lemmas:

\begin{lemma}Let $k \equiv 2 \pmod{6}$. We have the following bounds$:$\\
``When $(x / 180) \pi < \alpha_{7, k} < (y / 180) \pi$, we have $|R_{7, 1}^{*}| < 2 \cos({c_0}' \pi)$ for $k \geqslant k_0$ and $\theta_1 = \pi / 2 + \alpha_7 - t \pi / k$.''
\def\labelenumi{(\arabic{enumi})}
\begin{enumerate}
\item $(x, y) = (131.5, 135)$, $({c_0}', k_0, t) = (71 / 180, 212, 1 / 4)$.
\item $(x, y) = (130.1, 131.5)$, $({c_0}', k_0, t) = (2743 / 7200, 740, 83 / 400)$.
\item $(x, y) = (129.5, 130.1)$, $({c_0}', k_0, t) = (53 / 144, 872, 7 / 40)$.
\item $(x, y) = (129.18, 129.5)$, $({c_0}', k_0, t) = (541 / 1500, 1000, 47 / 300)$.
\item $(x, y) = (129, 129.18)$, $({c_0}', k_0, t) = (1063 / 3000, 1000, 71 / 500)$.
\item $(x, y) = (128.86, 129)$, $({c_0}', k_0, t) = (2519 / 7200, 1000, 263 / 2000)$.
\item $(x, y) = (128.77, 128.86)$, $({c_0}', k_0, t) = (6221 / 18000, 1000, 61 / 500)$.
\item $(x, y) = (128.71, 128.77)$, $({c_0}', k_0, t) = (30793 / 90000, 1000, 143 / 1250)$.
\item $(x, y) = (128.68, 128.71)$, $({c_0}', k_0, t) = (6113 / 18000, 1000, 109 / 1000)$.
\end{enumerate}
\def\labelenumi{\arabic{enumi}.}\label{lem-f7s3}
\end{lemma}

\begin{lemma}Let $k \equiv 2 \pmod{6}$. We have the following bounds$:$\\
``When $(x / 180) \pi < \alpha_{7, k} < (y / 180) \pi$, we have $|R_{7, 2}^{*}| < 2 \cos({c_0}' \pi)$ for $k \geqslant k_0$ and $\theta_2 = \alpha_7 - \pi / 6 + t \pi / k$.''
\def\labelenumi{(\arabic{enumi})}
\begin{enumerate}
\item $(x, y) = (120, 126.7)$, $({c_0}', k_0, t) = (971 / 3600, 194, 93 / 200)$.
\item $(x, y) = (126.7, 127.3)$, $({c_0}', k_0, t) = (379 / 1800, 500, 17 / 50)$.
\item $(x, y) = (127.3, 127.63)$, $({c_0}', k_0, t) = (3403 / 18000, 1000, 22 / 75)$.
\item $(x, y) = (127.63, 127.68)$, $({c_0}', k_0, t) = (259 / 1500, 1000, 13 / 50)$.
\end{enumerate}
\def\labelenumi{\arabic{enumi}.}\label{lem-f7s4}
\end{lemma}

We need one more zero between the last integer point for $A_{7, 1}^{*}$ and the first one for $A_{7, 2}^{*}$. For the remaining zero, we consider following cases.\\

\paragraph{(i)} ``The case $3 \pi / 4 < \alpha_{5, k} < \pi$''
\begin{itemize}
\item When $3 \pi / 4 < \alpha_{7, k} < 5 \pi / 6$, we can use Lemma \ref{lem-f7s2} (2).
\item When $5 \pi / 6 < \alpha_{7, k} < \pi$, we can use Lemma \ref{lem-f7s2} (1).
\end{itemize}
For each case, we consider the point such that $k \theta_1 / 2 = k (\pi / 2 + \alpha_7) / 2 - \alpha_{7, k} + \pi - {c_0}' \pi$. We have $\alpha_{7, k} - \pi + {c_0}' \pi > (t / 2) \pi$ and $|R_{7, 1}^{*}| < 2 \cos({c_0}' \pi)$, and we have $2 \cos(k \theta_1 / 2) = \pm 2 \cos({c_0}' \pi)$. Then, we have at least one zero between the second to last integer point for $A_{7, 1}^{*}$ and the point $k \theta_1 / 2$.\\

\paragraph{(ii)} ``The case $3217 \pi / 4500 < \alpha_{7, k} < 3 \pi / 4$''

We can use Lemma \ref{lem-f7s3}. For each case, we consider the point such that $k \theta_1 / 2 = k (\pi / 2 + \alpha_7) / 2 - (t / 2) \pi$. We have $\alpha_{7, k} - \pi + {c_0}' \pi > (t / 2) \pi$ and $|R_{7, 1}^{*}| < 2 \cos({c_0}' \pi)$, and we have $|2 \cos(k \theta_1 / 2)| > 2 \cos({c_0}' \pi)$. Then, we have at least one zero between the second to the last integer point for $A_{7, 1}^{*}$ and the point $k \theta_1 / 2$.\\

\paragraph{(iii)} ``The case $2 \pi / 3 < \alpha_{7, k} < 266 \pi / 375$'' 

Similar to the case (ii), we can use Lemma \ref{lem-f7s4}. We consider the point such that $k \theta_2 / 2 = k (\alpha_7 - \pi / 6) / 2 + (t / 2) \pi$ for each case.\\

In conclusion, we have the following proposition:
\begin{proposition}
Let $k \geqslant 4$ be an integer which satisfies $k \equiv 2 \pmod{6}$, and let $\alpha_{7, k} \in [0, \pi]$ be the angle that satisfies $\alpha_{7, k} \equiv k (\pi / 2 + \alpha_7) / 2 \pmod{\pi}$. If we have $\alpha_{7, k} < 266 \pi / 375$ or $3217 \pi / 4500 < \alpha_{7, k}$, then all of the zeros of $E_{k, 7}^{*}(z)$ in $\mathbb{F}^{*}(7)$ are on the arc $A_7^{*}$. Otherwise, all but at most one zero of $E_{k, 7}^{*}(z)$ in $\mathbb{F}^{*}(7)$ are on $A_7^{*}$\label{prop-g0s7-61}
\end{proposition}\quad

\subsection{The case $k \equiv 4 \pmod{6}$}\label{subsec-g0s7-62}

\subsubsection{The case $0 < \alpha_{5, k} < \pi / 3$}
Similar to the previous subsubsection, at the integer points for $A_{7, 1}^{*}$, we can use Lemma \ref{lem-f7s2} (1). When $\pi / 6 < \alpha_{7, k}$, we can prove at all of the integer points. When $\alpha_{7, k} < \pi / 6$, we consider the point $k \theta_1 / 2 = k (\pi / 2 + \alpha_7) / 2 - \alpha_{7, k} - \pi / 3$.

On the other hand, at the integer points for $A_{7, 2}^{*}$, we can use Lemma \ref{lem-f7s2} (4) and (6). By Lemma \ref{lem-f7s2} (4), we can prove at all of the integer points greater than the first one, and we consider the point $k \theta_2 / 2 = k (\alpha_7 - \pi / 6) / 2 + (\pi - \beta_{7, k}) + \pi / 6$ when $0 < \alpha_{7, k} < \pi / 4$. Further, by Lemma \ref{lem-f7s2} (6), we consider the point $k \theta_2 / 2 = k (\alpha_7 - \pi / 6) / 2 + (\pi - \beta_{7, k}) + \pi / 3$ when $\pi / 4 < \alpha_{7, k} < \pi / 3$.\\

\subsubsection{The case $\pi / 3 < \alpha_{5, k} < \pi$}
We have the following lemmas:
\begin{lemma}Let $k \equiv 4 \pmod{6}$. We have the following bounds$:$\\
``When $(x / 180) \pi < \alpha_{7, k} < (y / 180) \pi$, we have $|R_{7, 1}^{*}| < 2 \cos({c_0}' \pi)$ for $k \geqslant k_0$ and $\theta_1 = \pi / 2 + \alpha_7 - t \pi / k$.''
\def\labelenumi{(\arabic{enumi})}
\begin{enumerate}
\item $(x, y) = (127.6, 135)$, $({c_0}', k_0, t) = (1841 / 4500, 118, 59 / 250)$.
\item $(x, y) = (120, 127.6)$, $({c_0}', k_0, t) = (11 / 24, 100, 1 / 4)$.
\end{enumerate}
\def\labelenumi{\arabic{enumi}.}\label{lem-f7s5}
\end{lemma}

\begin{lemma}Let $k \equiv 4 \pmod{6}$. We have the following bounds$:$\\
``When $(x / 180) \pi < \alpha_{7, k} < (y / 180) \pi$, we have $|R_{7, 2}^{*}| < 2 \cos({c_0}' \pi)$ for $k \geqslant k_0$ and $\theta_2 = \alpha_7 - \pi / 6 + t \pi / k$.''
\def\labelenumi{(\arabic{enumi})}
\begin{enumerate}
\item $(x, y) = (65, 90)$, $({c_0}', k_0, t) = (11 / 30, 10, 2 / 5)$.
\item $(x, y) = (90, 100)$, $({c_0}', k_0, t) = (19 / 45, 28, 2 / 5)$.
\end{enumerate}
\def\labelenumi{\arabic{enumi}.}\label{lem-f7s6}
\end{lemma}

Now, we can write
\begin{gather*}
F_{k, 7, 1}^{*}(\theta_1) = 2 \cos\left( k \theta_1 / 2 \right) + 2 Re(3 e^{- i \theta_1 / 2} + \sqrt{7} e^{i \theta_1 / 2})^{- k} + {R_{7, 1}^{*}}'',\\
F_{k, 7, 2}^{*}(\theta_2) = 2 \cos\left( k \theta_2 / 2 \right) + 2^k \cdot 2 Re(3 e^{- i \theta_2 / 2} - \sqrt{7} e^{i \theta_2 / 2})^{- k} + {R_{7, 2}^{*}}''.
\end{gather*}

\begin{lemma}
When $k \equiv 4 \pmod{6}$ and $73 \pi / 120 < \alpha_{7, k} < 2 \pi / 3$,\\
 \quad we have $Sign\{\cos\left( k \theta_1 / 2 \right)\} = Sign\{Re(3 e^{- i \theta_1 / 2} + \sqrt{7} e^{i \theta_1 / 2})^{- k}\}$ and the following bounds$:$\\
``When $(x / 180) \pi < \alpha_{7, k} < (y / 180) \pi$, we have $|{R_{7, 1}^{*}}''| < |2 \cos\left( k \theta_1 / 2 \right) + 2 Re(3 e^{- i \theta_1 / 2} + \sqrt{7} e^{i \theta_1 / 2})^{- k}|$\\
 \qquad for $k \geqslant k_0$ and $\theta_1 = \pi / 2 + \alpha_7 - t \pi / k$.''
\def\labelenumi{(\arabic{enumi})}
\begin{enumerate}
\item $(x, y) = (111.6, 120)$, $(k_0, t) = (82, 23 / 150)$.
\item $(x, y) = (110.1, 111.6)$, $(k_0, t) = (742, 1 / 10)$.
\item $(x, y) = (109.65, 110.1)$, $(k_0, t) = (1000, 43 / 625)$.
\item $(x, y) = (109.5, 109.65)$, $(k_0, t) = (1000, 21 / 400)$.
\end{enumerate}
\def\labelenumi{\arabic{enumi}.}\label{lem-f7s7}
\end{lemma}

\begin{lemma}
When $k \equiv 4 \pmod{6}$ and $5 \pi / 9 < \alpha_{7, k} < 217 \pi / 360$,\\
 \quad we have $Sign\{\cos\left( k \theta_2 / 2 \right)\} = Sign\{Re(3 e^{- i \theta_2 / 2} - \sqrt{7} e^{i \theta_2 / 2})^{- k}\}$ and the following bounds$:$\\
``When $(x / 180) \pi < \alpha_{7, k} < (y / 180) \pi$, we have $|{R_{7, 2}^{*}}''| < |2 \cos\left( k \theta_2 / 2 \right) + 2^k \cdot 2 Re(3 e^{- i \theta_2 / 2} - \sqrt{7} e^{i \theta_2 / 2})^{- k}|$\\
 \qquad for $k \geqslant k_0$ and $\theta_2 = \alpha_7 - \pi / 6 + t \pi / k$.''
\def\labelenumi{(\arabic{enumi})}
\begin{enumerate}
\item $(x, y) = (100, 106)$, $(k_0, t) = (46, 3/10)$.
\item $(x, y) = (106, 107.7)$, $(k_0, t) = (196, 11 / 50)$.
\item $(x, y) = (107.7, 108.21)$, $(k_0, t) = (1000, 33 / 200)$.
\item $(x, y) = (108.21, 108.42)$, $(k_0, t) = (1000, 2 / 15)$.
\item $(x, y) = (108.42, 108.5)$, $(k_0, t) = (1000, 113 / 1000)$.
\end{enumerate}
\def\labelenumi{\arabic{enumi}.}\label{lem-f7s8}
\end{lemma}

We need one more zero between the last integer point for $A_{7, 1}^{*}$ and the first integer point for $A_{7, 2}^{*}$. For this remaining zero, we consider following cases.\\

\paragraph{(i)} ``The case $3 \pi / 4 < \alpha_{5, k} < \pi$''
\begin{itemize}
\item When $5 \pi / 6 < \alpha_{7, k} < \pi$, we can use Lemma \ref{lem-f7s2} (1).
\item When $29 \pi / 36 < \alpha_{7, k} < 5 \pi / 6$, we can use Lemma \ref{lem-f7s2} (3).
\item When $3 \pi / 4 < \alpha_{7, k} < 5 \pi / 6$, we can use Lemma \ref{lem-f7s2} (2).
\end{itemize}
Similar to the case (i) in the previous subsubsection, we consider the point such that $k \theta_1 / 2 = k (\pi / 2 + \alpha_7) / 2 - \alpha_{7, k} + \pi - {c_0}' \pi$.\\

\paragraph{(ii)} ``The case $\pi / 3 < \alpha_{7, k} < 13 \pi / 36$''

Similar to the case (i), we can use Lemma \ref{lem-f7s2} (7). We consider the point such that $k \theta_2 / 2 = k (\alpha_7 - \pi / 6) / 2 - \beta_{7, k} + 5 \pi / 18$ for each case.\\

\paragraph{(iii)} ``The case $2 \pi / 3 < \alpha_{7, k} < 3 \pi / 4$''

Similar to the case (ii) in the previous subsubsection, we can use Lemma \ref{lem-f7s5}. We consider the point such that $k \theta_1 / 2 = k (\pi / 2 + \alpha_7) / 2 - (t / 2) \pi$.\\

\paragraph{(iv)} ``The case $13 \pi / 36 < \alpha_{7, k} < 5 \pi / 9$'' 

Similar to the case (iii), we can use Lemma \ref{lem-f7s6}. We consider the point such that $k \theta_2 / 2 = k (\alpha_7 - \pi / 6) / 2 + (t / 2) \pi$ for each case.\\

\paragraph{(v)} ``The case $73 \pi / 120 < \alpha_{7, k} < 2 \pi / 3$'' 

We can use Lemma \ref{lem-f7s7}. For each case, we consider the point such that $k \theta_1 / 2 = k (\pi / 2 + \alpha_7) / 2 - (t / 2) \pi$. We have $Sign\{\cos\left( k \theta_1 / 2 \right)\} = Sign\{Re(3 e^{- i \theta_1 / 2} + \sqrt{7} e^{i \theta_1 / 2})^{- k}\}$ and $|{R_{7, 1}^{*}}''| < 2 \cos\left( k \theta_1 / 2 \right) + 2 Re(3 e^{- i \theta_1 / 2} + \sqrt{7} e^{i \theta_1 / 2})^{- k}$. Then, we can show $Sign\{\cos\left( k \theta_1 / 2 \right)\} = Sign\{F_{k, 7, 1}^{*}(\theta_1)\}$ , and then we have at least one zero between the second to last integer point for $A_{7, 1}^{*}$ and the point $k \theta_1 / 2$.\\

\paragraph{(vi)} ``The case $5 \pi / 9 < \alpha_{7, k} < 217 \pi / 360$'' 

Similar to the case (v), we can use Lemma \ref{lem-f7s8}.\\

In conclusion, we have the following proposition:
\begin{proposition}
Let $k \geqslant 4$ be an integer which satisfies $k \equiv 4 \pmod{6}$, and let $\alpha_{7, k} \in [0, \pi]$ be the angle that satisfies $\alpha_{7, k} \equiv k (\pi / 2 + \alpha_7) / 2 \pmod{\pi}$. If we have $\alpha_{7, k} < 217 \pi / 360$ or $73 \pi / 120 < \alpha_{7, k}$, then all of the zeros of $E_{k, 7}^{*}(z)$ in $\mathbb{F}^{*}(7)$ are on the arc $A_7^{*}$. Otherwise, all but at most one zero of $E_{k, 7}^{*}(z)$ in $\mathbb{F}^{*}(7)$ are on $A_7^{*}$\label{prop-g0s7-62}
\end{proposition}\quad

\subsection{The remaining cases ``$k \equiv 2 \pmod{6}$, $266 \pi / 375 < \alpha_{7, k} < 3217 \pi / 4500$'' and ``$k \equiv 4 \pmod{6}$, $217 \pi / 360 < \alpha_{7, k} < 73 \pi / 120$'' } \label{subsec-rest7}
Similar to Subsection \ref{subsec-rest5}, it is difficult to prove for these remaining cases, and these cannot be proved in a similar way. However, when $k$ is large enough, we have the following observation.

\subsubsection{\bfseries The case ``$k \equiv 2 \pmod{6}$ and $266 \pi / 375 < \alpha_{7, k} < 3217 \pi / 4500$''}
Let $t > 0$ be small enough, then we have $\pi / 2 < \alpha_{7, k} - (t / 2) \pi < \pi$, $\pi < 2 \pi / 3 + \alpha_{7, k} + d_{1, 1} (t / 2) \pi < {\alpha_{7, k, 1}}' < 2 \pi / 3 + \alpha_{7, k} + (t / 2) \pi < 3 \pi / 2$, and $2 \pi < 4 \pi / 3 + \alpha_{7, k} - t \pi < {\alpha_{7, k, 2}}' < 4 \pi / 3 + \alpha_{7, k} - d_{1, 2} (t / 2) \pi < 5 \pi / 2$. Thus, we have
\begin{allowdisplaybreaks}
\begin{align*}
&- \cos(\alpha_{7, k} - (t / 2) \pi)
 - \cos(2 \pi / 3 + \alpha_{7, k} + d_{1, 1} (t / 2) \pi) \cdot (1 + 2 \sqrt{3} t (\pi / k))^{- k / 2}\\
 &\qquad \qquad \qquad - \cos(4 \pi / 3 + \alpha_{7, k} - d_{1, 2} (t / 2) \pi) \cdot e^{- (3 \sqrt{3} / 2) \pi t}\\
 &\qquad > |\cos(k \theta_1 / 2)|
 + \left| Re\left\{ \left( 2 e^{i \theta_1 / 2} + \sqrt{7} e^{- i \theta_1 / 2} \right)^{- k} \right\} \right|\\
 &\qquad \qquad \qquad \qquad - \left| Re\left\{ \left( 3 e^{i \theta_1 / 2} + \sqrt{7} e^{- i \theta_1 / 2} \right)^{- k} \right\} \right|\\
 &\qquad > - \cos(\alpha_{7, k} - (t / 2) \pi)
 - \cos(2 \pi / 3 + \alpha_{7, k} + (3 t / 2) \pi) \cdot e^{- \sqrt{3} \pi t}\\
 &\qquad \qquad \qquad \qquad - \cos(4 \pi / 3 + \alpha_{7, k} - t \pi) \cdot (1 + 3 \sqrt{3} t (\pi / k))^{- k / 2}.
\end{align*}
\end{allowdisplaybreaks}
We denote the upper bound by $A_1$ and the lower bound by $B_1$. First, we have $A_1 |_{t = 0} = B_1 |_{t = 0} = 0$ and $\frac{\partial}{\partial t} A_1 |_{t = 0} = \frac{\partial}{\partial t} B_1 |_{t = 0} = 0$. Second, let $C_1 = (5 \sqrt{3} / 2) \pi^2 (- \cos \alpha_{7, k}) (\tan \alpha_{7, k} + 11 / (5 \sqrt{3}))$, then we have $\frac{\partial^2}{\partial t^2} A_1 |_{t = 0} = C_1 + 6 \pi^2 (- \cos(2 \pi / 3 + \alpha_{7, k})) / k$ and $\frac{\partial^2}{\partial t^2} B_1 |_{t = 0} = C_1 - (27 / 2) \pi^2 \cos(4 \pi / 3 + \alpha_{7, k}) / k$. Finally, since $\tan (3 \pi / 2 - 2 \alpha_7) + 11 / (5 \sqrt{3}) = 0$, we have $B_1 > 0$ if $\alpha_{7, k} > 3 \pi / 2 - 2 \alpha_7$, and we have $A_1 < 0$ if $\alpha_{7, k} < 3 \pi / 2 - 2 \alpha_7$ for $k$ large enough and for $t$ small enough.

Similarly, we have $0 < \beta_{7, k} + (t / 2) \pi < \pi / 2$, $\pi < 4 \pi / 3 + \beta_{7, k} - (3 t / 4) \pi < {\beta_{7, k, 1}}' < 4 \pi / 3 + \beta_{7, k} - d_{2, 1} (t / 2) \pi < 3 \pi / 2$, and $\pi / 2 < 2 \pi / 3 + \beta_{7, k} + d_{2, 2} (t / 2) \pi < {\beta_{7, k, 2}}' < 2 \pi / 3 + \beta_{7, k} + (t / 4) \pi < \pi$. Thus
\begin{allowdisplaybreaks}
\begin{align*}
&\cos(\beta_{7, k} + (t / 2) \pi)\\
 &\quad + \cos(4 \pi / 3 + \beta_{7, k} - d_{2, 1} (t / 2) \pi) \cdot (1 + (\sqrt{3} / 2) t (\pi / k) + (1 / 2) t^2 (\pi^2 / k^2))^{- k / 2}\\
 &\quad \qquad + \cos(2 \pi / 3 + \beta_{7, k} + d_{2, 2} (t / 2) \pi) \cdot e^{- (3 \sqrt{3} / 4) \pi t}\\
 &\qquad > |\cos(k \theta_2 / 2)|
 - \left| Re\left\{ 2^k \cdot \left( e^{i \theta_2 / 2} - \sqrt{7} e^{- i \theta_2 / 2} \right)^{- k} \right\} \right|\\
 &\qquad \qquad \qquad \qquad - \left| Re\left\{ 2^k \cdot \left( 3 e^{i \theta_2 / 2} - \sqrt{7} e^{- i \theta_2 / 2} \right)^{- k} \right\} \right|\\
 &\qquad > \cos(\beta_{7, k} + (t / 2) \pi)
 + \cos(4 \pi / 3 + \beta_{7, k} - (3 t / 4) \pi) \cdot e^{- (\sqrt{3} / 4) \pi t}\\
 &\qquad \qquad \qquad \qquad + \cos(2 \pi / 3 + \beta_{7, k} + (t / 4) \pi) \cdot (1 + (3 \sqrt{3} / 2) t (\pi / k))^{- k / 2}.
\end{align*}
\end{allowdisplaybreaks}
We denote the upper bound by $A_2$ and the lower bound by $B_2$. First, we have $A_2 |_{t = 0} = B_2 |_{t = 0} = 0$ and $\frac{\partial}{\partial t} A_2 |_{t = 0} = \frac{\partial}{\partial t} B_2 |_{t = 0} = 0$. Second, let $C_2 = \sqrt{3} \pi^2 \cos \beta_{7, k} (\sqrt{3} / 12 - \tan \beta_{7, k})$, then we have $\frac{\partial^2}{\partial t^2} A_2 |_{t = 0} = C_2 + (1 / 8) \pi^2 (- \cos(4 \pi / 3 + \beta_{7, k})) / k$ and $\frac{\partial^2}{\partial t^2} B_2 |_{t = 0} = C_2 - (27 / 8) \pi^2 (- \cos(2 \pi / 3 + \beta_{7, k})) / k$. Finally, since $\sqrt{3} / 12 - \tan (5 \pi / 6 - 2 \alpha_7) = 0$ and $\alpha_{7, k} = \beta_{7, k} + 2 \pi / 3$, we have $B_2 > 0$ if $\alpha_{7, k} < 3 \pi / 2 - 2 \alpha_7$, and we have $A_2 < 0$ if $\alpha_{7, k} > 3 \pi / 2 - 2 \alpha_7$ for $k$ large enough and for $t$ small enough.

In conclusion, if $k$ is large enough, then $|{R_{7, 1}^{*}}'|$ and $|{R_{7, 2}^{*}}'|$ is small enough, and then we have one more zero on the arc $A_{7, 1}^{*}$ when $\alpha_{7, k} > 3 \pi / 2 - 2 \alpha_7$, and we have one more zero on the arc $A_{7, 2}^{*}$ when $\alpha_{7, k} < 3 \pi / 2 - 2 \alpha_7$. However, if $k$ is small, the proof not clear.\\

\subsubsection{\bfseries The case ``$k \equiv 4 \pmod{6}$ and $217 \pi / 360 < \alpha_{7, k} < 73 \pi / 120$''}
Let $t > 0$ be small enough, then we have
\begin{allowdisplaybreaks}
\begin{align*}
&- \cos(\alpha_{7, k} - (t / 2) \pi)
 - \cos(4 \pi / 3 + \alpha_{7, k} + (3 t / 2) \pi) \cdot e^{- \sqrt{3} \pi t}\\
 &\qquad \qquad \qquad \qquad - \cos(2 \pi / 3 + \alpha_{7, k} - t \pi) \cdot (1 + 3 \sqrt{3} t (\pi / k))^{- k / 2}\\
 &\qquad > |\cos(k \theta_1 / 2)|
 - \left| Re\left\{ \left( 2 e^{i \theta_1 / 2} + \sqrt{7} e^{- i \theta_1 / 2} \right)^{- k} \right\} \right|\\
 &\qquad \qquad \qquad \qquad + \left| Re\left\{ \left( 3 e^{i \theta_1 / 2} + \sqrt{7} e^{- i \theta_1 / 2} \right)^{- k} \right\} \right|\\
 &\qquad > - \cos(\alpha_{7, k} - (t / 2) \pi)
 - \cos(4 \pi / 3 + \alpha_{7, k} + d_{1, 1} (t / 2) \pi) \cdot (1 + 2 \sqrt{3} t (\pi / k))^{- k / 2}\\
 &\qquad \qquad \qquad - \cos(2 \pi / 3 + \alpha_{7, k} - d_{1, 2} (t / 2) \pi) \cdot e^{- (3 \sqrt{3} / 2) \pi t}.
\end{align*}
\end{allowdisplaybreaks}
We denote the upper bound by $A_1$ and the lower bound by $B_1$. First, we have $A_1 |_{t = 0} = B_1 |_{t = 0} = 0$. Second, we have $\frac{\partial}{\partial t} A_1 |_{t = 0} = \frac{\partial}{\partial t} B_1 |_{t = 0} = (3 \pi / 2) (- \cos \alpha_{7, k}) (5 / \sqrt{3} + \tan \alpha_{7, k})$. Finally, we have $B_1 > 0$ if $\alpha_{7, k} > \pi - \alpha_7$, and we have $A_1 < 0$ if $\alpha_{7, k} < \pi - \alpha_7$ for $t$ small enough.

Similarly, we have
\begin{allowdisplaybreaks}
\begin{align*}
&\cos(\beta_{7, k} + (t / 2) \pi)\\
 &\quad + \cos(4 \pi / 3 + \beta_{7, k} - (3 t / 4) \pi) \cdot (1 + (\sqrt{3} / 2) t (\pi / k) + (1 / 2) t^2 (\pi^2 / k^2))^{- k / 2}\\
 &\qquad \qquad \qquad \qquad + \cos(2 \pi / 3 + \beta_{7, k} + (t / 4) \pi) \cdot (1 + (3 \sqrt{3} / 2) t (\pi / k))^{- k / 2}\\
 &\qquad > |\cos(k \theta_2 / 2)|
 - \left| Re\left\{ 2^k \cdot \left( e^{i \theta_2 / 2} - \sqrt{7} e^{- i \theta_2 / 2} \right)^{- k} \right\} \right|\\
 &\qquad \qquad \qquad \qquad + \left| Re\left\{ 2^k \cdot \left( 3 e^{i \theta_2 / 2} - \sqrt{7} e^{- i \theta_2 / 2} \right)^{- k} \right\} \right|\\
 &\qquad > \cos(\beta_{7, k} + (t / 2) \pi)
 + \cos(4 \pi / 3 + \beta_{7, k} - d_{2, 1} (t / 2) \pi) \cdot e^{- (\sqrt{3} / 4) \pi t}\\
 &\quad \qquad + \cos(2 \pi / 3 + \beta_{7, k} + d_{2, 2} (t / 2) \pi) \cdot e^{- (3 \sqrt{3} / 4) \pi t}.
\end{align*}
\end{allowdisplaybreaks}
We denote the upper bound by $A_2$ and the lower bound by $B_2$. First, we have $A_2 |_{t = 0} = B_2 |_{t = 0} = 0$. Second, we have $\frac{\partial}{\partial t} A_2 |_{t = 0} = \frac{\partial}{\partial t} B_2 |_{t = 0} = (3 \pi / 2) \cos \beta_{7, k} (2 / \sqrt{3} - \tan \beta_{7, k})$. Finally, since $2 / \sqrt{3} - \tan (2 \pi / 3 - \alpha_7) = 0$ and $\alpha_{7, k} = \beta_{7, k} + 4 \pi / 3$, we have $B_2 > 0$ if $\alpha_{7, k} < \pi - \alpha_7$, and we have $A_2 < 0$ if $\alpha_{7, k} > \pi - \alpha_7$ for $t$ small enough.

In conclusion, if $k$ is large enough, then we have one more zero on the arc $A_{7, 1}^{*}$ when $\alpha_{7, k} > \pi - \alpha_7$, and we have one more zero on the arc $A_{7, 2}^{*}$ when $\alpha_{7, k} < \pi - \alpha_7$.\\

The remaining subsections in this section detail the proofs of lemmas \ref{lem-f7s1}, ..., \ref{lem-f7s8}\\

\subsection{Proof of Lemma \ref{lem-f7s1}} \label{subsec-pf71}

\begin{trivlist}
\item[(1)]
Let $k \geqslant 4$ be an even integer divisible by $4$.

First, we consider the case $N = 1$. Then, we can write:
\begin{equation*}
F_{k, 7}^{*} (\pi / 2) = F_{k, 7, 1}^{*}(\pi / 2) = 2 \cos(k \pi / 4) + R_{7, \pi / 2}^{*}
\end{equation*}
where $R_{7, \pi / 2}^{*}$ denotes the remaining terms.

We have $v_{k}(c, d, \pi / 2) = 1 / (c^2 + 7 d^2)^{k/2}$. Now we will consider the following cases: $N = 2$, $5$, $10$, $13$, $17$, and $N \geqslant 25$. We have
\begin{allowdisplaybreaks}
\begin{align*}
&\text{When $N = 2$,}&
v_k(1, 1, \pi / 2) &\leqslant (1 / 8)^{k/2}.\\
&\text{When $N = 5$,}&
v_k(1, 2, \pi / 2) &\leqslant (1 / 29)^{k/2},
&v_k(2, 1, \pi / 2) &\leqslant (1 / 11)^{k/2}.\\
&\text{When $N = 10$,}&
v_k(1, 3, \pi / 2) &\leqslant (1 / 8)^{k},
&v_k(3, 1, \pi / 2) &\leqslant (1 / 4)^{k}.\\
&\text{When $N = 13$,}&
v_k(2, 3, \pi / 2) &\leqslant (1 / 69)^{k/2},
&v_k(3, 2, \pi / 2) &\leqslant (1 / 37)^{k/2}.\\
&\text{When $N = 17$,}&
v_k(1, 4, \pi / 2) &\leqslant (1 / 113)^{k/2},
&v_k(4, 1, \pi / 2) &\leqslant (1 / 23)^{k/2}.\\
&\text{When $N \geqslant 25$,}&
c^2 + 7 d^2 &\geqslant N,
\end{align*}
\end{allowdisplaybreaks}
and the number of terms with $c^2 + d^2 = N$ is not more than $(144 / 35) N^{1/2}$ for $N \geqslant 25$. Then
\begin{equation*}
|R_{7, \pi / 2}^{*}|_{N \geqslant 25}
\leqslant \frac{288}{35(k-3)} \left(\frac{1}{24}\right)^{(k - 3) / 2}.
\end{equation*}
Furthermore,
\begin{align*}
|R_{7, \pi / 2}^{*}|
 &\leqslant
 4 \left(\frac{1}{8}\right)^{k/2} + 4 \left(\frac{1}{11}\right)^{k/2}
 + \cdots + 4 \left(\frac{1}{113}\right)^{k/2}
 + \frac{288}{35(k-3)} \left(\frac{1}{24}\right)^{(k - 3) / 2},\\
 &\leqslant 1.80820... \quad (k \geqslant 4)
\end{align*}\quad

\item[(2)]

First, we consider the case $N = 1$. Then, we can write:
\begin{equation*}
F_{4, 7}^{*} (5 \pi / 6) = 2 \cos(10 \pi / 3) + R_{7, 4}^{*} = 1 + R_{7, 4}^{*}
\end{equation*}
where
\begin{equation*}
R_{7, 4}^{*}
 \leqslant \frac{1}{2} \sum_{\begin{subarray}{c} (c,d)=1\\ 7 \nmid d \; N > 1\end{subarray}}(c e^{i 5 \pi / 12} + \sqrt{7} d e^{-i 5 \pi / 12})^{-4} + \frac{1}{2} \sum_{\begin{subarray}{c} (c,d)=1\\ 7 \nmid d \; N > 1\end{subarray}}(c e^{- i 5 \pi / 12} + \sqrt{7} d e^{i 5 \pi / 12})^{-4}.
\end{equation*}

We want to prove $F_{4, 7}^{*} (5 \pi / 6) > 0$, but it is too difficult to prove that $|R_{7, 4}^{*}| < 1$. However, we have only to prove $R_{7, 4}^{*} > -1$.

Let $u_0(c, d) := (c e^{i 5 \pi / 12} + \sqrt{7} d e^{-i 5 \pi / 12})^{-4} + (c e^{- i 5 \pi / 12} + \sqrt{7} d e^{i 5 \pi / 12})^{-4}$, and let $u(c, d) := u_0(c, d) + u_0(c, -d) + u_0(d, c) + u_0(d, -c)$ for every pair $(c, d)$ such that $c \neq d$ and $7 \nmid c$, and $u_1(c, d) := u_0(c, d) + u_0(c, -d)$ for every pair $(c, d)$ such that $c = d$ or $7 \mid d$. Now we will consider the following cases: $N = 2$, $5$, $10$,..., $197$, and $N \geqslant 202$. We have the following:
\begin{allowdisplaybreaks}
\begin{align*}
&\text{When $N = 2$,}&u_1(1, 1) &\geqslant -0.08151.\\
&\text{When $N = 5$,}&u(1, 2) &\geqslant -0.19373.&
&\text{When $N = 10$,}&u(1, 3) &\geqslant 0.24147.\\
&\text{When $N = 13$,}&u(2, 3) &\geqslant -0.02162.&
&\text{When $N = 17$,}&u(1, 4) &\geqslant -0.07736.\\
&\text{When $N = 25$,}&u(3, 4) &\geqslant -0.00313.&
&\text{When $N = 26$,}&u(1, 5) &\geqslant -0.02262.\\
&\text{When $N = 29$,}&u(2, 5) &\geqslant 0.03569.&
&\text{When $N = 34$,}&u(3, 5) &\geqslant -0.00503.\\
&\text{When $N = 37$,}&u(1, 6) &\geqslant -0.00586.&
&\text{When $N = 41$,}&u(4, 5) &\geqslant -0.00083.\\
&\text{When $N = 50$,}&u_1(1, 7) &\geqslant 0.00000.&
&\text{When $N = 53$,}&u_1(2, 7) &\geqslant 0.00000.\\
&\text{When $N = 58$,}&u_1(3, 7) &\geqslant 0.00000.&
&\text{When $N = 61$,}&u(5, 6) &\geqslant -0.00033.\\
&\text{When $N = 65$,}&u(1, 8) &\quad + u_1(4, 7)& \geqslant &-0.00052.\\
&\text{When $N = 73$,}&u(3, 8) &\geqslant 0.00692.&
&\text{When $N = 74$,}&u_1(5, 7) &\geqslant -0.00006.\\
&\text{When $N = 82$,}&u(1, 9) &\geqslant -0.00014.\\
&\text{When $N = 85$,}&u(2, 9) &\quad + u_1(6, 7)& \geqslant &-0.00282.\\
&\text{When $N = 89$,}&u(5, 8) &\geqslant -0.00064.&
&\text{When $N = 97$,}&u(4, 9) &\geqslant 0.00099.\\
&\text{When $N = 101$,}&u(1, 10) &\geqslant -0.00003.&
&\text{When $N = 106$,}&u(5, 9) &\geqslant -0.00064.\\
&\text{When $N = 109$,}&u(3, 10) &\geqslant -0.00014.&
&\text{When $N = 113$,}&u_1(8, 7) &\geqslant -0.00007.\\
&\text{When $N = 122$,}&u(1, 11) &\geqslant 0.00002.&
&\text{When $N = 125$,}&u(2, 11) &\geqslant -0.00073.\\
&\text{When $N = 130$,}&u(3, 11) &\quad + u_1(9, 7)& \geqslant &-0.00115.\\
&\text{When $N = 137$,}&u(4, 11) &\geqslant 0.00195.\\
&\text{When $N = 145$,}&u(1, 12) &\quad + u(8, 9)& \geqslant &-0.00003.\\
&\text{When $N = 146$,}&u(5, 11) &\geqslant 0.00025.&
&\text{When $N = 149$,}&u_1(10, 7) &\geqslant -0.00014.\\
&\text{When $N = 157$,}&u(6, 11) &\geqslant -0.00030.&
&\text{When $N = 169$,}&u(1, 2) &\geqslant 0.00077.\\
&\text{When $N = 170$,}&u(1, 13) &\quad + u_1(11, 7)& \geqslant &-0.00015.\\
&\text{When $N = 173$,}&u(2, 13) &\geqslant -0.00021.&
&\text{When $N = 178$,}&u(3, 13) &\geqslant -0.00068.\\
&\text{When $N = 181$,}&u(9, 10) &\geqslant -0.00004.\\
&\text{When $N = 185$,}&u(4, 13) &\quad + u(8, 11)& \geqslant &0.00003.\\
&\text{When $N = 193$,}&u_1(12, 7) &\geqslant -0.00018.&
&\text{When $N = 194$,}&u(5, 13) &\geqslant 0.00093.\\
&\text{When $N = 197$,}&u_1(1, 14) &\geqslant 0.00000.
\end{align*}
\end{allowdisplaybreaks}
When $N \geqslant 202$, for the case of $c d < 0$, we put $X + Y i = (c e^{i 5 \pi / 12} + \sqrt{7} d e^{-i 5 \pi / 12})^2 = - (\sqrt{3} / 2) (c^2 + 7 d^2) + 2 \sqrt{7} c d + (1 / 2) (c^2 - 7 d^2) i$. Then, $|Y| - |X| < - \frac{\sqrt{3} - 1}{2} c^2 - \frac{\sqrt{3} - 1}{2} d^2 + 2 \sqrt{7} c d < 0$. Thus, $(c e^{i 5 \pi / 12} + \sqrt{7} d e^{-i 5 \pi / 12})^{-4} + (c e^{- i 5 \pi / 12} + \sqrt{7} d e^{i 5 \pi / 12})^{-4} > 0$.

For the case of $c d > 0$, we have $c^2 + 7 d^2 - \sqrt{21}|c d| > (2 / 9) N$, and the number of terms with $c^2 + d^2 = N$ is not more than $(13 / 7) N^{1/2}$ for $N \geqslant 144$. However, this bound is too large. We must consider some sub-cases. 

For the case of $|c| < |d|$, we have $c^2 + 7 d^2 - \sqrt{21}|c d| > (3 / 2) N$ and $|c| < (1 / \sqrt{2}) N^{1 / 2}$, $1 / \sqrt{2} > 7 / 10$. For the case of $|d| \leqslant |c| < (6 / \sqrt{21}) |d|$, we have $c^2 + 7 d^2 - \sqrt{21}|c d| > N$ and $|c| < (6 / \sqrt{57}) N^{1 / 2}$, $6 / \sqrt{57} > 7 / 9$. For the case of $(6 / \sqrt{21}) |d| \leqslant |c| < \sqrt{7 / 3} |d|$, we have $c^2 + 7 d^2 - \sqrt{21}|c d| > (1 / 2) N$ and $|c| < \sqrt{7 / 10} N^{1 / 2}$, $\sqrt{7 / 10} > 5 / 6$. For the case of $c d > 0$ and $\sqrt{7 / 3} |d| \leqslant |c| < (22 / 3 \sqrt{21})|d|$, we have $c^2 + 7 d^2 - \sqrt{21}|c d| > (1 / 3) N$ and $|c| < (22 / \sqrt{673}) N^{1 / 2}$, $22 / \sqrt{673} > 22 / 25$.

In conclusion, we have
\begin{allowdisplaybreaks}
\begin{align*}
R_{7, 4}^{*} &\Big|_{N \geqslant 202, c d > 0}\\
 &\geqslant - \frac{13}{7} \Big( \frac{7}{10} N^{1 / 2} \sum_{N \geqslant 202} \left(\frac{3}{2} N\right)^{- k / 2} + \frac{7}{90} N^{1 / 2} \sum_{N \geqslant 202} N^{- k / 2} + \frac{1}{18} N^{1 / 2} \sum_{N \geqslant 202} \left(\frac{1}{2} N\right)^{- k / 2}\\
 &\quad + \frac{7}{150} N^{1 / 2} \sum_{N \geqslant 202} \left(\frac{1}{3} N\right)^{- k / 2} + \frac{3}{25} N^{1 / 2} \sum_{N \geqslant 202} \left(\frac{2}{9} N\right)^{- k / 2} \Big),\\
 &= - \frac{13}{7} \left( \frac{7}{10} \cdot \frac{4}{9} + \frac{7}{90} \cdot 1 + \frac{1}{18} \cdot 4 + \frac{7}{150} \cdot 9 + \frac{3}{25} \cdot \frac{81}{4} \right) \sum_{N \geqslant 202} N^{(1 - k) / 2},\\
 &= - \frac{7579}{630 \sqrt{201}}.
\end{align*}
\end{allowdisplaybreaks}

Furthermore,
\begin{equation*}
R_{7, 4}^{*}
 \geqslant -0.13164 - \frac{7579}{630 \sqrt{201}}
 = - 0.98018... \quad (k \geqslant 4)
\end{equation*}\quad

\item[(3)]
Let $k \geqslant 6$ be an even integer.

First, we consider the case $N = 1$. Since $2 \pi / 3 < \pi / 2 + \alpha_7$, we can write:
\begin{equation*}
F_{k, 7}^{*} (\theta) = F_{k, 7, 1}^{*}(\theta) = 2 \cos(k \theta / 2) + R_{7, 2 \pi / 3}^{*} \quad \text{for} \: \theta \in [\pi / 2, 2 \pi / 3],
\end{equation*}
where $R_{7, 2 \pi / 3}^{*}$ denotes the remaining terms. Now we will consider the following cases: $N = 2$ and $N \geqslant 5$. Considering $- 1 / 2 \leqslant \cos\theta \leqslant 0$, we have
\begin{align*}
|R_{7, 2 \pi / 3}^{*}|
 &\leqslant
 2 \left(\frac{1}{8 - \sqrt{7}}\right)^{k/2}
 + 2 \left(\frac{1}{8}\right)^{k/2}
 + \frac{576}{7 (k - 3)} \left(\frac{7}{20}\right)^{(k - 3) / 2},\\
 &\leqslant 1.19293... \quad (k \geqslant 6)
\end{align*}\quad

\item[(4)]
Let $k \geqslant 6$ be an even integer.

First, we consider the case $N = 1$. Then, we can write:
\begin{equation*}
F_{k, 7}^{*} (\theta) = F_{k, 7, 2}^{*}(\theta - 2 \pi / 3) = 2 \cos(k (\theta - 2 \pi / 3) / 2) + R_{7, \pi}^{*} \quad \text{for} \: \theta \in [\pi, 7 \pi / 6],
\end{equation*}
where $R_{7, \pi}^{*}$ denotes the remaining terms. Now we will consider the following cases: $N = 2$, $5$, $10$,..., $82$, and $N \geqslant 85$. Considering $0 \leqslant \cos\theta \leqslant 1 / 2$, we have
\begin{allowdisplaybreaks}
\begin{align*}
&\text{When $N = 2$,}&
2^k \cdot v_k(1, 1, \theta) &\leqslant (1 / 2)^{k/2}, \quad
&2^k \cdot v_k(1, - 1, \theta) &\leqslant \left(4 / \left(8 - \sqrt{7}\right)\right)^{k/2}.\\
&\text{When $N = 5$,}&
v_k(1, 2, \theta) &\leqslant (1 / 29)^{k/2},
&v_k(1, - 2, \theta) &\leqslant 1 / \left(29 - 2 \sqrt{7}\right)^{k/2},\\
&&v_k(2, 1, \theta) &\leqslant (1 / 11)^{k/2},
&v_k(2, - 1, \theta) &\leqslant 1 / \left(11 - 2 \sqrt{7}\right)^{k/2}.\\
&\text{When $N = 10$,}&
2^k \cdot v_k(1, 3, \theta) &\leqslant (1 / 4)^{k},
&2^k \cdot v_k(1, - 3, \theta) &\leqslant \left(4 / \left(64 - 3 \sqrt{7}\right)\right)^{k/2},\\
&&2^k \cdot v_k(3, 1, \theta) &\leqslant (1 / 2)^{k},
&2^k \cdot v_k(3, - 1, \theta) &\leqslant \left(4 / \left(16 - 3 \sqrt{7}\right)\right)^{k/2}.\\
&\text{When $N = 13$,}&
v_k(2, 3, \theta) &\leqslant (1 / 69)^{k/2},
&v_k(2, - 3, \theta) &\leqslant 1 / \left(69 - 6 \sqrt{7}\right)^{k/2},\\
&&v_k(3, 2, \theta) &\leqslant (1 / 37)^{k/2},
&v_k(3, - 2, \theta) &\leqslant 1 / \left(37 - 6 \sqrt{7}\right)^{k/2}.\\
&\text{When $N = 17$,}&
v_k(1, 4, \theta) &\leqslant (1 / 113)^{k/2},
&v_k(1, - 4, \theta) &\leqslant 1 / \left(113 - 4 \sqrt{7}\right)^{k/2},\\
&&v_k(4, 1, \theta) &\leqslant (1 / 23)^{k/2},
&v_k(4, - 1, \theta) &\leqslant 1 / \left(23 - 4 \sqrt{7}\right)^{k/2}.\\
&\text{When $N = 25$,}&
v_k(3, 4, \theta) &\leqslant (1 / 11)^{k},
&v_k(3, - 4, \theta) &\leqslant 1 / \left(121 - 12 \sqrt{7}\right)^{k/2},\\
&&v_k(4, 3, \theta) &\leqslant (1 / 79)^{k/2},
&v_k(4, - 3, \theta) &\leqslant 1 / \left(79 - 12 \sqrt{7}\right)^{k/2}.\\
&\text{When $N = 26$,}&
2^k \cdot v_k(1, 5, \theta) &\leqslant (1 / 44)^{k/2},
&2^k \cdot v_k(1, - 5, \theta) &\leqslant \left(4 / \left(176 - 5 \sqrt{7}\right)\right)^{k/2},\\
&&2^k \cdot v_k(5, 1, \theta) &\leqslant (1 / 8)^{k/2},
&2^k \cdot v_k(5, - 1, \theta) &\leqslant \left(4 / \left(32 - 5 \sqrt{7}\right)\right)^{k/2}.\\
&\text{When $N = 29$,}&
v_k(2, 5, \theta) &\leqslant (1 / 179)^{k/2},
&v_k(2, - 5, \theta) &\leqslant 1 / \left(179 - 10 \sqrt{7}\right)^{k/2},\\
&&v_k(5, 2, \theta) &\leqslant (1 / 53)^{k/2},
&v_k(5, - 2, \theta) &\leqslant 1 / \left(53 - 10 \sqrt{7}\right)^{k/2}.\\
&\text{When $N = 34$,}&
2^k \cdot v_k(3, 5, \theta) &\leqslant (1 / 46)^{k/2},
&2^k \cdot v_k(3, - 5, \theta) &\leqslant \left(4 / \left(184 - 15 \sqrt{7}\right)\right)^{k/2},\\
&&2^k \cdot v_k(5, 3, \theta) &\leqslant (1 / 22)^{k},
&2^k \cdot v_k(5, - 3, \theta) &\leqslant \left(4 / \left(88 - 15 \sqrt{7}\right)\right)^{k/2}.\\
&\text{When $N = 37$,}&
v_k(1, 6, \theta) &\leqslant (1 / 253)^{k/2},
&v_k(1, - 6, \theta) &\leqslant 1 / \left(253 - 6 \sqrt{7}\right)^{k/2},\\
&&v_k(6, 1, \theta) &\leqslant (1 / 43)^{k/2},
&v_k(6, - 1, \theta) &\leqslant 1 / \left(43 - 6 \sqrt{7}\right)^{k/2}.\\
&\text{When $N = 41$,}&
v_k(4, 5, \theta) &\leqslant (1 / 191)^{k/2},
&v_k(4, - 5, \theta) &\leqslant 1 / \left(191 - 20 \sqrt{7}\right)^{k/2},\\
&&v_k(5, 4, \theta) &\leqslant (1 / 137)^{k/2},
&v_k(5, - 4, \theta) &\leqslant 1 / \left(137 - 20 \sqrt{7}\right)^{k/2}.\\
&\text{When $N = 50$,}&
2^k \cdot v_k(1, 7, \theta) &\leqslant (1 / 86)^{k/2},
&2^k \cdot v_k(1, - 7, \theta) &\leqslant \left(4 / \left(344 - 7 \sqrt{7}\right)\right)^{k/2},\\
&\text{When $N = 53$,}&
v_k(2, 7, \theta) &\leqslant (1 / 347)^{k/2},
&v_k(2, - 7, \theta) &\leqslant 1 / \left(347 - 14 \sqrt{7}\right)^{k/2},\\
&\text{When $N = 58$,}&
2^k \cdot v_k(3, 7, \theta) &\leqslant (1 / 88)^{k/2},
&2^k \cdot v_k(3, - 7, \theta) &\leqslant \left(4 / \left(352 - 21 \sqrt{7}\right)\right)^{k/2},\\
&\text{When $N = 61$,}&
v_k(5, 6, \theta) &\leqslant (1 / 277)^{k/2},
&v_k(5, - 6, \theta) &\leqslant 1 / \left(277 - 30 \sqrt{7}\right)^{k/2},\\
&&v_k(6, 5, \theta) &\leqslant (1 / 211)^{k/2},
&v_k(6, - 5, \theta) &\leqslant 1 / \left(211 - 30 \sqrt{7}\right)^{k/2}.\\
&\text{When $N = 65$,}&
v_k(1, 8, \theta) &\leqslant (1 / 449)^{k/2},
&v_k(1, - 8, \theta) &\leqslant 1 / \left(449 - 8 \sqrt{7}\right)^{k/2},\\
&&v_k(8, 1, \theta) &\leqslant (1 / 71)^{k/2},
&v_k(8, - 1, \theta) &\leqslant 1 / \left(71 - 8 \sqrt{7}\right)^{k/2},\\
&&v_k(4, 7, \theta) &\leqslant (1 / 359)^{k/2},
&v_k(4, - 7, \theta) &\leqslant 1 / \left(359 - 28 \sqrt{7}\right)^{k/2},\\
&\text{When $N = 73$,}&
v_k(3, 8, \theta) &\leqslant (1 / 457)^{k/2},
&v_k(3, - 8, \theta) &\leqslant 1 / \left(457 - 24 \sqrt{7}\right)^{k/2},\\
&&v_k(8, 3, \theta) &\leqslant (1 / 127)^{k/2},
&v_k(8, - 3, \theta) &\leqslant 1 / \left(127 - 24 \sqrt{7}\right)^{k/2}.\\
&\text{When $N = 74$,}&
2^k \cdot v_k(5, 7, \theta) &\leqslant (1 / 92)^{k/2},
&2^k \cdot v_k(5, - 7, \theta) &\leqslant \left(4 / \left(368 - 35 \sqrt{7}\right)\right)^{k/2},\\
&\text{When $N = 82$,}&
2^k \cdot v_k(1, 9, \theta) &\leqslant (1 / 142)^{k/2},
&2^k \cdot v_k(1, - 9, \theta) &\leqslant \left(4 / \left(568 - 9 \sqrt{7}\right)\right)^{k/2},\\
&&2^k \cdot v_k(9, 1, \theta) &\leqslant (1 / 22)^{k/2},
&2^k \cdot v_k(9, - 1, \theta) &\leqslant \left(4 / \left(88 - 9 \sqrt{7}\right)\right)^{k/2}.\\
&\text{When $N \geqslant 85$,}&
|c e^{i \theta / 2} \pm \sqrt{7} d &e^{-i \theta / 2}|^2 \geqslant 5 N / 7,
\end{align*}
\end{allowdisplaybreaks}
and the number of terms with $c^2 + d^2 = N$ is not more than $(27 / 7) N^{1/2}$ for $N \geqslant 64$. Then
\begin{equation*}
|R_{7, \pi}^{*}|_{N \geqslant 85}
\leqslant \frac{1296 \sqrt{21}}{k - 3} \left(\frac{1}{15}\right)^{(k - 3) / 2}.
\end{equation*}
Furthermore,
\begin{align*}
|R_{7, \pi}^{*}|
 &\leqslant
 2 \left(\frac{4}{8 - \sqrt{7}}\right)^{k/2}
 + \cdots + 2 \left(\frac{1}{457}\right)^{k/2}
 + \frac{1296 \sqrt{21}}{k - 3} \left(\frac{1}{15}\right)^{(k - 3) / 2},\\
 &\leqslant 1.98681... \quad (k \geqslant 6).
\end{align*}\quad

\item[(5)]
Let $k \geqslant 8$ be an even integer.

First, we consider the case of $N = 1$. Since $5 \pi / 6 < \pi / 2 + \alpha_5$, we can write:
\begin{equation*}
F_{k, 7}^{*} (\theta) = F_{k, 7, 1}^{*}(\theta) = 2 \cos(k \theta / 2) + R_{7, 5 \pi / 6}^{*} \quad \text{for} \: \theta \in [\pi / 2, 5 \pi / 6],
\end{equation*}
where $R_{7, 5 \pi / 6}^{*}$ denotes the remaining terms. Now we will consider the following cases: $N = 2, 5, 10$, and $N \geqslant 13$. Considering $- \sqrt{3} / 2 \leqslant \cos\theta \leqslant 0$, we have
\begin{allowdisplaybreaks}
\begin{align*}
&\text{When $N = 2$,}&
v_k(1, 1, \theta) &\leqslant 1 / \left(8 - \sqrt{21}\right)^{k/2}, \quad
&v_k(1, - 1, \theta) &\leqslant (1 / 8)^{k/2}.\\
&\text{When $N = 5$,}&
v_k(1, 2, \theta) &\leqslant 1 / \left(29 - 2 \sqrt{21}\right)^{k/2},
&v_k(1, - 2, \theta) &\leqslant (1 / 29)^{k/2},\\
&&v_k(2, 1, \theta) &\leqslant 1 / \left(11 - 2 \sqrt{21}\right)^{k/2},
&v_k(2, - 1, \theta) &\leqslant (1 / 11)^{k/2}.\\
&\text{When $N = 10$,}&
v_k(1, 3, \theta) &\leqslant 1 / \left(64 - 3 \sqrt{21}\right)^{k/2},
&v_k(1, - 3, \theta) &\leqslant (1 / 8)^{k},\\
&&v_k(3, 1, \theta) &\leqslant 1 / \left(16 - 3 \sqrt{21}\right)^{k/2},
&v_k(3, - 1, \theta) &\leqslant (1 / 4)^{k}.\\
&\text{When $N \geqslant 13$,}&
|c e^{i \theta / 2} \pm \sqrt{7} d &e^{-i \theta / 2}|^2 \geqslant 2 N / 9,
\end{align*}
\end{allowdisplaybreaks}
and the number of terms with $c^2 + d^2 = N$ is not more than $(36 / 7) N^{1/2}$ for $N \geqslant 13$. Then
\begin{equation*}
|R_{7, 5 \pi / 6}^{*}|_{N \geqslant 13}
\leqslant \frac{1728 \sqrt{3}}{7(k-3)} \left(\frac{3}{8}\right)^{(k - 3) / 2}.
\end{equation*}
Furthermore,
\begin{align*}
|R_{7, \pi / 2}^{*}|
 &\leqslant
 2 \left(\frac{1}{11 - 2 \sqrt{21}}\right)^{k/2}
 + \cdots + 2 \left(\frac{1}{8}\right)^{k}
 + \frac{1728 \sqrt{3}}{7(k-3)} \left(\frac{3}{8}\right)^{(k - 3) / 2},\\
 &\leqslant 1.96057... \quad (k \geqslant 8).
\end{align*}
\end{trivlist}
\begin{flushright}$\square$\end{flushright}\quad

In the proofs of the remaining lemmas, we will use the algorithm of the subsection \ref{subsec-alg}. Furthermore, we have $X_1 = v_k(2, 1, \theta_1)^{- 2 / k} \geqslant 1 + 2 \sqrt{3} t (\pi / k)$ and $X_2 = v_k(3, 1, \theta_1)^{- 2 / k} \geqslant 1 + 3 \sqrt{3} t (\pi / k)$ in the proof of lemmas \ref{lem-f7s2} (1)--(3), \ref{lem-f7s3}, and \ref{lem-f7s5}, and we have $X_1 = (1 / 4) \; v_k(1, -1, \theta_2)^{- 2 / k} \geqslant 1 + (\sqrt{3} / 2) t (\pi / k)$ and\\ $X_2 = (1 / 4) \; v_k(3, -1, \theta_2)^{- 2 / k} \geqslant 1 + (3 \sqrt{3} / 2) t (\pi / k)$ in the proof of lemmas \ref{lem-f7s2} (4)--(7), \ref{lem-f7s4}, and \ref{lem-f7s6}.\\

\subsection{Proof of Lemma \ref{lem-f7s2}}\quad\vspace{-0.1in} \\

\noindent
We have $c_0 \leqslant \cos({c_0}' \pi)$.
\begin{trivlist}
\item[(1)]
Let $(t, b, s, k_0) = (1 / 3, 28, 3, 10)$, then we can define $a_1 := 53/10$. Furthermore, let\\ $({c_1}', c_{0, 1}, a_{1, 1}, u_1) = (1, 5 / 14, 53 / 14, 2 \sqrt{3} / 3)$, then we can define $c_1 := 14 / 5$, $a_{2, 1} := 34$, and then we have $Y_1 = 0.98063... > 0$. Also, let $({c_2}', c_{0, 2}, a_{1, 2}, u_2) = (1, 1 / 7, 53 / 35, \sqrt{3})$, then we can define $c_2 := 7$, $a_{2, 2} := 84$, and then we have $Y_2 = 0.043547... > 0$.

\item[(2)]
Let $(t, b, s, k_0) = (8 / 25, 7, 3, 80)$, then $a_1 := 1 / 100$. Furthermore, let\\ $({c_1}', c_{0, 1}, a_{1, 1}, u_1) = (1, 11 / 60, 1 / 200, 16 \sqrt{3} / 25)$, then $c_1 := 60 / 11$, $a_{2, 1} := 1 / 5$, and then $Y_1 = 0.014518... > 0$. Also, let $({c_2}', c_{0, 2}, a_{1, 2}, u_2) = (1, 11 / 120, 1 / 200, 24 \sqrt{3} / 25)$, then $c_2 := 120 / 11$, $a_{2, 2} := 3 / 5$, and then $Y_2 = 0.29298... > 0$.

\item[(3)]
Let $(t, b, s, k_0) = (1 / 3, 13 / 20, 3, 22)$, then $a_1 := 1 / 100$. Furthermore, let\\ $({c_1}', c_{0, 1}, a_{1, 1}, u_1) = (1, 4 / 15, 1 / 200, 2 \sqrt{3} / 3)$, then $c_1 := 15 / 4$, $a_{2, 1} := 1 / 10$, and then $Y_1 = 0.80485... > 0$. Also, let $({c_2}', c_{0, 2}, a_{1, 2}, u_2) = (1, 2 / 15, 1 / 200, \sqrt{3})$, then $c_2 := 15 / 2$, $a_{2, 2} := 3 / 10$, and then $Y_2 = 0.96810... > 0$.\\

\item[(4)]
When $k = 8$, let $(t, b, s, k_0) = (1 / 2, 583 / 100, 2, 8)$, then $a_1 := 2363 / 500$. Furthermore, let\\ $({c_1}', c_{0, 1}, a_{1, 1}, u_1) = (1, 3 \sqrt{3} / 8, 177 / 50, \sqrt{3} / 4)$, then $c_1 := 8 / (3 \sqrt{3})$, $a_{2, 1} := 1063 / 100$, and then $X_1 - (8 / (3 \sqrt{3}))^{2 / k} \{ 1 + (2 \cdot 112^2 a_{2, 1}) / (8 k^3 / (3 \sqrt{3})) \} \geqslant 0.00012586... > 0$. Also, let\\ $({c_2}', c_{0, 2}, a_{1, 2}, u_2) = (1, \sqrt{3} / 8, 593 / 500, 3 \sqrt{3} / 4)$, then $c_2 := 8 / \sqrt{3}$, $a_{2, 2} := 321 / 10$, and then $X_2 - (8 / \sqrt{3})^{2 / k} \{ 1 + (2 \cdot (1 / 2)^2 a_{2, 1}) / (8 k^3 / \sqrt{3}) \} \geqslant 0.00031002... > 0$.

When $k \geqslant 10$, let $(t, b, s, k_0) = (1 / 2, 16 / 5, 2, 10)$, then $a_1 := 21 / 10$. Furthermore, let\\ $({c_1}', c_{0, 1}, a_{1, 1}, u_1) = (1, 3 \sqrt{3} / 8, 7 / 5, \sqrt{3} / 4)$, then $c_1 := 8 / (3 \sqrt{3})$, $a_{2, 1} := 4$, and then $Y_1 = 0.31693... > 0$. Also, we let\\ $({c_2}', c_{0, 2}, a_{1, 2}, u_2) = (1, \sqrt{3} / 8, 7 / 10, 3 \sqrt{3} / 4)$, then $c_2 := 8 / \sqrt{3}$, $a_{2, 2} := 17$, and then $Y_2 = 0.13820... > 0$.

\item[(5)]
Let $(t, b, s, k_0) = (2 / 3, 203 / 100, 2, 26)$, then $a_1 := 1 / 50$. Furthermore, let\\ $({c_1}', c_{0, 1}, a_{1, 1}, u_1) = (1, 5 / 12, 1 / 60, \sqrt{3} / 3)$, then $c_1 := 12 / 5$, $a_{2, 1} := 1 / 10$, and then $X_1 - (12 / 5)^{2 / k} \{ 1 + (2 \cdot (2 / 3)^2 a_{2, 1}) / (12 k^3 / 5) \} \geqslant 0.0032434... > 0$ for $k = 26$ and $Y_1 = 0.026412... > 0$ for $k \geqslant 44$. Also, let $({c_2}', c_{0, 2}, a_{1, 2}, u_2) = (1, 1 / 12, 1 / 300, \sqrt{3})$, then $c_2 := 12$, $a_{2, 2} := 1 / 2$, and then $X_2 - 12^{2 / k} \{ 1 + (2 \cdot (2 / 3)^2 a_{2, 2}) / (12 k^3) \} \geqslant 0.0081241... > 0$ for $k = 26$ and $Y_2 = 0.15714... > 0$ for $k \geqslant 44$.

\item[(6)]
Let $(t, b, s, k_0) = (1 / 2, 41 /20, 2, 70)$, then $a_1 := 1 / 100$. Furthermore, let\\ $({c_1}', c_{0, 1}, a_{1, 1}, u_1) = (1, 14 / 25, 1 / 200, \sqrt{3} / 4)$, then $c_1 := 25 / 14$, $a_{2, 1} := 1 / 50$, and then $Y_1 = 0.19093... > 0$. Also, let $({c_2}', c_{0, 2}, a_{1, 2}, u_2) = (1, 7 / 50, 1 / 200, 3 \sqrt{3} / 4)$, then $c_2 := 50 / 7$, $a_{2, 2} := 3 / 10$, and then $Y_2 = 0.031950... > 0$.

\item[(7)]
Let $(t, b, s, k_0) = (1 / 2, 41 / 20, 2, 200)$, then $a_1 := 1 / 100$. Furthermore, let\\ $({c_1}', c_{0, 1}, a_{1, 1}, u_1) = (1, 4883 / 9600, 1 / 200, \sqrt{3} / 4)$, then $c_1 := 9600 / 4883$, $a_{2, 1} := 1 / 50$, and then $Y_1 = 0.0037408... > 0$. Also, let $({c_2}', c_{0, 2}, a_{1, 2}, u_2) = (1, 257 / 1920, 1 / 200, 3 \sqrt{3} / 4)$, then $c_2 := 3 / 10$, $a_{2, 2} := 1920 / 257$, and then $Y_2 = 0.017771... > 0$.
\end{trivlist}
\begin{flushright} $\square$ \end{flushright}\quad

\subsection{Proofs of lemmas \ref{lem-f7s3}, \ref{lem-f7s4}, \ref{lem-f7s5}, and \ref{lem-f7s6}}\quad \vspace{-0.1in} \\

\noindent
{\it Proof of Lemma \ref{lem-f7s3}}\quad
We have $c_0 \leqslant \cos({c_0}' \pi) = - \cos((x / 180) \pi - (t / 2) \pi)$. We define $d_{1, 1} := 75 / 26$ and $d_{1, 2} := 100 / 51$ for $k \geqslant 80$ and $t < 1 / 2$. Furthermore, when $3217 \pi / 4500 \leqslant (x / 180) \pi < \alpha_{7, k} < (y / 180) \pi \leqslant 3 \pi / 4$, we have $\pi < 2 \pi / 3 + (x / 180) \pi + d_{1, 1} (t / 2) \pi < {\alpha_{7, k, 1}}' < 2 \pi / 3 + (y / 180) \pi + (3 t / 2) \pi < 2 \pi$ and $3 \pi / 2 < 4 \pi / 3 + (x / 180) \pi - t \pi < {\alpha_{7, k, 2}}' < 4 \pi / 3 + (y / 180) \pi - d_{1, 2} (t / 2) \pi < 2 \pi$ for the $t$ that appears in the Lemma. Thus, we can define ${c_1}'$ such that ${c_1}' \geqslant Max\{ 0, \cos(2 \pi / 3 + (y / 180) \pi + (3 t / 2) \pi) \}$ and ${c_2}'$ such that ${c_2}' \geqslant \cos(4 \pi / 3 + (y / 180) \pi - d_{1, 2} (t / 2) \pi)$.

For every item, we have $(b, s) = (7, 3)$, then we can define $a_1 := 1 / 100$.
\begin{trivlist}
\item[(1)]
Let $({c_1}', c_{0, 1}, a_{1, 1}, u_1) = (396677 / 500000, (23 / 36) (80639 / 250000), 1 / 200, \sqrt{3} / 2)$,\\
 then $c_1 := 7140486 / 1854697$, $a_{2, 1} := 1 / 10$, and then $Y_1 = 0.0072158... > 0$.

Also, let $({c_2}', c_{0, 2}, a_{1, 2}, u_2) = (873623 / 1000000, (13 / 36) (80639 / 250000), 1 / 200, 3 \sqrt{3} / 4)$,\\
 then $c_2 := 7862607 / 1048307$, $a_{2, 2} := 2 / 5$, and then $Y_2 = 0.012126... > 0$.

\item[(2)]
Let $({c_1}', c_{0, 1}, a_{1, 1}, u_1) = (152277 / 250000, (13 / 24) (365283 / 1000000), 1 / 200, 83 \sqrt{3} / 200)$,\\
 then $c_1 := 1624288 / 527631$, $a_{2, 1} := 1 / 10$, and then $Y_1 = 0.0058972... > 0$.

Also, let $({c_2}', c_{0, 2}, a_{1, 2}, u_2) = (905439 / 1000000, (11 / 24) (365283 / 1000000), 1 / 200, 249 \sqrt{3} / 400)$,\\
 then $c_2 := 2414504 / 446757$, $a_{2, 2} := 1 / 5$, and then $Y_2 = 0.0037228... > 0$.

\item[(3)]
Let $({c_1}', c_{0, 1}, a_{1, 1}, u_1) = (18377 / 40000, (53 / 120) (201373 / 500000), 1 / 200, 7 \sqrt{3} / 20)$,\\
 then $c_1 := 27565500 / 10672769$, $a_{2, 1} := 1 / 200$, and then $Y_1 = 0.0046808... > 0$.

Also, let $({c_2}', c_{0, 2}, a_{1, 2}, u_2) = (116867 / 125000, (67 / 120) (201373 / 500000), 1 / 200, 21 \sqrt{3} / 40)$,\\
 then $c_2 := 56096160 / 13491991$, $a_{2, 2} := 1 / 200$, and then $Y_2 = 0.0021192... > 0$.

\item[(4)]
Let $({c_1}', c_{0, 1}, a_{1, 1}, u_1) = (46421 / 125000, (449 / 1200) (423883 / 1000000), 1 / 200, 47 \sqrt{3} / 150)$,\\
 then $c_1 := 445641600 / 190323467$, $a_{2, 1} := 1 / 10$, and then $Y_1 = 0.0019413... > 0$.

Also, let $({c_2}', c_{0, 2}, a_{1, 2}, u_2) = (950261 / 1000000, (751 / 1200) (423883 / 1000000), 1 / 200, 47 \sqrt{3} / 100)$,\\
 then $c_2 := 1140313200 / 318336133$, $a_{2, 2} := 1 / 10$, and then $Y_2 = 0.0022918... > 0$.

\item[(5)]
Let $({c_1}', c_{0, 1}, a_{1, 1}, u_1) = (301039 / 1000000, (227 / 720) (441819 / 1000000), 1 / 200, 71 \sqrt{3} / 250)$,\\
 then $c_1 := 24083120 / 11143657$, $a_{2, 1} := 1 / 10$, and then $Y_1 = 0.0028859... > 0$.

Also, let $({c_2}', c_{0, 2}, a_{1, 2}, u_2) = (961843 / 1000000, (493 / 720) (441819 / 1000000), 1 / 200, 213 \sqrt{3} / 500)$,\\
 then $c_2 := 156080 / 49091$, $a_{2, 2} := 1 / 10$, and then $Y_2 = 0.0019675... > 0$.

\item[(6)]
Let $({c_1}', c_{0, 1}, a_{1, 1}, u_1) = (50093 / 200000, (971 / 3600) (454379 / 1000000), 1 / 200, 263 \sqrt{3} / 1000)$,\\
 then $c_1 := 901674000 / 441202009$, $a_{2, 1} := 1 / 10$, and then $Y_1 = 0.00056409... > 0$.

Also, let $({c_2}', c_{0, 2}, a_{1, 2}, u_2) = (969421 / 1000000, (2629 / 3600) (454379 / 1000000), 1 / 200, 789 \sqrt{3} / 2000)$,\\
 then $c_2 := 3489915600 / 1194562391$, $a_{2, 2} := 1 / 10$, and then $Y_2 = 0.00013254... > 0$..

\item[(7)]
Let $({c_1}', c_{0, 1}, a_{1, 1}, u_1) = (204497 / 1000000, (407 / 1800) (58279 / 125000), 1 / 200, 61 \sqrt{3} / 250)$,\\
 then $c_1 := 46011825 / 23719553$, $a_{2, 1} := 1 / 10$, and then $Y_1 = 0.0016244... > 0$.

Also, let $({c_2}', c_{0, 2}, a_{1, 2}, u_2) = (243913 / 250000, (1393 / 1800) (58279 / 125000), 1 / 200, 183 \sqrt{3} / 500)$,\\
 then $c_2 := 219521700 / 81182647$, $a_{2, 2} := 1 / 10$, and then $Y_2 = 0.000069553... > 0$.

\item[(8)]
Let $({c_1}', c_{0, 1}, a_{1, 1}, u_1) = (167767 / 1000000, (341 / 1800) (475839 / 1000000), 1 / 200, 143 \sqrt{3} / 625)$,\\
 then $c_1 := 33553400 / 18029011$, $a_{2, 1} := 1 / 5$, and then $Y_1 = 0.0019071... > 0$.

Also, let $({c_2}', c_{0, 2}, a_{1, 2}, u_2) = (980209 / 1000000, (1459 / 1800) (475839 / 1000000), 1 / 200, 429 \sqrt{3} / 1250)$,\\
 then $c_2 := 196041800 / 77138789$, $a_{2, 2} := 1 / 10$, and then $Y_2 = 0.00030133... > 0$.

\item[(9)]
Let $({c_1}', c_{0, 1}, a_{1, 1}, u_1) = (141593 / 1000000, (13 / 80) (482823 / 1000000), 1 / 200, 109 \sqrt{3} / 500)$,\\
 then $c_1 := 11327440 / 6276699$, $a_{2, 1} := 1 / 5$, and then $Y_1 = 0.0047590... > 0$.

Also, let $({c_2}', c_{0, 2}, a_{1, 2}, u_2) = (39327 / 40000, (67 / 80) (482823 / 1000000), 1 / 200, 327 \sqrt{3} / 1000)$,\\
 then $c_2 := 26218000 / 10783047$, $a_{2, 2} := 1 / 10$, and then $Y_2 = 0.00081355... > 0$. \qquad $\square$
\end{trivlist}\quad

\noindent
{\it Proof of Lemma \ref{lem-f7s4}}\quad
We have $\cos({c_0}' \pi) = \cos((y / 180) \pi - 2 \pi / 3 - (t / 2) \pi)$. We define $d_{2, 1} := 750 / 511$ and $d_{2, 2} := 250 / 533$ for $k \geqslant 62$ and $t < 1$. Furthermore, when $2 \pi / 3 \leqslant (x / 180) \pi < \alpha_{7, k} < (y / 180) \pi \leqslant 266 \pi / 375$, we have $0 \leqslant (x / 180) \pi - 2 \pi / 3 < \beta_{7, k} < (y / 180) \pi - 2 \pi / 3 \leqslant 16 \pi / 375$, $\pi / 2 < 4 \pi / 3 + (x / 180) \pi - 2 \pi / 3 - (3 t / 4) \pi < {\beta_{7, k, 1}}' < 4 \pi / 3 + (y / 180) \pi - 2 \pi / 3 - d_{2, 1} (t / 2) \pi < 3 \pi / 2$, and $\pi / 2 < 2 \pi / 3 + (x / 180) \pi - 2 \pi / 3 + d_{2, 2} (t / 2) \pi < {\beta_{7, k, 2}}' < 2 \pi / 3 + (y / 180) \pi - 2 \pi / 3 + (t / 4) \pi < \pi$. Thus, we can define ${c_2}'$ such that ${c_2}' \geqslant - \cos(2 \pi / 3 + (y / 180) \pi - 2 \pi / 3 + (t / 4) \pi)$.

For each item, let $(b, s) = (3, 2)$, then we can define $a_1 := 1 / 100$.
\begin{trivlist}
\item[(1)]
Let $({c_1}', c_{0, 1}, a_{1, 1}, u_1) = (1, (193 / 240) (330983 / 500000), 1 / 200, 93 \sqrt{3} / 400)$,\\
 then $c_1 := 120000000 / 63879719$, $a_{2, 1} := 1 / 50$, and then $Y_1 = 0.000019362... > 0$..

Also, let $({c_2}', c_{0, 2}, a_{1, 2}, u_2) = (422281 / 500000, (47 / 240) (330983 / 500000), 1 / 200, 279 \sqrt{3} / 400)$,\\
 then $c_2 := 101347440 / 15556201$, $a_{2, 2} := 3 / 10$, and then $Y_2 = 0.010264... > 0$.

For the following items, we can define ${c_1}'$ such that ${c_1}' \geqslant - \cos(4 \pi / 3 + (x / 180) \pi - 2 \pi / 3 - (3 t / 4) \pi)$.
\item[(2)]
Let $({c_1}', c_{0, 1}, a_{1, 1}, u_1) = (467413 / 500000, (269 / 360) (197271 / 250000), 1 / 200, 17 \sqrt{3} / 100)$,\\
 then $c_1 := 9348260 / 5896211$, $a_{2, 1} := 1 / 10$, and then $Y_1 = 0.0024253... > 0$.

Also, let $({c_2}', c_{0, 2}, a_{1, 2}, u_2) = (158883 / 200000, (91 / 360) (197271 / 250000), 1 / 200, 51 \sqrt{3} / 100)$,\\
 then $c_2 := 7944150 / 1994629$, $a_{2, 2} := 1 / 10$, and then $Y_2 = 0.0034755... > 0$.

\item[(3)]
Let $({c_1}', c_{0, 1}, a_{1, 1}, u_1) = (442697 / 500000, (717 / 1000) (103593 / 125000), 1 / 200, 11 \sqrt{3} / 75)$,\\
 then $c_1 := 110674250 / 74276181$, $a_{2, 1} := 1 / 50$, and then $Y_1 = 0.00015159... > 0$.

Also, let $({c_2}', c_{0, 2}, a_{1, 2}, u_2) = (193819 / 250000, (283 / 1000) (103593 / 125000), 1 / 200, 11 \sqrt{3} / 25)$,\\
 then $c_2 := 96909500 / 29316819$, $a_{2, 2} := 1 / 10$, and then $Y_2 = 0.00011691... > 0$. 

\item[(4)]
Let $({c_1}', c_{0, 1}, a_{1, 1}, u_1) = (843111 / 1000000, (83 / 120) (856447 / 1000000), 1 / 200, 13 \sqrt{3} / 100)$,\\
 then $c_1 := 101173320 / 71085101$, $a_{2, 1} := 1 / 50$, and then $Y_1 = 0.0012177... > 0$.
 
Also, let $({c_2}', c_{0, 2}, a_{1, 2}, u_2) = (151809 / 200000, (37 / 120) (856447 / 1000000), 1 / 200, 39 \sqrt{3} / 100)$,\\
 then $c_2 := 91085400 / 31688539$, $a_{2, 2} := 1 / 10$, and then $Y_2 = 0.0082259... > 0$. \qquad $\square$
\end{trivlist}\quad

\noindent
{\it Proof of Lemma \ref{lem-f7s5}}\quad
We have $\cos({c_0}' \pi) = - \cos((x / 180) \pi - (t / 2) \pi)$. Furthermore, when $2 \pi / 3 \leqslant (x / 180) \pi < \alpha_{7, k} < (y / 180) \pi \leqslant 3 \pi / 4$, we have $2 \pi < 4 \pi / 3 + (x / 180) \pi + d_{1, 1} (t / 2) \pi < {\alpha_{7, k, 1}}' < 4 \pi / 3 + (y / 180) \pi + (3 t / 2) \pi < 3 \pi$ and $\pi / 2 < 2 \pi / 3 + (x / 180) \pi - t \pi < {\alpha_{7, k, 2}}' < 2 \pi / 3 + (y / 180) \pi - d_{1, 2} (t / 2) \pi < 3 \pi / 2$. Thus, we can define ${c_2}' = 0$.

For every item, we have $(b, s) = (7, 3)$, then we can define $a_1 := 1 / 100$.
\begin{trivlist}
\item[(1)]
Let $({c_1}', c_{0, 1}, a_{1, 1}, u_1) = (1, 176 / 625, 1 / 100, 59 \sqrt{3} / 125)$, then $c_1 := 625 / 176$, $a_{2, 1} := 1 / 5$, and then $Y_1 = 0.0059892... > 0$.

We define $d_{1, 1} := 30 / 11$ for $k \geqslant 50$ and $t < 3 / 4$. Then, we can define ${c_1}'$ such that ${c_1}' \geqslant - \cos(4 \pi / 3 + (x / 180) \pi + d_{1, 1} (t / 2) \pi)$.
\item[(2)]
Let $({c_1}', c_{0, 1}, a_{1, 1}, u_1) = (12 / 25, 13 / 100, 1 / 100, \sqrt{3} / 2)$, then $c_1 := 48 / 13$, $a_{2, 1} := 3 / 10$, and then $Y_1 = 0.073156... > 0$.\qquad $\square$
\end{trivlist}\quad

\noindent
{\it Proof of Lemma \ref{lem-f7s6}}\quad
We have $\cos({c_0}' \pi) = \cos((y / 180) \pi - \pi / 3 + (t / 2) \pi)$. Furthermore, when $13 \pi / 36 \leqslant (x / 180) \pi < \alpha_{7, k} < (y / 180) \pi \leqslant 5 \pi / 9$, we have $\pi / 36 \leqslant (x / 180) \pi - 2 \pi / 3 < \beta_{7, k} < (y / 180) \pi - 2 \pi / 3 \leqslant 2 \pi / 9$, $0 < 2 \pi / 3 + (x / 180) \pi - 2 \pi / 3 - (3 t / 4) \pi < {\beta_{7, k, 1}}' < 2 \pi / 3 + (y / 180) \pi - 2 \pi / 3 - d_{2, 1} (t / 2) \pi < \pi$, and $\pi < 4 \pi / 3 + (x / 180) \pi - 2 \pi / 3 + d_{2, 2} (t / 2) \pi < {\beta_{7, k, 2}}' < 4 \pi / 3 + (y / 180) \pi - 2 \pi / 3 + (t / 4) \pi < 2 \pi$. Thus, we can define ${c_1}'$ such that ${c_1}' \geqslant - \cos(2 \pi / 3 + (y / 180) \pi - 2 \pi / 3 - d_{2, 1} (t / 2) \pi)$, and ${c_2}'$ such that ${c_2}' \geqslant Max\{0, - \cos(4 \pi / 3 + (x / 180) \pi - 2 \pi / 3 + d_{2, 2} (t / 2) \pi) \}$.

\begin{trivlist}
\item[(1)]
We define $d_{2, 1} := 30 / 23$ and $d_{2, 2} := 10 / 29$ for $k \geqslant 10$ and $t < 17 / 36$. Let $(t, b, s, k_0) = (2 / 5, 16 / 5, 2, 10)$, then $a_1 := 16 / 5$.

Let $({c_1}', c_{0, 1}, a_{1, 1}, u_1) = (2257 / 10000, (3 / 4) (203 / 500), (3 / 4) (16 / 5), \sqrt{3} / 5)$, then $c_1 := 2257 / 3045$, $a_{2, 1} := 67 / 10$, and then $Y_1 = 1.4014... > 0$.

Also, let $({c_2}', c_{0, 2}, a_{1, 2}, u_2) = (109 / 500, (1 / 4) (203 / 500), (1 / 4) (16 / 5), 3 \sqrt{3} / 5)$, then $c_2 := 436 / 203$, $a_{2, 2} := 20$, and then $Y_2 = 1.2571... > 0$.

\item[(2)]
We define $d_{2, 1} := 10 / 7$ and $d_{2, 2} := 10 / 23$ for $k \geqslant 28$ and $t < 2 / 3$. Let $(t, b, s, k_0) = (2 / 5, 41 / 20, 2, 28)$, then $a_1 := 1 / 25$. Now, we can define ${c_2}' = 0$.

Let $({c_1}', c_{0, 1}, a_{1, 1}, u_1) = (637 / 2000, 2419 / 10000, 1 / 25, \sqrt{3} / 5)$, then $c_1 := 3185 / 2419$, $a_{2, 1} := 3 / 10$, and then $Y_1 = 0.53163... > 0$.\qquad $\square$
\end{trivlist}\quad

\subsection{Proofs of lemmas \ref{lem-f7s7} and Lemma \ref{lem-f7s8}}\quad \vspace{-0.1in} \\

\noindent
{\it Proof of Lemma \ref{lem-f7s7}}\quad
We have $X_1 = v_k(2, 1, \theta_1)^{- 2 / k} \geqslant 1 + 2 \sqrt{3} t (\pi / k)$, $X_2 = v_k(3, 1, \theta_1)^{- 2 / k} \leqslant e^{- 3 \sqrt{3} (t / 2) \pi}$. We define $d_{1, 1} := 75 / 26$ and $d_{1, 2} := 100 / 51$ for $k \geqslant 46$ and $t < 1 / 3$. Furthermore, when $73 \pi / 120 \leqslant (x / 180) \pi < \alpha_{7, k} < (y / 180) \pi \leqslant 2 \pi / 3$, we have $2 \pi < 4 \pi / 3 + (x / 180) \pi + d_{1, 1} (t / 2) \pi < {\alpha_{7, k, 1}}' < 4 \pi / 3 + (y / 180) \pi + (3 t / 2) \pi < 5 \pi / 2$ and $\pi < 2 \pi / 3 + (x / 180) \pi - t \pi < {\alpha_{7, k, 2}}' < 2 \pi / 3 + (y / 180) \pi - d_{1, 2} (t / 2) \pi < 3 \pi / 2$. Then, we have $Sign\{\cos\left( k \theta_1 / 2 \right)\} = Sign\{Re(3 e^{- i \theta_1 / 2} + \sqrt{7} e^{i \theta_1 / 2})^{- k}\}$. Thus, we can define $c_0$ such that $c_0 \leqslant - \cos((x / 180) \pi - (t / 2) \pi) + {c_2}'' e^{- 3 \sqrt{3} (t / 2) \pi}$ and ${c_1}'$ such that ${c_1}' \geqslant \cos(4 \pi / 3 + (x / 180) \pi + d_{1, 1} (t / 2) \pi)$, where ${c_2}'' \leqslant - \cos(2 \pi / 3 + (y / 180) \pi - d_{1, 2} (t / 2) \pi)$.

For every item, we have $(b, s) = (7, 3)$, then we can define $a_1 := 1 / 100$.
\begin{trivlist}
\item[(1)]
Let $({c_1}', c_{0, 1}, a_{1, 1}, u_1) = (854381 / 1000000, 23487 / 62500, 1 / 100, 23 \sqrt{3} / 75)$,\\
 then $c_1 := 854381 / 375792$, $a_{2, 1} := 1 / 10$, and then $Y_1 = 0.0092228... > 0$.

\item[(2)]
Let $({c_1}', c_{0, 1}, a_{1, 1}, u_1) = (480483 / 500000, 139811 / 250000, 1 / 100, \sqrt{3} / 5)$,\\
 then $c_1 := 480483 / 279622$, $a_{2, 1} := 1 / 10$, and then $Y_1 = 0.0047822... > 0$.

\item[(3)]
Let $({c_1}', c_{0, 1}, a_{1, 1}, u_1) = (991419 / 1000000, 170601 / 250000, 1 / 100, 86 \sqrt{3} / 625)$,\\
 then $c_1 := 330473 / 227468$, $a_{2, 1} := 1 / 10$, and then $Y_1 = 0.0014263... > 0$.

\item[(4)]
Let $({c_1}', c_{0, 1}, a_{1, 1}, u_1) = (998509 / 1000000, 751257 / 1000000, 1 / 100, 21 \sqrt{3} / 200)$,\\
 then $c_1 := 998509 / 751257$, $a_{2, 1} := 1 / 10$, and then $Y_1 = 0.0021540... > 0$.
\end{trivlist}\quad

\noindent
{\it Proof of Lemma \ref{lem-f7s8}}\quad
We have $X_1 = (1 / 4) \; v_k(1, -1, \theta_2)^{- 2 / k} \geqslant 1 + t (\pi / k)$, $X_2 = (1 / 4) \; v_k(3, -1, \theta_2)^{- 2 / k} \leqslant e^{- (3 \sqrt{3} / 2) (t / 2) \pi}$. We define $d_{2, 1} := 25 / 17$ and $d_{2, 2} := 25 / 53$ for $k \geqslant 46$ and $t < 5 / 18$. Furthermore, when $5 \pi / 9 \leqslant (x / 180) \pi < \alpha_{7, k} < (y / 180) \pi \leqslant 217 \pi / 360$, we have $2 \pi / 9 \leqslant (x / 180) \pi - \pi / 3 < \beta_{7, k} < (y / 180) \pi - \pi / 3 \leqslant 97 \pi / 360$, $\pi / 2 < 2 \pi / 3 + (x / 180) \pi - \pi / 3 - (3 t / 4) \pi < {\beta_{7, k, 1}}' < 2 \pi / 3 + (y / 180) \pi - \pi / 3 - d_{2, 1} (t / 2) \pi < \pi$, and $3 \pi / 2 < 4 \pi / 3 + (x / 180) \pi - \pi / 3 + d_{2, 2} (t / 2) \pi < {\beta_{7, k, 2}}' < 4 \pi / 3 + (y / 180) \pi - \pi / 3 + (t / 4) \pi < 2 \pi$ for the $t$ appearing in the Lemma. Thus, we can define $c_0$ such that $c_0 \leqslant \cos((y / 180) \pi - \pi / 3 + (t / 2) \pi) + {c_2}'' e^{- (3 \sqrt{3} / 2) (t / 2) \pi}$ and ${c_1}'$ such that ${c_1}' \geqslant \cos(2 \pi / 3 + (y / 180) \pi - \pi / 3 - d_{2, 1} (t / 2) \pi)$, where ${c_2}'' \geqslant \cos(4 \pi / 3 + (x / 180) \pi - \pi / 3 + d_{2, 2} (t / 2) \pi)$.

For every item, we have $(b, s) = (3, 2)$, then we can define $a_1 := 1 / 100$.
\begin{trivlist}
\item[(1)]
Let $({c_1}', c_{0, 1}, a_{1, 1}, u_1) = (591931 / 1000000, 405981 / 1000000, 1 / 100, 3 \sqrt{3} / 20)$,\\
 then $c_1 := 591931 / 405981$, $a_{2, 1} := 1 / 10$, and then $Y_1 = 0.055699... > 0$.

\item[(2)]
Let $({c_1}', c_{0, 1}, a_{1, 1}, u_1) = (187477 / 250000, 111413 / 200000, 1 / 100, 11 \sqrt{3} / 100)$,\\
 then $c_1 := 749908 / 557065$, $a_{2, 1} := 1 / 10$, and then $Y_1 = 0.0031102... > 0$.

\item[(3)]
Let $({c_1}', c_{0, 1}, a_{1, 1}, u_1) = (832649 / 1000000, 666201 / 1000000, 1 / 100, 33 \sqrt{3} / 400)$,\\
 then $c_1 := 832649 / 666201$, $a_{2, 1} := 1 / 10$, and then $Y_1 = 0.0027742... > 0$.

\item[(4)]
Let $({c_1}', c_{0, 1}, a_{1, 1}, u_1) = (218173 / 250000, 728483 / 1000000, 1 / 100, \sqrt{3} / 15)$,\\
 then $c_1 := 872692 / 728483$, $a_{2, 1} := 1 / 10$, and then $Y_1 = 0.0014577... > 0$.

\item[(5)]
Let $({c_1}', c_{0, 1}, a_{1, 1}, u_1) = (447639 / 500000, 38427 / 50000, 1 / 100, 113 \sqrt{3} / 2000)$,\\
 then $c_1 := 149213 / 128090$, $a_{2, 1} := 1 / 10$, and then $Y_1 = 0.0021090... > 0$.
\end{trivlist}\quad

\quad\\

\begin{center}
{\large Acknowledgement.}
\end{center}
I would like to thank Professor Eiichi Bannai for suggesting these problems as a doctoral course project.

\quad\\

\end{document}